\documentclass{article}

\usepackage{amssymb,latexsym,amsmath}
\usepackage[pdftex]{graphicx}

\usepackage{hyperref}

\begin{document}

\pagestyle{plain}

\newtheorem{theorem}{Theorem}[section]

\newtheorem{proposition}[theorem]{Proposition}

\newtheorem{lema}[theorem]{Lemma}

\newtheorem{corollary}[theorem]{Corollary}

\newtheorem{definition}[theorem]{Definition}

\newtheorem{remark}[theorem]{Remark}

\newtheorem{exempl}{Example}[section]

\newenvironment{exemplu}{\begin{exempl}  \em}{\hfill $\square$

\end{exempl}}

\newcommand{\ea}{\mbox{{\bf a}}}

\newcommand{\eu}{\mbox{{\bf u}}}

\newcommand{\ueu}{\underline{\eu}}

\newcommand{\ueo}{\overline{u}}

\newcommand{\oeu}{\overline{\eu}}

\newcommand{\ew}{\mbox{{\bf w}}}

\newcommand{\ef}{\mbox{{\bf f}}}

\newcommand{\eF}{\mbox{{\bf F}}}

\newcommand{\eC}{\mbox{{\bf C}}}

\newcommand{\en}{\mbox{{\bf n}}}

\newcommand{\eT}{\mbox{{\bf T}}}

\newcommand{\eL}{\mbox{{\bf L}}}

\newcommand{\eR}{\mbox{{\bf R}}}

\newcommand{\eV}{\mbox{{\bf V}}}

\newcommand{\eU}{\mbox{{\bf U}}}

\newcommand{\ev}{\mbox{{\bf v}}}

\newcommand{\eve}{\mbox{{\bf e}}}

\newcommand{\uev}{\underline{\ev}}

\newcommand{\eY}{\mbox{{\bf Y}}}

\newcommand{\eK}{\mbox{{\bf K}}}

\newcommand{\eP}{\mbox{{\bf P}}}

\newcommand{\eS}{\mbox{{\bf S}}}

\newcommand{\eJ}{\mbox{{\bf J}}}

\newcommand{\eB}{\mbox{{\bf B}}}

\newcommand{\eH}{\mbox{{\bf H}}}

\newcommand{\leb}{\mathcal{ L}^{n}}

\newcommand{\eI}{\mathcal{ I}}

\newcommand{\eE}{\mathcal{ E}}

\newcommand{\hen}{\mathcal{H}^{n-1}}

\newcommand{\eBV}{\mbox{{\bf BV}}}

\newcommand{\eA}{\mbox{{\bf A}}}

\newcommand{\eSBV}{\mbox{{\bf SBV}}}

\newcommand{\eBD}{\mbox{{\bf BD}}}

\newcommand{\eSBD}{\mbox{{\bf SBD}}}

\newcommand{\ecs}{\mbox{{\bf X}}}

\newcommand{\eg}{\mbox{{\bf g}}}

\newcommand{\paromega}{\partial \Omega}

\newcommand{\gau}{\Gamma_{u}}

\newcommand{\gaf}{\Gamma_{f}}

\newcommand{\sig}{{\bf \sigma}}

\newcommand{\gac}{\Gamma_{\mbox{{\bf c}}}}

\newcommand{\deu}{\dot{\eu}}

\newcommand{\dueu}{\underline{\deu}}

\newcommand{\dev}{\dot{\ev}}

\newcommand{\duev}{\underline{\dev}}

\newcommand{\weak}{\stackrel{w}{\approx}}

\newcommand{\mild}{\stackrel{m}{\approx}}

\newcommand{\lrightarrow}{\stackrel{L}{\rightarrow}}

\newcommand{\rrightarrow}{\stackrel{R}{\rightarrow}}

\newcommand{\strong}{\stackrel{s}{\approx}}

\newcommand{\weakdown}{\rightharpoondown}

\newcommand{\opg}{\stackrel{\mathfrak{g}}{\cdot}}

\newcommand{\opunu}{\stackrel{1}{\cdot}}
\newcommand{\opdoi}{\stackrel{2}{\cdot}}

\newcommand{\opn}{\stackrel{\mathfrak{n}}{\cdot}}
\newcommand{\opx}{\stackrel{x}{\cdot}}

\newcommand{\tr}{\ \mbox{tr}}

\newcommand{\Ad}{\ \mbox{Ad}}

\newcommand{\ad}{\ \mbox{ad}}

\renewcommand{\contentsname}{ }

\title{Computing with space: a tangle formalism for chora and difference}

\author{Marius Buliga \\ 
\\
Institute of Mathematics, Romanian Academy \\
P.O. BOX 1-764, RO 014700\\
Bucure\c sti, Romania\\
{\footnotesize Marius.Buliga@imar.ro}}

\date{This version:  18.04.2011}

\maketitle

\begin{abstract}
What is space computing, simulation, or understanding? Converging from several 
sources, this seems to be something more primitive than what is usually meant 
by computation, something that was along with us since antiquity 
(the word "choros", "chora", denotes "space" or "place" and is seemingly the 
most mysterious notion from Plato, described in Timaeus 48e - 53c) which has 
to do with cybernetics and with the understanding of the front end visual 
system. It may have some unexpected applications, also. Here, inspired by 
Bateson (see Supplementary Material), I explore from the mathematical side 
the point of view that there is no difference between the map and the 
territory, but instead the transformation of one into another can be 
understood by using a formalism of tangle diagrams. 

This paper continues arXiv:1009.5028 "What is a space? Computations in emergent 
algebras and the front end visual system" and the 
arXiv:1007.2362 "Introduction to metric spaces with 
dilations". 
\end{abstract}

\newpage

%\tableofcontents

%\newpage

%\thispagestyle{empty}

%This is a paper based on the following: 
%\begin{enumerate}
%\item[-] \url{http://arxiv.org/abs/1103.6007}
%\item[-] \url{http://arxiv.org/abs/1009.5028}
%\item[-] \url{http://arxiv.org/abs/1011.4485}
%\end{enumerate}

%\vfill

%\centerline{\includegraphics[angle=270, width=0.7\textwidth]{r3_4.pdf}}

%\vfill

%\newpage

\tableofcontents

\newpage

\section{The map is the territory}

There is a whole 
discussion around the key phrases "The map is not the territory" and "The map 
is the territory". The map-territory relation has a wikipedia 
entry \cite{wikimap}  
 where we find the following relevant citation. 

\vspace{.5cm}

"The expression "the map is not the territory" first appeared in print in the  
paper \cite{korzybski} that Alfred Korzybski gave at a meeting of the American 
Association for the Advancement of Science in New Orleans, Louisiana in 1931:
\begin{enumerate}
\item[A)] A map may have a structure similar or dissimilar to the structure of 
the territory, 
\item[B)] A map is not the territory.
\end{enumerate}

Korzybski's dictum "the map is not the territory" is [...] used 
to signify that individual people in fact do not in general have 
access to absolute knowledge of reality, but in fact only have access 
to a set of beliefs they have built up over time, about reality. 
So it is considered important to be aware that people's beliefs about 
reality and their awareness of things (the "map") are not reality itself or 
everything they could be aware of ("the
territory")."

\vspace{.5cm}

Here, inspired by Bateson (see Supplementary Material), I explore from 
the mathematical side the point 
of view that there is no difference between the map and the territory, 
but instead the transformation of one into another can be understood 
by using a tangle diagram, of one of the types figured below.

\vspace{.5cm}

\centerline{\includegraphics[angle=270, width=0.7\textwidth]{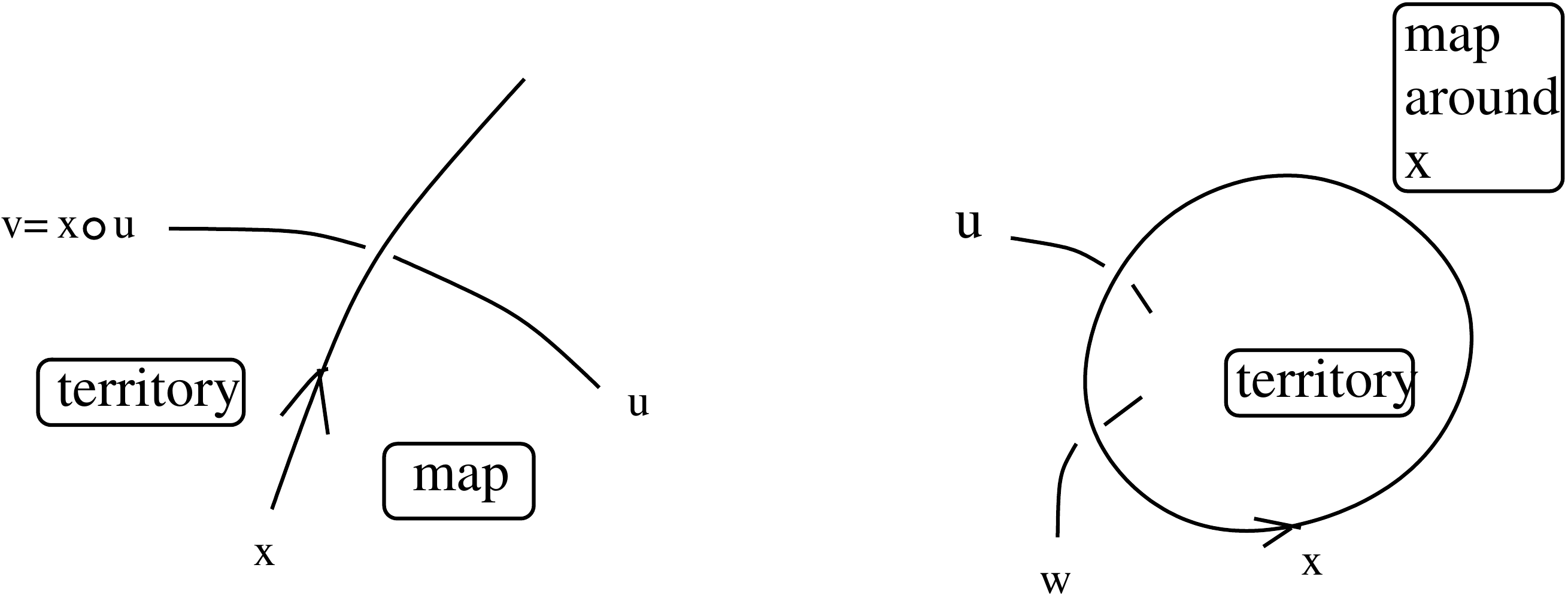}}
%\caption{}
%\label{elemchoracircuit2}

\vspace{.5cm}

The figure at the right, where the arc decorated by "$x$" closes itself, plays
 a special role. For such diagrams I shall use the name "place", or "chora",
 inspired by the discussion about space, or place, or chora, from Plato' Timaeus
(see Supplementary Material about this dialogue).

More precisely, we may imagine that the exploration of the territory provides 
us with an atlas, a collection of maps, mathematically understood as 
 a family of two operations (see later "emergent algebra"). A more precise
 figure than the previous one from the right would be the following.

\vspace{.5cm}

\centerline{\includegraphics[angle=270, width=0.8\textwidth]{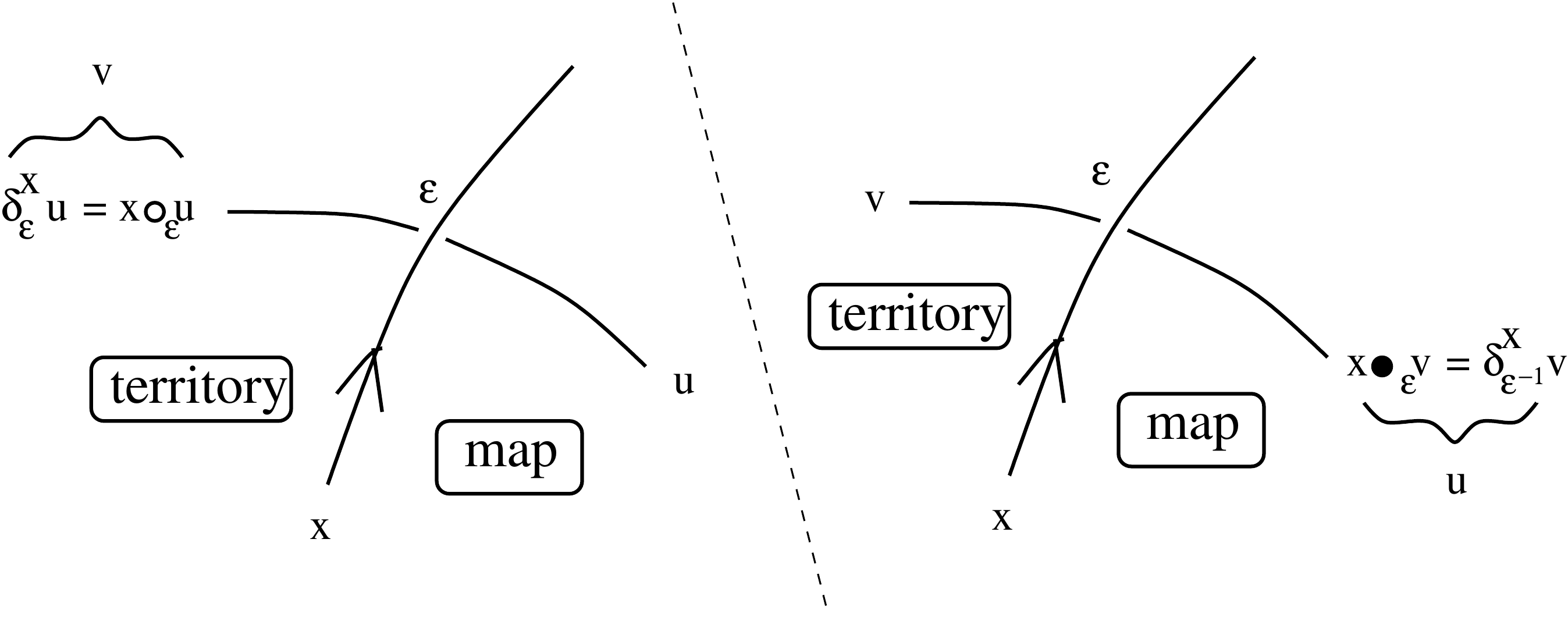}}
%\caption{}
%\label{elemchoracircuit2}

\vspace{.5cm}

The "$\varepsilon$" which decorates the crossing represents the scale of the map. 
In the figure at the left, the point "$v$" from the territory 
is represented by the "pixel" $u$  from the "map space".

 By accepting that "the map is the territory", the point 
$x$ (as well as the other point $v$ and the pixel $u$) is at the same time on
the  map and in the territory. The decorated crossing from the left of the
previous figure represents a map, at scale $\varepsilon$, of the territory 
(near $x$), which is a function, denoted by $\displaystyle
\delta^{x}_{\varepsilon}$ which takes $u$ to $v$. Equivalently, we may see such
a map as a binary operation, denoted by $\displaystyle \circ_{\varepsilon}$, 
with inputs $x$ and $u$ and output $v$. We describe all this by the algebraic 
relations: 
$$v\, = \, x \circ_{\varepsilon} u \, = \, \delta^{x}_{\varepsilon} u$$ 

The inverse transformation, or operation, of passing from the territory to the
map, is described by the decorated crossing from the right of the previous
figure. Algebraically, we write
$$u\, = \, x \bullet_{\varepsilon} v \, = \, \delta^{x}_{\varepsilon^{-1}} v$$ 

Imagine now a complex diagram, with lots of crossings, decorated by various
scale parameters, and segments decorated with points from a space $X$ which 
is seen both as territory (to explore) and map (of it). In such a diagram 
the convention map-territory can be only local, around each crossing. 

There is
though a diagram which could unambiguously  serve as a symbol for 
"the place (near) $x$ at scale $\varepsilon$".

%\vspace{.5cm}

\centerline{\includegraphics[angle=270, width=0.45\textwidth]{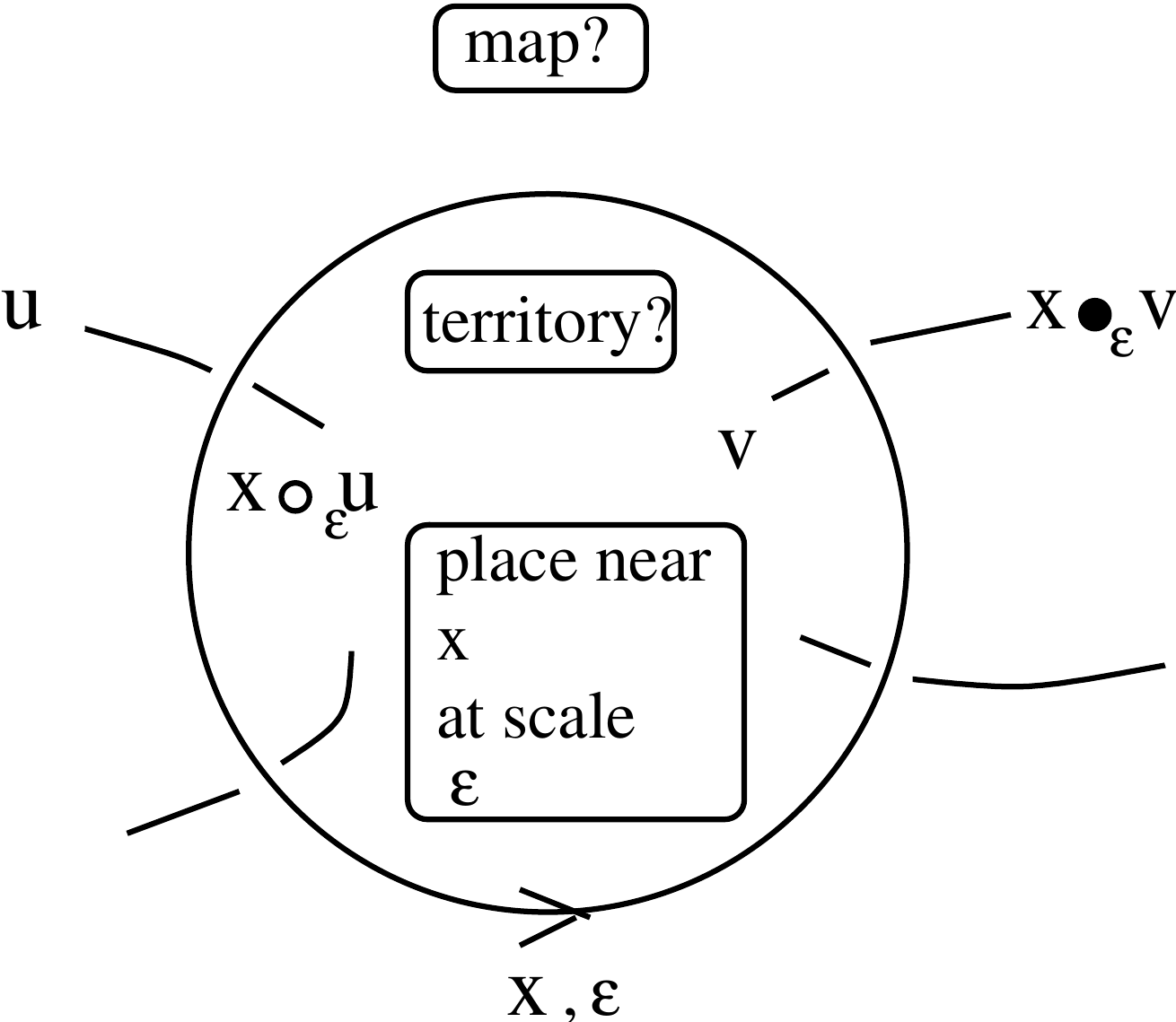}}
%\caption{}
%\label{elemchoracircuit2}

\vspace{.5cm}

In this diagram, all crossings which are not decorated have   
$\varepsilon$ as a decoration, but this decoration can be unambiguously 
placed near the decoration "$x$" of the closed arc. Such a diagram will bear the name "place (or chora) $x$ at scale $\varepsilon$".

There is another important diagram, called the "difference".

\vspace{.5cm}

\centerline{\includegraphics[angle=270, width=0.45\textwidth]{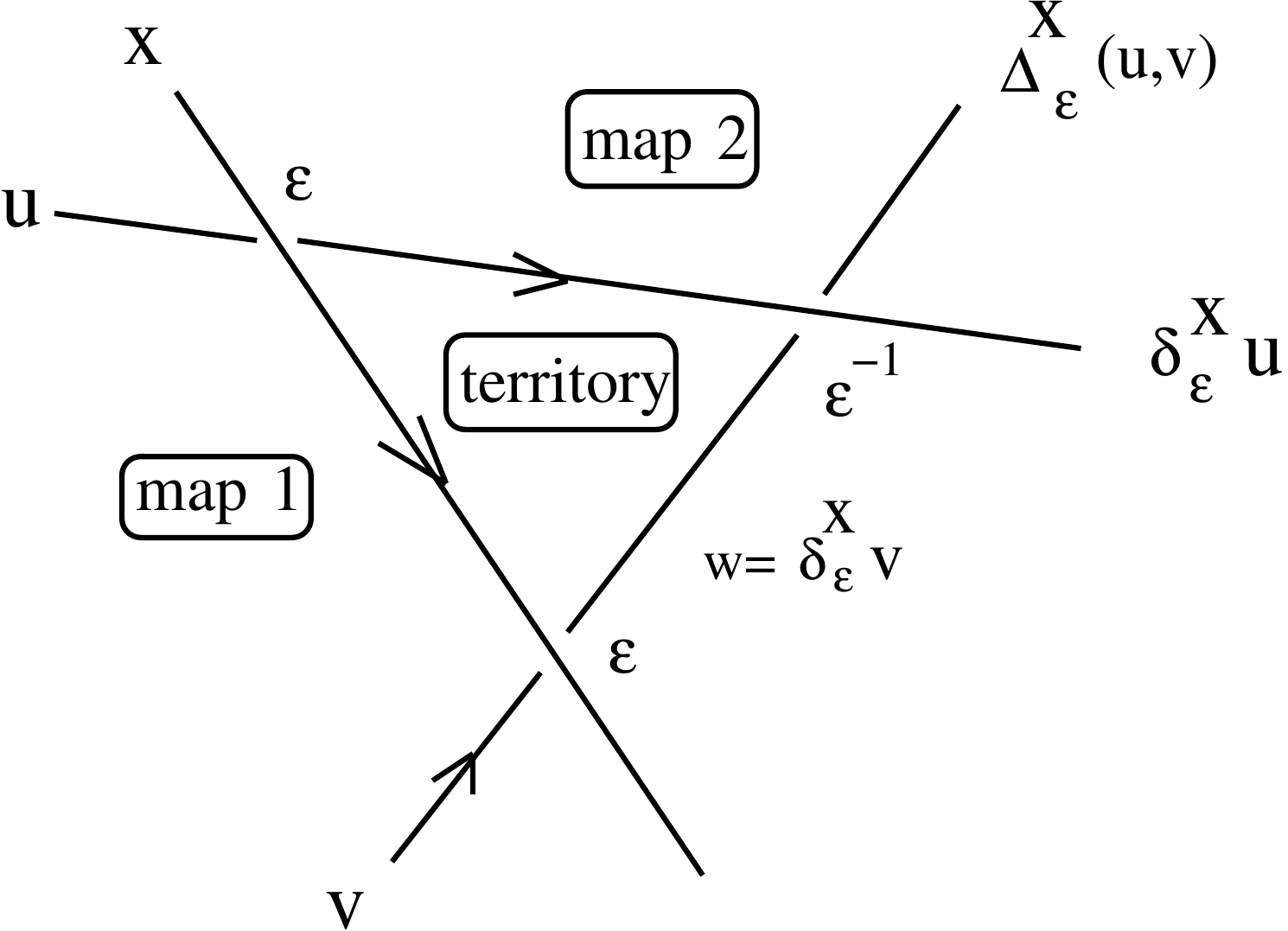}}
%\caption{}
%\label{elemchoracircuit2}

\vspace{.5cm}

Indeed, a point "$w$" from a territory, is 
seen in two maps. The "map 1" represents the map of the territory made around
$x$. In this map the point $w$ is represented by the pixel "$v$", therefore 
$$w = \delta_{\varepsilon}^{x} v$$
In the same map 1 choose another pixel, denoted by "$u$". In the territory,
this corresponds to the point $\displaystyle \delta_{\varepsilon}^{x} u$. But
this point admits a map made around it, this is "map 2". With this map, the
same point $w$ from the territory is seen at the pixel $\displaystyle 
\Delta_{\varepsilon}^{x}(u,v)$. 

There is a  function which takes  the pixel $v$, the image of $w$ in the map 1, 
to the pixel  $\displaystyle \Delta_{\varepsilon}^{x}(u,v)$, the image of the
same $w$ in the map 2. This transition function represents the 
{\it  difference} between the two maps. In algebraic terms, this difference 
has the expression: 
$$\Delta_{\varepsilon}^{x}(u,v) \, = \, \left( x\circ_{\varepsilon}u\right)  
\bullet_{\varepsilon} \left(x\circ_{\varepsilon}v\right) \, = \,
\delta_{\varepsilon^{-1}}^{\delta^{x}_{\varepsilon} u} \delta^{x}_{\varepsilon}
v$$
The function $\displaystyle (u,v) \mapsto \Delta_{\varepsilon}^{x}(u,v)$ is 
called the difference at scale $\varepsilon$, with respect to $x$.

\section{Computing with space}

  What is space computing,simulation, or understanding? Converging from 
  several sources, this seems to be something more primitive than what is 
  meant nowadays by
  "computation"\footnote{\url{http://www.changizi.com/viscomp.pdf}}, 
  something that was along with us since 
  antiquity\footnote{The word "choros", "chora",  denotes "space" or "place" and is  seemingly the most mysterious notion 
from Plato, described in Timaeus 48e - 53c, and appears as  a third class needed in Plato's description 
of the
universe, \url{http://www.ellopos.net/elpenor/physis/plato-timaeus/space.asp}},  which has to do with cybernetics and with the understanding of 
  the visual system. 
  
  It may have some unexpected applications, also. 
  
  My own interest into the subject of simulating space emerged from the 
  realization that most of differential calculus (the background of 
  physics) can be understood as a kind of construction and manipulation 
  of assemblies of dilations as elementary gates (transistor-like). Then, 
  I ask myself, what is the computer which may be constructed with such 
  transistors? For the stem of the idea of space simulation, see the page: 
  
  \url{http://imar.ro/~mbuliga/buliga_sim.html}

\subsection{Computing with tangles}

  Computations based on manipulation of tangle diagrams are in fact explored 
   ground, for example in Kauffman and Lomonaco Topological Quantum
Computing\footnote{\url{http://www.math.uic.edu/~kauffman/Quanta.pdf}} 
which is basically computation with braids, which may be implemented (or not) 
by anyons\footnote{\url{http://en.wikipedia.org/wiki/Topological_quantum_computer}}.   
This type of braid computation branded as quantum could be relevant for 
 "space simulation" via the braids formalism first explained in  
 \cite{buligafrontend}. See also the recent paper by Meredith and Snyder
 \cite{meredith} "Knots as processes: a new kind of invariant".  
 Meredith is working towards a kind of space computation, see this 
link\footnote{\url{http://biosimilarity.blogspot.com/2009/12/secret-life-of-space.html}}.

\subsection{Computations in the front-end visual system  as a paradigm}

In mathematics "spaces" come in many flavours. There are vector spaces, affine spaces, symmetric spaces, groups and so on. We usually take such objects as the stage where the plot of reasoning is laid. 
But in fact what we use, in many instances, are properties of particular spaces which, I claim, can be seen as coming from a particular class of computations.  

There is though a "space" which is "given" almost beyond doubt, namely the physical space where we all live. But as it regards perception of this space, we know now that things are not so simple. 
As I am writing these notes, here in Baixo Gavea, my eyes are attracted by the 
wonderful complexity 
of a tree near my window. The nature of the tree is foreign to me, as are the 
other smaller beings growing on or around the tree. I can make some 
educated guesses about what they are: some are orchids, there is a smaller, 
iterated version of the big tree. However, somewhere in my brain,  at a very 
fundamental level, the visible space is constructed in my head, before the 
stage where I am capable of recognizing and naming the objects or beings that 
I see. I cite from Koenderink  \cite{koen}, p. 126: 
\vspace{.5cm}

"The brain can organize {\em itself} through information obtained via interactions with the physical world into 
an embodiment of geometry, it becomes a veritable {\em geometry engine}. [...] 

There may be a point in holding that many of the better-known brain processes 
are most easily understood in terms of differential geometrical calculations 
running on massively parallel processor arrays whose nodes can be understood 
quite directly in terms of multilinear operators (vectors, tensors, etc). 
In this view brain processes in fact are space."

\vspace{.5cm}
 
In the paper \cite{koen2} Koenderink, Kappers and van Doorn study the "front end visual system" starting from general invariance  principles.  They define the "front end visual system" as "the interface between the light field and those parts of the brain nearest to the transduction stage".  After explaining that 
the "exact limits of the interface are essentially arbitrary", the authors propose the following 
characterization of what the front end visual system does, section 1 \cite{koen2}: 
\begin{enumerate}
\item[1.] "the front end is a "machine" in the sense of a syntactical transformer; 
\item[2.] there is no semantics. The front end processes structure;
\item[3.] the front end is precategorical, thus -- in a way -- the front end does not compute anything; 
\item[4.] the front end operates in a bottom-up fashion. Top down commands based upon semantical interpretations are not considered to be part of the front end proper;
\item[5.] the front end is a deterministic machine; ... all output depends causally on the (total) input from the immediate past." 
\item[6.] "What is not explicitly encoded by the front end is irretrievably lost. Thus the front end should be 
universal (undedicated) and yet should provide explicit data structures (in order to sustain fast processing past the front end) without sacrificing completeness (everything of potential importance to the survival of the agent has to be represented somehow)."
\end{enumerate}
 
 The authors arrive at the conclusion that the (part of the brain dealing with 
vision) is a "geometry engine", working somehow with the help of many elementary
circuits which implement (a discretization of?) differential calculus using
partial derivatives up to order 4. 
 I cite from the last paragraph of section 1, just 
 before the beginning of sections 1.1 \cite{koen2}.

\vspace{.5cm}

"In a local representation one can do without extensive (that is spatial, or 
geometrical) properties and represent everything in terms of intensive 
properties. This obviates the need for explicit geometrical expertise. 
The local representation of geometry is the typical tool of differential 
geometry. ... The columnar organization of representation in primate visual 
cortex suggests exactly such a structure."

\vspace{.5cm}

So, the front end does perform a kind of a computation, although a  very 
particular one. It is not, at first sight, a logical, boolean type of 
computation. I think this is what Koenderink and coauthors want to say in  
point 3. of  the characterization of the front end visual system.

 We may imagine that there is an abstract mathematical "front end" which, if 
 fed with the definition of  a "space",  then spews out a "data structure" 
 which is used for "past processing", that is for mathematical reasoning in 
 that space. (In fact, when we say "let $M$ be a manifold", for example, we don't "have" that manifold,  only some properties of it, together with some very general abstract nonsense concerning 
 "legal" manipulations in the universe of "manifolds". All these can be liken with the image that we get 
 past the "front end" , in the sense that, like a real perceived image, we see it all, but we are incapable 
 of really enumerating and naming all that we see.)
 
 Even more, we may think that  the physical space can be understood, at some very fundamental level,   
 as the input of a "universal front end", and physical observers are "universal front ends". That is, maybe biology uses at a different scale an embodiment of a fundamental mechanism of the nature.

Thus, the  biologically inspired viewpoint is  that observers are like universal front ends looking at 
the same (but otherwise unknown) space. Interestingly this  may give a link between the problem of 
"local sign" in neuroscience and  the problem of understanding the properties of the physical space as emerging from some non-geometrical, more  fundamental structure, like a net, a foam, a graph... 

Indeed, in this physics research, one wants to obtain geometrical structure of the space (for example 
that locally, at the macroscopic scale,  it looks like $\mathbb{R}^n$) from a  non spatial like  structure, 
"seen from afar" (not unlike Gromov does with metric spaces). But in fact the brain does this all the time: 
from a class of intensive quantities (like the electric impulses sent by the neurons in the retina) the front 
end visual system reconstructs the space, literally as we see it. How it does it without "geometrical 
expertise" is called in neuroscience the problem of the "local sign" or of the "homunculus".

In \cite{buligafrontend}  is proposed an equivalence between the problem
 of "local sign" in vision and the problem of constructing
 space as an emergent reality, coming from a non-spatial
 substrate.  In \cite{buligamore} I formulated things more clearly and proposed
  the following equivalence:
\begin{enumerate}
\item[(A)] reality emerges from a more primitive, non-geometrical, reality
\item[ ] in the same way as 
\item[(B)] the brain construct (understands, simulates, transforms, encodes or decodes) the image 
of reality, starting from intensive properties (like  a bunch of spiking signals sent by receptors in the retina), without any use of extensive (i.e. spatial or geometric)  properties. 
\end{enumerate}

The problem (B) is known in life sciences as the problem of "local sign" \cite{koen} \cite{koen2}. Indeed, any use of extensive 
properties would lead to the "homunculus
fallacy"\footnote{\url{http://en.wikipedia.org/wiki/Homunculus_argument}}.

These equivalent problems are difficult and wonderfully simple: 
\begin{enumerate}
\item[-] we don't know how to solve completely problem (A) (in physics) or problem (B) (in neuroscience),  
\item[-] but our ventral/dorsal streams and cerebellum do this all the time in about 150 ms. Moreover, 
any fly does it, as witnessed by its spectacular flight capacities  \cite{fly}. 
\end{enumerate}

The backbone of the argument is  that, as boolean logic is based on a primitive 
gate (like NAND), differential calculus and differential geometry are also based on a primitive gate (a dilation gate), appearing naturally from the least 
sophisticated strategy of exploration which a binocular creature might have,  namely 
jumping randomly from a place to another, orienting herself by comparing what she 
sees with her two eyes. 

\subsection{Exploring space}
I shall discuss about the simplest strategy to explore an unknown territory  
$X$.   For this we send an explorer to look around.  The explorer will make 
charts of parts of $X$ {\em into} $X$, then use these charts in order to plan its 
movements. 

I shall suppose that we can put a distance on the set $X$, that is a function 
$d: X \times X \rightarrow [0,+\infty)$ which satisfies the following 
requirement:  for any three points $x,y,z \in X$ there is a bijective 
correspondence with a triple A,B,C in the plane  such that the sizes (lengths) 
of AB, BC, AC are equal respectively with $d(x,y)$, $d(y,z)$, $d(z,x)$. 
Basically, we accept that we can represent in the plane any three points from the 
space $X$. An interpretation of the distance $d(x,y)$ is the following:
the explorer has a ruler and $d(x,y)$ is the numerical value shown by the ruler 
when streched between points marked with "$x$" and "$y$". (Then the explorer has
to use somehow these numbers in order to make a chart of $X$.)

\vspace{.5cm}

\centerline{\includegraphics[angle=270, width=0.5\textwidth]{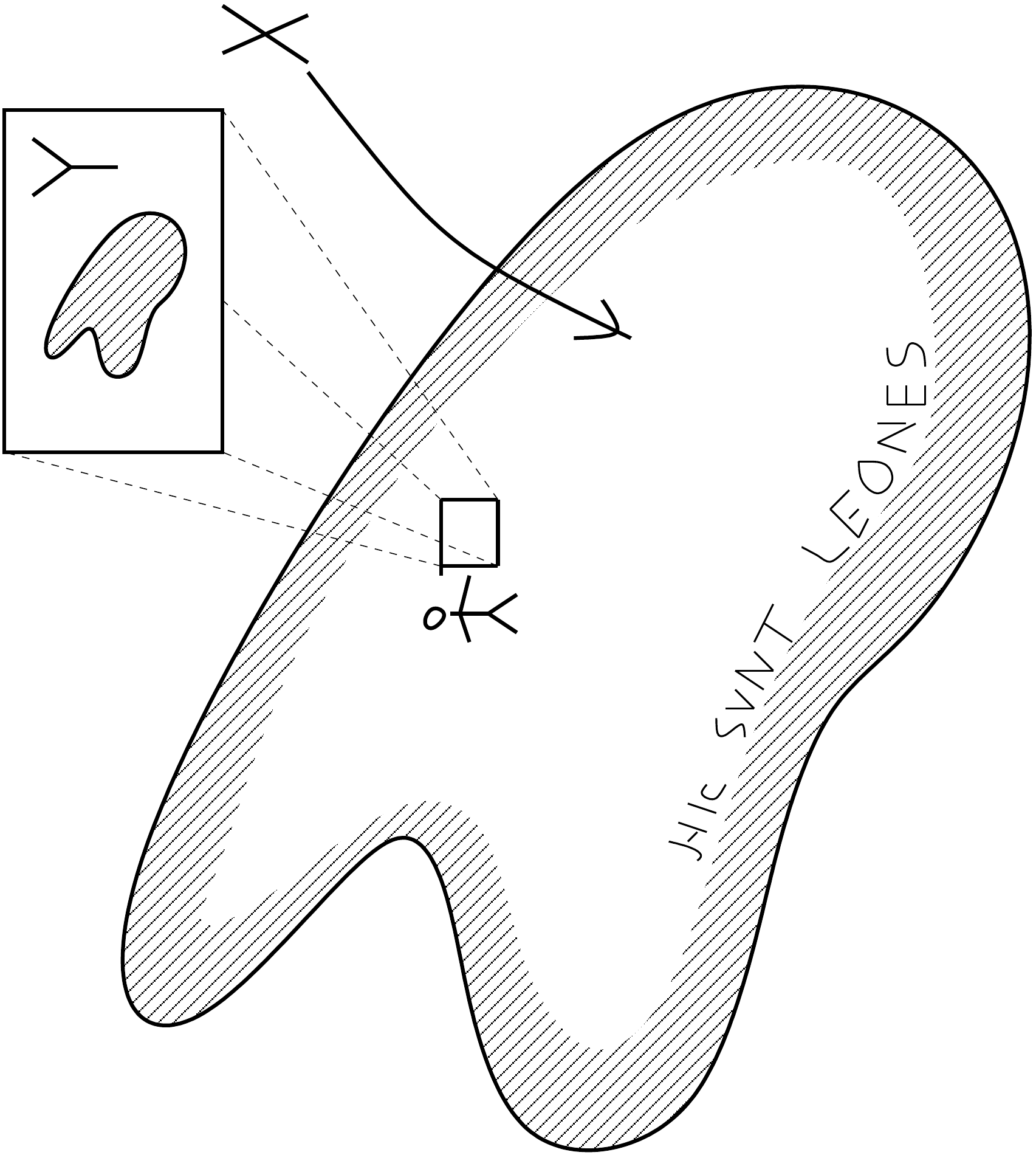}}
%\caption{}
%\label{hicsuntleones1}

\vspace{.5cm}

\paragraph{A too simple model.}
The explorer (call her Alice)  wants to make a chart of 
a newly discovered land $X$ on the piece of paper $Y$ -- or by using a mound 
of clay $Y$,  or by using any mean of recording her discoveries on charts 
done in the metric space $(Y,D)$.  
"Understanding"  the space $X$ (with respect 
to the choice of the "gauge" function $d$) into the terms of the more familiar 
space $Y$   means  making a chart $f:X \rightarrow Y$ 
of $X$ into $Y$ which is not deforming distances too much. Ideally, a 
perfect chart has to be Lipschitz, that is the distances between points in 
$X$ are transformed by the chart into distances between points in $Y$, with a 
precision independent of the scale: the chart  
$f$ is (bi)Lipschitz if there are positive numbers 
$c < C$ such that for any two points 
$x,y \in X$  
$$ c \, d(x,y) \, \leq \, D(f(x),f(y)) \, \leq \, C \, d(x,y)$$
This would be a very  good chart, because it would tell how $X$ is at all scales. 
There are two difficulties related to this model. First, it is impossible to 
make such a chart in practice. What we can do 
instead, is to sample the space $X$ (take a $\varepsilon$-dense subset of $X$ 
with respect to the distance $d$) and try to represent as good as possible this 
subspace in $Y$. Mathematically this is like asking for the chart function $f$ to have the following property: there are supplementary positive constants $a, A$ such that 
for any two points 
$x,y \in X$  
$$ c \, d(x,y) - a  \, \leq \, D(f(x),f(y)) \, \leq \, C \, d(x,y) + A$$
The second difficulty is that such a chart might not exist at all, from
mathematical reasons (there is no quasi-isometry between the metric spaces 
$(X,d)$ and $(Y,D)$). Such a chart exists of course if we want to make charts of 
regions with bounded distance, but remark that all details are erased at small
scale. The remedy would be to make better and better charts, at smaller and
smaller scales, eventually obtaining something resembling a road atlas, with 
charts of countries, regions, counties, cities, charts which have to be compatible 
one with another in a clear sense. 

\subsection{The metaphor of the binocular explorer} 

 Let us go into details of the following exploration program of the space $X$: the explorer Alice jumps randomly in the metric space $(X,d)$, making drafts of  maps at scale $\varepsilon$ and simultaneously orienting herself by using these draft maps. 

I shall explain each part of this exploration problem. 

\paragraph{Making maps at scale $\varepsilon$.} I shall suppose that Alice, while 
at $x \in X$, makes a map at scale $\varepsilon$ of a neighbourhood of 
$x$, called $\displaystyle V_{\varepsilon}(x)$, into another neighbourhood of 
$x$, called $\displaystyle U_{\varepsilon}(x)$. The map is a function: 
$$\delta^{x}_{\varepsilon} : U_{\varepsilon}(x) \rightarrow V_{\varepsilon}(x)$$ 
The "distance on the map" is just: 
\begin{equation}
d^{x}_{\varepsilon} (u,v) \, = \, \frac{1}{\varepsilon} d( \delta^{x}_{\varepsilon} u \, , \, \delta^{x}_{\varepsilon} v)
\label{drel}
\end{equation}
so the map $\displaystyle \delta^{x}_{\varepsilon}$ is indeed a rescaling map, with scale $\varepsilon$.

We suppose of course that the map $\displaystyle \delta^{x}_{\varepsilon}$ is a bijection and we call the inverse of it by the name  $\displaystyle \delta^{x}_{\varepsilon^{-1}}$.  

We may suppose that the map of Alice, while being at point $x$, is a piece of paper 
laying down on the ground, centered such that the "$x$" on paper coincides with the place in $X$ marked by "$x$". In mathematical terms, I ask that 
\begin{equation}
\delta_{\varepsilon}^{x} x \, = \,  x
\label{comp0}
\end{equation}

\paragraph{Jumping randomly in the space $X$.} For random walks we need random  walk kernels, therefore we shall suppose that metric balls in $(X,d)$ have finite, non zero, (Hausdorff) measure. Let $\mu$ be  the associated  Hausdorff measure. The 
$\varepsilon$-random walk is then random  jumping from $x$ into the ball $B(x,\varepsilon)$, so the random walk kernel is  
\begin{equation}
m_{x}^{\varepsilon} \, = \,  \frac{1}{\mu(B(x,\varepsilon))} \,  \mu_{|_{B(x,\varepsilon)}} 
\label{drwk}
\end{equation}
but any Borel probability $\displaystyle m_{x}^{\varepsilon} \in \mathcal{P}(X)$ would be good. 

\paragraph{Compatibility considerations.} At his point I want to add some relations which give a more precise meaning to the scale $\varepsilon$.  

I shall first introduce a standard notation. If $\displaystyle f: X_{1} \rightarrow X_{2}$ is a Borel map and $\displaystyle \mu \in \mathcal{P}(X_{1})$ is a probability measure in 
the space $\displaystyle X_{1}$  then the push-forward of 
$\mu$ through $f$ is the probability measure $\displaystyle f \sharp \mu \in
\mathcal{P}(X_{2}$ defined by: for any $\displaystyle B\in\mathcal{B}(X_{2})$ 
$$\left(f \sharp \mu \right) (B) \, = \, \mu(f^{-1}(B))$$ 

For example,  notice that the random walk kernel is transported  into a random walk kernel by the inverse of the map $\displaystyle \delta^{x}_{\varepsilon^{-1}}$.  I shall impose  that for any (open) set $\displaystyle A \subset  U_{\varepsilon}(x)$ we have 
\begin{equation}
m_{x}^{\varepsilon}(\delta^{x}_{\varepsilon} (A)) \, = \, m_{x}(A) + \mathcal{O}(\varepsilon)
\label{comp1}
\end{equation}
where $\mathcal{O}(\varepsilon)$ is a function (independent on $x$ and $A$) going to zero as $\varepsilon$ is going to zero and $\displaystyle m_{x}$ is a random walk kernel of the map representation $\displaystyle (U_{\varepsilon}(x), d^{x}_{\varepsilon})$, more precisely this is a Borel probability. With the previously induced notation, relation (\ref{comp1}) appears as: 
$$\delta_{\varepsilon^{-1}}^{x} \, \sharp \, m_{x}^{\varepsilon} \, (A) \, = \, 
m_{x}(A) + \mathcal{O}(\varepsilon)$$

Another condition to be imposed is   that  $\displaystyle V_{\varepsilon}$ is approximately the ball $B(x,\varepsilon)$, in the following sense:
\begin{equation}
m_{x}(U_{\varepsilon}(x) \setminus \delta_{\varepsilon^{-1}}^{x} B(x,\varepsilon)) =  \mathcal{O}(\varepsilon)
\label{comp2}
\end{equation}

\paragraph{Multiple drafts reality.} We may see the  atlas that Alice draws as a  "multiple drafts" theory of $(X,d)$. Indeed, the reality at scale $\varepsilon$ 
(and centered at $x \in X$), according to Alice's exploration, can be seen as the triple: 
$$(U_{\varepsilon}(x), d^{x}_{\varepsilon}, \delta_{\varepsilon^{-1}}^{x} \, \sharp \, m_{x}^{\varepsilon})$$
Moreover, we may think that the data used to construct Alice's atlas are: the distance $d$ and for any $x \in X$ and $\varepsilon > 0$,  the functions $\displaystyle \delta_{\varepsilon}^{x}$ and the probability measures 
$\displaystyle m_{x}^{\varepsilon}$. 

For any $x \in X$ and $\varepsilon > 0$ these data may be transported to the 
"reality at scale $\varepsilon$, centered at $x \in X$". Indeed, transporting all 
the data by the function $\displaystyle \delta_{\varepsilon}^{x}$, we get 
$\displaystyle d^{x}_{\varepsilon}$ instead of $d$, and  for any 
$u \in X$ (sufficiently close to $x$ w.r.t. the distance $d$) and any  $\mu >0$,  
we define the relative map making functions: 
\begin{equation}
\delta^{x, u}_{\varepsilon, \mu} \, = \, \delta^{x}_{\varepsilon^{-1}} \, 
\delta_{\mu}^{\delta^{x}_{\varepsilon} u} \, \delta^{x}_{\varepsilon} 
\label{reldil}
\end{equation}
and the relative kernels of random walks: 
\begin{equation}
m_{x, u}^{\varepsilon, \mu} \, = \, \delta^{x}_{\varepsilon^{-1}} \, \sharp \, 
m_{u}^{\mu} 
\label{relrwk}
\end{equation}
With these data I define $reality(x,\varepsilon)$ as: 
\begin{equation}
reality(x,\varepsilon) \, = \, \left(x, \varepsilon, d^{x}_{\varepsilon} , (u,\mu) \mapsto \left( \delta^{x, u}_{\varepsilon, \mu} , m_{x, u}^{\varepsilon, \mu} \right) \right)
\label{real1}
\end{equation}
Remark that now we may repeat this construction and define $reality(x,\varepsilon)(y,\mu)$, and so on. A sufficient condition for having a "multiple drafts" like reality: 
$$reality(x, \varepsilon)(x, \mu) \, = \, reality(x, \varepsilon \mu)$$ 
is to suppose that for any $x \in X$ and any $\varepsilon , \mu > 0$ we have 
\begin{equation}
\delta_{\varepsilon}^{x} \,  \delta^{x}_{\mu} \, = \, \delta^{x}_{\varepsilon \mu}
\label{semigroup}
\end{equation}
For more discussion about this (but without considering random walk kernels) see 
section 2.3 \cite{buligaintro} and references therein about the notion of "metric profiles".

\paragraph{Binocular orientation.} Here I explain how Alice orients herself by using the draft maps from her atlas. Remember that Alice jumps from $x$ to $\displaystyle \delta_{\varepsilon}^{x} u$ (the point $u$ on her map centered at $x$). Once arrived at 
$\displaystyle \delta_{\varepsilon}^{x} u$, Alice draws another map and then she uses a binocular approach to understand what she sees from her new location. She compares the maps 
simultaneously. Each map is like a telescope, a microscope, or an eye. Alice has two eyes. With one she sees $reality(x,\varepsilon)$ and, with the other, $\displaystyle 
reality(\delta_{\varepsilon}^{x}u , \varepsilon)$. These two "retinal images" are compatible in a very specific sense, which gives to Alice 
the sense of her movements. Namely Alice uses the following mathematical fact, which was explained before several times, in the language of: "dilatation structures" \cite{buligadil1} \cite{buligadil2}, "metric spaces with dilations" \cite{buligaintro} (but not using groupoids), 
"emergent algebras" \cite{buligairq} (showing that the distance is not necessary 
for having the result, and also not using groupoids), \cite{buligagr} (using 
normed groupoids). 

\begin{theorem}
There is a groupoid (a small category with invertible arrows) $Tr(X,\varepsilon)$ which has as objects $reality(x,\varepsilon)$, defined at (\ref{real1}), for all 
$x \in X$, and  as arrows the functions: 
\begin{equation}
\Sigma_{\varepsilon}^{x}(u, \cdot) : \, reality(\delta_{\varepsilon}^{x}u , \varepsilon) \rightarrow \, reality(x,\varepsilon) 
\label{real2}
\end{equation}
defined for any $\displaystyle u, v \in U_{\varepsilon}(x)$ by: 
\begin{equation}
\Sigma_{\varepsilon}^{x}(u,v) \, = \, \delta^{x}_{\varepsilon^{-1}} \, \delta^{\delta^{x}_{\varepsilon} u}_{\varepsilon}  v
\label{sig1}
\end{equation}
 Moreover, all arrows are 
isomorphisms (they are isometries and transport in the right way the structures defined at (\ref{reldil}) and (\ref{relrwk}), from the source of the arrrow to the 
target of it). 
\label{t1}
\end{theorem}

Therefore Alice, by comparing the maps she has, orients herself by using the 
"approximate translation by $u$" 
$\displaystyle \Sigma_{\varepsilon}^{x}(u, \cdot)$ 
(\ref{sig1}), which appears as an isomorphism between the two realities, cf. (\ref{real2}). This function deserves the name "approximate translation", as proved elsewhere \cite{buligadil1} \cite{buligadil2} \cite{buligairq}, because of 
the following reason: if we suppose that $reality(x,\varepsilon)$ converges 
as the scale $\varepsilon$ goes to $0$, then $\displaystyle \Sigma_{\varepsilon}^{x}(u, \cdot)$ should converge to the translation by $u$, as seen  in the tangent space as $x$. Namely we have the following mathematical definition and theorem  (for the proof see the cited references; the only new thing concerns the convergence of the random walk kernel, but this is an easy consequence of the compatibility conditions (\ref{comp1}) (\ref{comp2})). 

\begin{definition}
A metric space with dilations  and random walk (or dilatation structure with a random walk) is a structure $(X,d,\delta,m)$ such that:  
\begin{enumerate}
\item[(a)] $(X,d,\delta)$ is a normed uniform idempotent right quasigroup  (or  dilatation structure), 
cf. \cite{buligabraided} Definition 7.1. 
\item[(b)] $m$ is a random walk kernel, that is for every $\varepsilon > 0$ 
we have a measurable function $\displaystyle x \in X \, \mapsto \, m_{x}^{\varepsilon} \in \mathcal{P}(X)$ (a transversal function on the pair groupoid 
$X \times X$, in the sense of \cite{connes} p. 35) which satisfies the compatibility conditions (\ref{comp1}) (\ref{comp2}). 
\end{enumerate}
\label{de1}
\end{definition}

\begin{theorem}
Let $(X,d,\delta,m)$ be a   dilatation structure with a random walk. Then for any 
$x \in X$, as the scale $\varepsilon$ goes to $0$ the structure $reality(x,\varepsilon)$ converges to $reality(x,0)$ defined by: 
\begin{equation}
reality(x,0) \, = \, \left(x, 0, d^{x} , (u,\mu) \mapsto \left( \delta^{x, u}_{ \mu} , m_{x, u}^{ \mu} \right) \right)
\label{realtang}
\end{equation}
in the sense that $\displaystyle d^{x}_{\varepsilon}$ converges uniformly 
 to $\displaystyle d^{x}$, dilations $\displaystyle \delta^{x, u}_{\varepsilon,  \mu}$ converge uniformly to dilations $\displaystyle \delta^{x, u}_{ \mu}$, 
 probabilities  $\displaystyle m_{x, u}^{\varepsilon, \mu}$  converge simply  
 to probabilities $\displaystyle  m_{x, u}^{ \mu}$. Finally $\displaystyle 
 \Sigma_{\varepsilon}^{x}(\cdot , \cdot)$ converges uniformly to a conical group 
 operation. 
 \label{t2}
 \end{theorem}
 
Conical groups are a (noncommutative and vast) generalization of real vector 
spaces. See for this \cite{buligadil2}. 

\paragraph{Alice knows the differential geometry of $X$.} All is hidden in dilations (or maps made by Alice) $\displaystyle \delta_{\varepsilon}^{x}$. They encode the approximate (theorem \ref{t1}) and 
exact (theorem \ref{t2}) differential AND algebraic structure of the tangent 
bundle of $X$.

\paragraph{Alice communicates with Bob.} Now Alice wants to send Bob her knowledge. 
Bob lives in the metric space $(Y,D)$, which was well explored previously, therefore 
himself knows well the differential geometry and calculus in the space $Y$. 

Based on Alice's informations, Bob hopes to construct a Lipschitz function 
from an open set in $X$ with image an open set in $Y$. If he succeeds, then he 
will know all that  Alice knows (by transport of all relevant structure using 
his Lipschitz function). But he might fail, because of a strong mathematical 
result . Indeed, as a consequence of Rademacher theorem 
(see Pansu \cite{pansu}, for a 
Rademacher theorem for Lipschitz function between Carnot groups, a type of 
conical groups), if such a function would exist then it would be 
differentiable almost everywhere. The differential (at a point) is a 
morphism of conical groups which 
commutes with Alice's and Bob's maps (their respective dilations). Or, it 
might happen that there is no non-trivial such morphism and we arrrive at a 
contradiction! 

This is not a matter of topology (continuous or discrete, or 
complex topological features at all scales), but a matter having to do with the 
construction of the "realities" (in the sense of relation (\ref{real1}), 
having  all to do with differential calculus and differential geometry 
in the most fundamental sense.

 Therefore, another strategy of Alice would be to communicate to Bob not 
 all details of her map, but the relevant algebraic identities that her maps 
 (dilations) satisfy. Then Bob may try to simulate what Alice saw. This is  
 a path first suggested in  \cite{buligafrontend}.

\subsection{Simulating spaces}
Let us consider  some examples of spaces, like: the real world, a virtual world of a game, mathematical spaces as manifolds, fractals, symmetric spaces, groups, linear spaces ...
     
All these  spaces  may be characterized by the class of algebraic/differential computations which are possible, like: zoom into details, look from afar, describe velocities and perform other differential calculations needed for describing the physics of such a space, perform reflexions (as in symmetric spaces), linear combinations (as in linear spaces), do affine or projective geometry constructions and so on.
 
 Suppose that on a set $X$ (called "a space") there is an operation
$$ (x,y) \mapsto x \circ_{\varepsilon} y $$
which is dilation-like. Here $\varepsilon$ is a parameter belonging to a 
commutative group $\Gamma$, for simplicity let us take $\Gamma = (0,+\infty)$ 
with multiplication as the group operation. 

By "dilation-like" I mean the following: 
for any $x \in X$  the function 
$$ y \mapsto \delta^{x}_{\varepsilon} y \, = \, x \circ_{\varepsilon} y$$
behaves like a dilation in a vector space, that is 
\begin{enumerate}
\item[(a)] for any $\varepsilon, \mu \in (0,+\infty)$ we have 
$\displaystyle  \delta^{x}_{\varepsilon} \,  \delta^{x}_{\mu} \, = \,  \delta^{x}_{\varepsilon \mu}$ and $\displaystyle  \delta^{x}_{1} \, = \, id$; 
\item[(b)] the limit as $\varepsilon$ goes to $0$ of $\displaystyle 
 \delta^{x}_{\varepsilon} y$ is $x$, uniformly with respect to $x,y$. 
\end{enumerate} 

Then the dilation operation is the basic building block of both the 
algebraic structure of the space (operations in the tangent spaces) and the 
differential calculus in the space, as it will be explained further. This leads to he introduction 
of emergent algebras \cite{buligairq}. 

Something amazing happens if we take compositions of dilations, like this 
ones 
$$\Delta^{x}_{\varepsilon} (u,v) \, = \,  \delta_{\varepsilon^{-1}}^{\delta^{x}_{\varepsilon} u} \, \delta^{x}_{\varepsilon} u$$
$$\Sigma^{x}_{\varepsilon} (u,v) \, = \, \delta^{x}_{\varepsilon^{-1}} \, 
\delta_{\varepsilon}^{\delta^{x}_{\varepsilon} u} v $$
called the approximate difference, respectively approximate sum operations based 
at $x$. If we suppose that $\displaystyle(x,u,v) \mapsto \Delta^{x}_{\varepsilon} (u,v)$ and $\displaystyle(x,u,v) \mapsto \Sigma^{x}_{\varepsilon} (u,v)$ converge 
uniformly as $\varepsilon \rightarrow 0$ to $\displaystyle(x,u,v) \mapsto \Delta^{x} (u,v)$ and $\displaystyle(x,u,v) \mapsto \Sigma^{x} (u,v)$ then, out of apparently 
nothing, we get that $\displaystyle \Sigma^{x}$ is a group operation with $x$ as neutral element, which can be interpreted as the operation of vector addition 
in the "tangent space at $x$" (even if there is no properly defined such space). 

To convince you about this just look at the following example: $X = \mathbb{R~}^{n}$ and 
$$x \circ_{\varepsilon} y \,  = \, \delta^{x}_{\varepsilon} y \,  = \, x + \varepsilon (-x+y)$$ 
Then the approximate difference and sum operations based at $x$ have the expressions: 
\begin{align*}
\Delta_{\varepsilon}^{x}(u,v) 
&= x+\varepsilon(-x+u) + (-u+v)  & \rightarrow & & x -u+v \\
\Sigma_{\varepsilon}^{x}(u,v) 
&=  u+ \varepsilon(-u+x) + (-x+v)  & \rightarrow & & u -x+v 
\end{align*}

Moreover, with the same dilation operation we may define something resembling very 
much with differentiation. Take a function $f: X \rightarrow X$, then define 
$$ D_{\varepsilon} f (x) u \, = \, \delta_{\varepsilon^{-1}}^{f(x)} f \delta^{x}_{\varepsilon} (u)$$  
In the particular example used previously we get 
$$ D_{\varepsilon} f (x) u \, = \, 
f(x) + \frac{1}{\varepsilon} ( -f(x)+f(x+\varepsilon(-x+u)))$$
which shows that the limit as $\varepsilon$ goes to $0$ of $\displaystyle 
u \mapsto  D_{\varepsilon} f (x) u$ is a kind of differential.

Such computations are finite or virtually infinite "recipes", which can be implemented by some class of circuits made by very simple gates based on dilation 
operations (as in boolean computing, where transistors are universal gates for computing boolean functions).
     
(Computation in) a space is then described by emergent algebras \cite{buligairq},  which are inspired by the considerations about a formal calculus with binary decorated planar trees in relation with dilatation structures \cite{buligadil1}:
     
A - a class of transistor-like gates, with in/out ports labelled by points of the space and a internal state variable which can be interpreted as "scale". I propose dilations as such gates (basically these are idempotent right quasigroup operations).
     
B - a class of elementary circuits made of such gates (these are the "generators" of the space). The elementary circuits have the property that the output converges as the scale goes to zero, uniformly with respect to the input. 

C - a class of equivalence rules saying that some simple assemblies of elementary circuits have equivalent function, or saying that relations between those simple assemblies converge to relations of the space as the scale goes to zero.

Seen like this, "simulating a space" means: give a set of transistors (and maybe some non-emergent operations and relations), elementary circuits and relations which are sufficient to generate any interesting computation in this space.
     
For the moment I know how to simulate affine spaces \cite{buligadil2},  sub-riemannian  or Carnot-Carath\'eodory spaces \cite{buligadil3}, riemannian \cite{buligairq} or sub-riemannian symmetric spaces \cite{buligabraided}.

\subsection{Axial maps in terms of choroi and carriers}

There is yet another field which may be relevant for space computation:
architecture.  Here several names 
come to mind, like Christopher Alexander and Bill Hillier\footnote{ for an intro
see 
Diagrammatic Transformation of Architectural Space, Kenneth J. Knoespel
Georgia Institute of Technology,
\url{http://www.lcc.gatech.edu/~knoespel/Publications/DiagTransofArchSpace.pdf}}.
Especially the notion of "axial map", which could be just an embodiment of 
a dilation structure, is relevant \cite{hil1}, \cite{hil2}. The notion of 
"axial map" still escapes a rigorous mathematical definition, but seems to be
highly significant in order to understand emergent social 
behaviour\footnote{\url{http://en.wikipedia.org/wiki/Space_syntax}} 
.

From the  article \cite{hil2}: 

"Hillier and Hanson (1984) noted that urban space in
particular seems to comprise two fundamental elements: "stringiness" 
and "beadiness" 
such that the space of the systems tends to resemble beads on a 
string. They write: "We can define "stringiness" as being to do with 
the extension of space in one dimension, whereas "beadiness" is to do 
with the extension of space in two dimensions" (page 91).

Hence, the epistemology of their methodology involves the investigation 
of how space is constructed in terms of configurations of interconnected 
beads and strings. To this end, they develop a more formal definition of 
the elements in which strings become "axial lines" and beads "convex spaces". 
The definition they give is one that is easily understood by human 
researchers, but which, it has transpired, is difficult to translate 
into a computational approach: "An axial map of the open space structure 
of the settlement will be the least set of [axial] lines which pass 
through each convex space and makes all axial links"
(Hillier and Hanson, 1984, pages 91-92)."

In my opinion, a good embodiment of an axial map is a dilation structure. In
fact, a weaker version of a dilation structure would be enough. In this
formalism the axial lines are differences and the beads are choroi. 

In this case a tangle diagram superimposed on a geographical map of  a place 
(like  a city) would explain also how to move in this place. Differences are 
like carriers along the architect' axial lines and beads are convex places 
(because of the mathematically obvious result that a subset of the plane 
is closed with respect to taking dilations of coefficients smaller than 1 if and
only if it is convex). An axial map is thus a skeleton for the freedom of
movement in a place.

\subsection{Spacebook}
 
(Name coined by Mark Changizi in a mail exchange, after seeing the first version
of this paper). 

How to make your library my library? see Changizi article "The Problem With the Web and
E-Books Is That There's No Space for 
Them"\footnote{\url{http://www.psychologytoday.com/blog/nature-brain-and-culture/201102/the-problem-the-web-and-e-books-is-there-s-no-space-them}}
The problem in fact is that to me the brain-spatial interface of Changizi 
library is largely incomprehensible. I have to spend time in order to 
reconstruct it in my head. The idea is then to do a "facebook" for space
competences. How to share spatial computations?

%\newpage

\section{Colorings of tangle diagrams}

The  idempotent right quasigroups are related to algebraic structures 
appearing in knot theory.   J.C. Conway and G.C. Wraith, in their unpublished 
correspondence  from 1959,  used the name "wrack" for  a  self-distributive 
right quasigroup generated by a link diagram. Later, Fenn and Rourke \cite{fennrourke} proposed the name "rack" instead. Quandles are particular case of racks, namely self-distributive idempotent right quasigroups. They were introduced by  Joyce \cite{joyce}, as a distillation  of the  Reidemeister moves. 

The axioms of a (rack ;  quandle ; irq) correspond respectively to the (2,3 ; 1,2,3 ; 1,2) Reidemeister moves. That is why we shall use decorated braids diagrams in order to explain what emergent algebras are.

The basic idea of racks and quandles is that these are algebraic operations 
related to the coloring of tangles diagrams. 

\subsection{Oriented tangle diagrams and trivalent graphs}

Visually, a oriented tangle diagram is the result of a regular projection 
on a plane of a properly embedded in the 3-dimensional space, oriented, 
one dimensional  manifold, together with additional over- and 
under-information at crossings (adapted from the "Tangle, relative
link" article from Encyclopaedia of mathematics. Supplement. Vol. III. Edited by
M. Hazewinkel. Kluwer Academic Publishers, Dordrecht, 2001, page 395). 

Because the tangle diagram is oriented, there are two types of crossings, 
indicated in the next figure. 

\vspace{.5cm}

\centerline{\includegraphics[angle=270, width=0.5\textwidth]{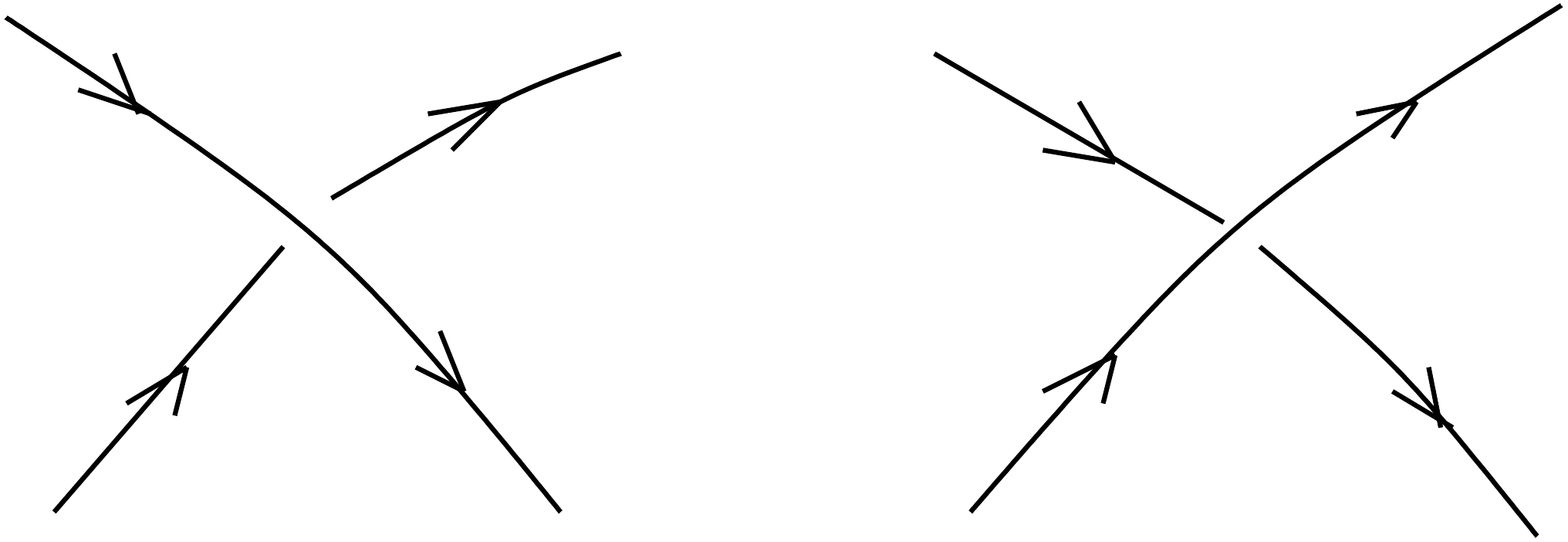}}

\vspace{.5cm}

The tangle which projects to the tangle diagram is to be seen as a
"parameterization" of the tangle diagram. In this sense, by using 
the image of the tangle diagram  with over- and under-crossings, it is 
easy to define an "arc" of the diagram as the projection of a (part of a) 
1-dimensional embedded manifold from the tangle. Arcs can be open or closed. 
An arc decomposes into connected parts (again by using the first image of 
a tangle diagram) which are called segments. 

An input segment is a segment which enters (with respect to its orientation) 
 into a crossing but it does not exit from a crossing. Likewise, an output segment is one which
exits from a crossing but there is no crossing where this segment enters. 

Oriented tangle diagrams are considered only up to continuous deformations of
the plane. 

Further I shall use, if necessary, the name "first image of a tangle diagram",
if the oriented tangle diagram is represented with this convention. 
There is a "second image", which I explain next.  
In a sense, the true image of a tangle diagram  is the second one, mainly 
because in this paper the fundamental object is the oriented tangle 
diagram, NOT the tangle which projects on the plane to the oriented tangle 
diagram.

For the second image I use  a chord diagram, or a Gauss diagram type of 
indication of crossings. See \cite{duzhin} \cite{kassel} for more on the mathematical aspects of chord diagrams. 

We may see a crossing like a gate, or black box
 with two inputs and two outputs. We open the black box and inside we find a
 combination of two simpler gates, among the following three available: the 
 FAN-OUT, the $\circ$ and the $\bullet$ gate. 

\vspace{.5cm}

\centerline{\includegraphics[angle=270, width=0.7\textwidth]{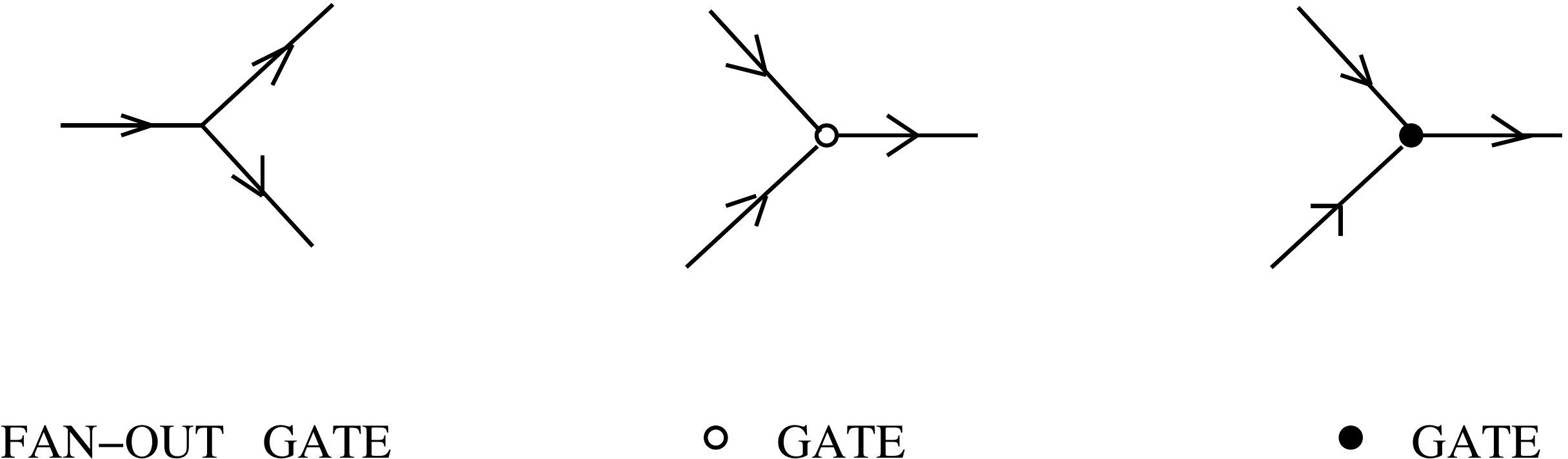}}

\vspace{.5cm}

By using these gates we may transform the oriented tangle diagram into a 
oriented planar graph with 3-valent 2-valent or 1-valent nodes, that is into 
a circuit made only with these gates, connected by wires which could cross 
(crossings of wires in this graph has no meaning). The graph is planar in 
the sense that each trivalent node, which represents one of these 3 gates, 
inherits an orientation coming from the plane. A trivalent node is either 
undecorated, if it belongs to a FAN-OUT gate, or decorated by a  $\circ$ or  
a $\bullet$, if it belongs to a $\circ$ gate or a $\bullet$ gate respectively. 

This graph is obtained by the following procedure. 1-valent nodes represent  
input or output tangle segments. 2-valent 
nodes are used for closing an arc.  
These 2-valent nodes can be replaced by 3-valent nodes 
(corresponding to FAN-OUT
gates), with the price of introducing also a 1-valent node. 
Crossings are replaced by combinations of trivalent nodes, as explained in 
the next figure. 
  
\vspace{.5cm}

\centerline{\includegraphics[angle=270, width=0.7\textwidth]{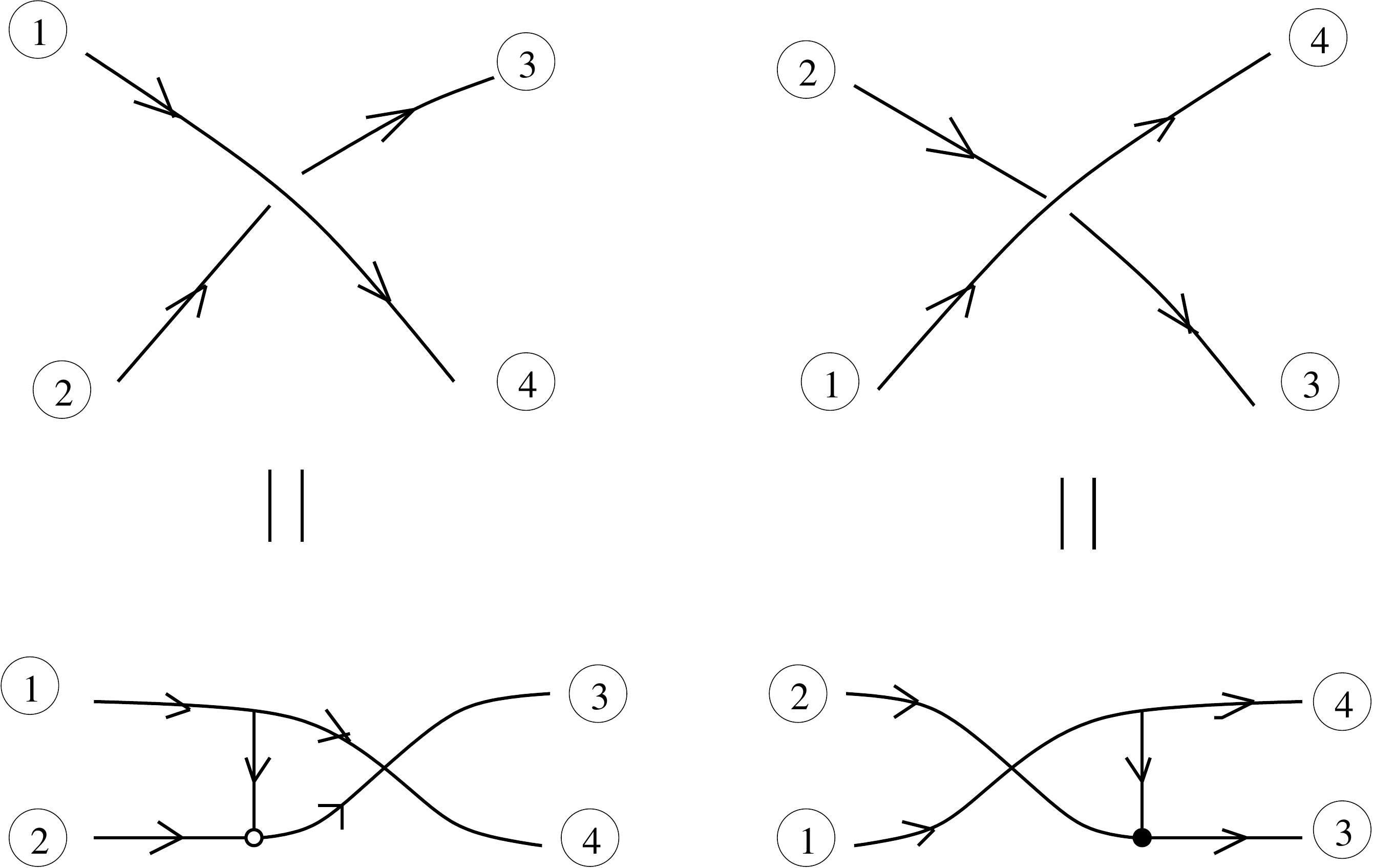}}

\vspace{.5cm}

In the second image of a tangle diagram, all segments are wires connecting the 
nodes of the graph, but some wires are not segments, namely the ones which connect a FAN-OUT with the
corresponding $\circ$ gate or $\bullet$ gate (i.e. with the gate which
constitutes together with the respective FAN-OUT a coding for a over- or
under-crossing). The wires which are not segments are called "chords".

Arcs appear as connected unions of segments with compatible orientations (such
that we can choose a segment of the arc and then  walk along the whole arc, 
by following the local  directions indicated by each segment). By walking along 
an arc, we can recognize the crossings: undecorated nodes correspond to
over-crossings and decorated nodes to under-crossings. 

\vspace{.5cm}

\centerline{\includegraphics[angle=270, width=0.7\textwidth]{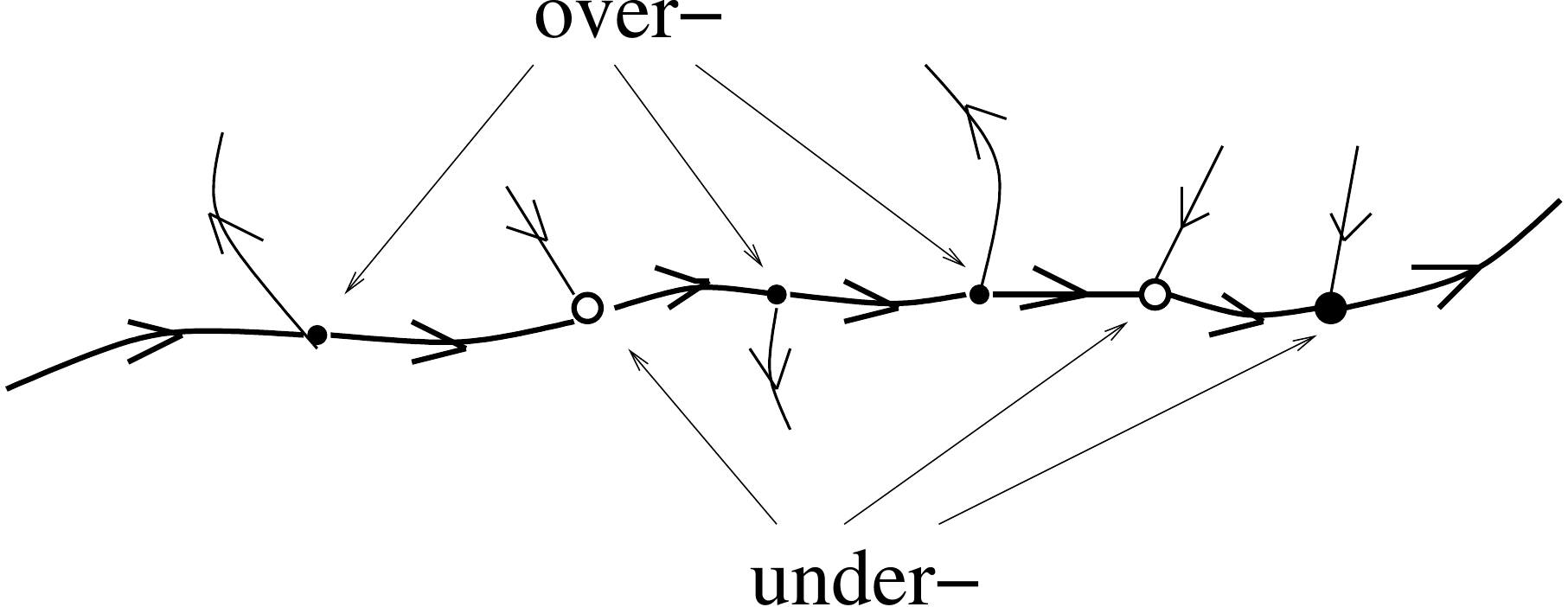}}

\vspace{.5cm}

As a circuit made by gates (in the second image), a tangle diagram appears as 
an ordered list of its crossings gates, each crossing gate being given as a pair of FAN-OUT and one of the other 
two gates.  1-valent nodes appear as input nodes 
(in a separate INPUT list), or as output nodes (in the OUTPUT list). 
The wiring is given as a matrix $M$ of connectivity, namely 
the element $\displaystyle M_{ij}$ corresponding to the pair of  nodes $(i,j)$ 
(from the trivalent and 1-valent graph, independently on the pairing of nodes 
given by the crossings) is equal to $1$ if there is a wire oriented from $i$ to 
to $j$, otherwise is equal to $0$. The definition is unambiguous because from 
$i$ to $j$ can be at most one oriented wire.

\subsection{Colorings with idempotent right quasigroups}

Let $X$ be a set of colors which will be used to 
decorate the  segments in a (oriented) tangle diagram.  There are two binary operations on $S$ related 
to the coloring, as shown in the next figure. 

\vspace{.5cm}

\centerline{\includegraphics[angle=270, width=0.7\textwidth]{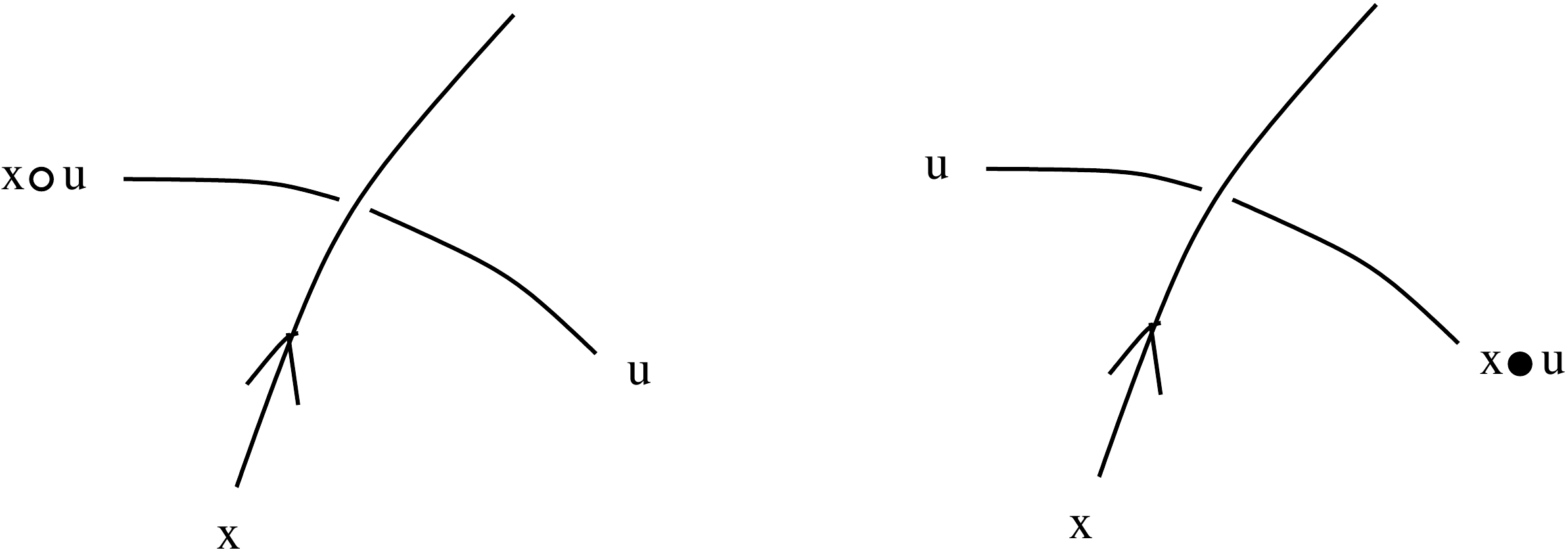}}

\vspace{.5cm}

Notice that only matters the orientation of the arc which passes over.

We have therefore a set $X$ endowed with two operations $\circ$ and $\bullet$.  
We want these operations to satisfy some conditions which ensure that the
decoration of the segments of the tangle diagram rest unchanged after 
performing a Reidemeister move I or II on the tangle diagram. This is explained
further. We show only the part of a (larger) diagram which changes during a 
Reidemeister move, with the convention that what is not shown will not change
after the Reidemeister move is done.

The first condition, related to the Reidemeister move I,
 is depicted in the next figure. 

%\vspace{.5cm}

\centerline{\includegraphics[angle=270, width=0.5\textwidth]{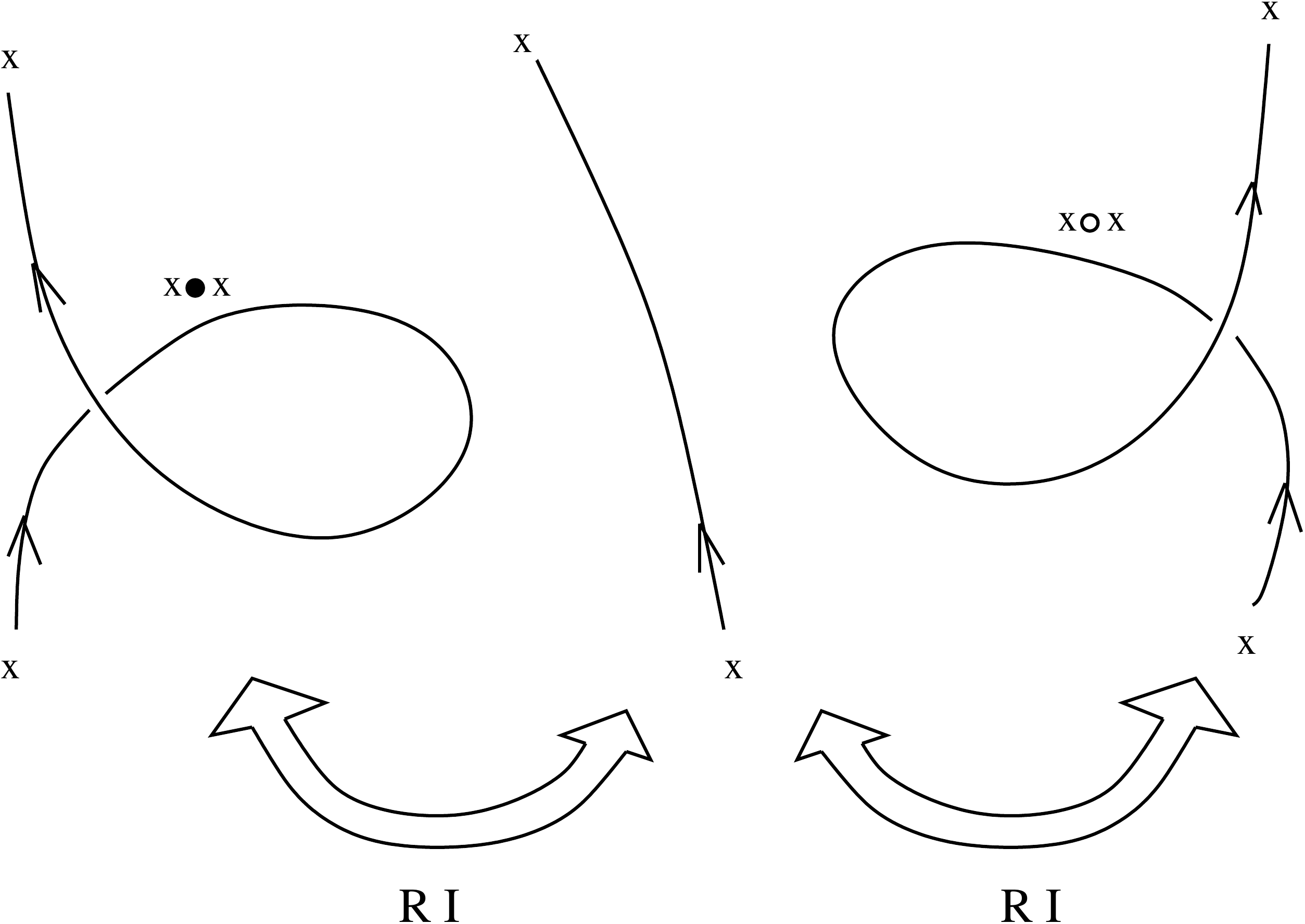}}

\vspace{.5cm}

It means that we can decorate "tadpoles" such that we may remove them 
(by using the Reidemeister move I) afterwards. In algebraic terms, this
condition means that we want the operations $\circ$ and $\bullet$ to be 
idempotent: 
$$x \circ x \, = \, x \bullet x \, = \, x$$
for all $x \in X$.

The second condition is related to the Reidemeister II move. It means that we 
can decorate the segments of a pair of arcs as shown in the following 
picture, in such a way that we can perform the Reidemeister II move and
eliminate a pair of "opposite" crossings.

This condition translates in algebraic terms into saying that 
 $(X,\circ,\bullet)$ is  a right quasigroup. Namely we want  that 
$$x \circ \left( x \bullet y\right) \, = \, x \bullet \left( x \circ y \right) \, = \, y$$
for all $x, y \in X$. This is the same as asking that for any $a$ and $b$ in $X$, the equation $a \circ x = b$ has a solution, which is unique, then denote the solution by $x = a \bullet b$.  
All in all, a set $(X,\circ,\bullet)$ which has the properties related to the first two Reidemeister 
moves is called an idempotent right quasigroup, or irq for short.

\vspace{.5cm}

\centerline{\includegraphics[angle=270, width=0.5\textwidth]{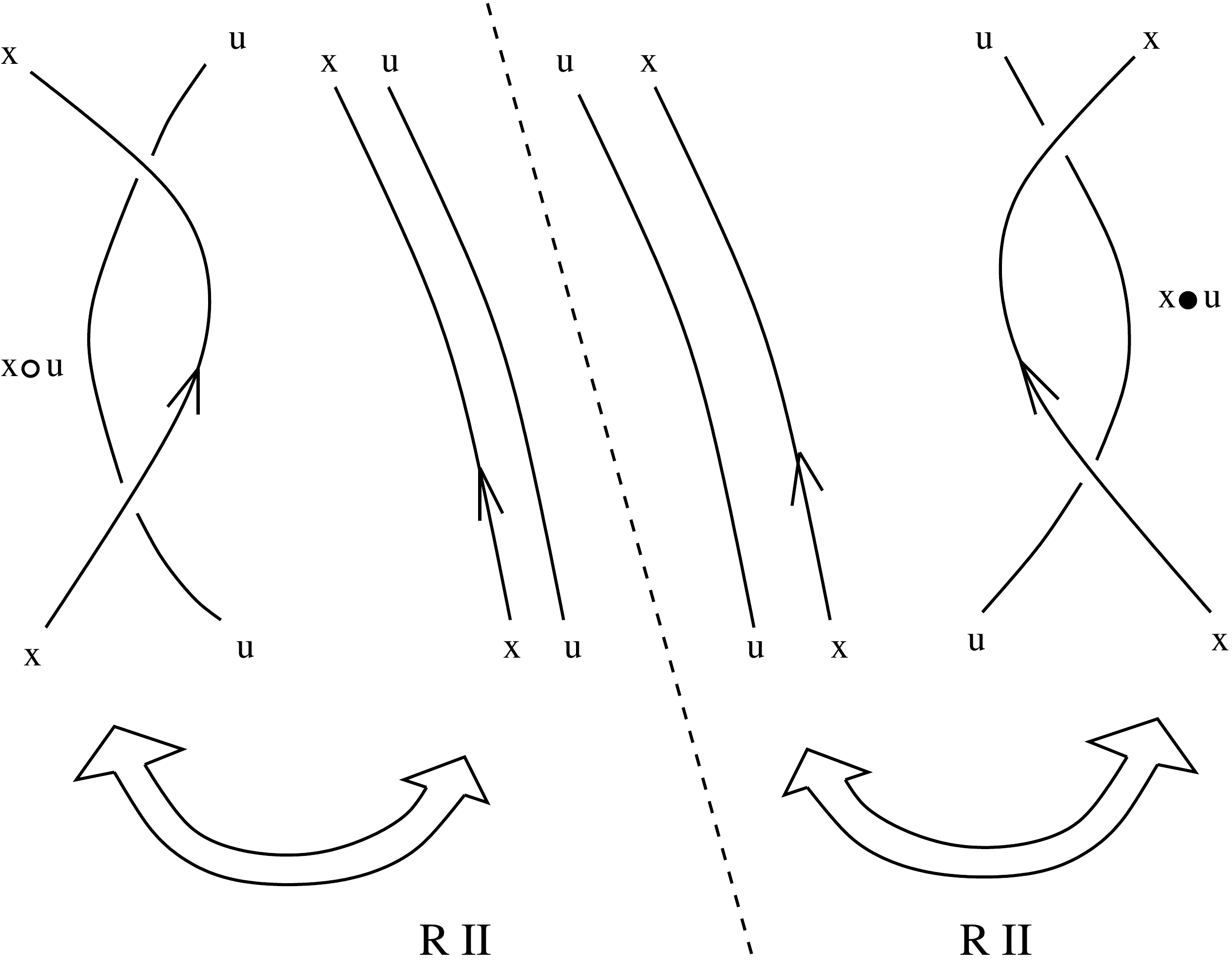}}

\vspace{.5cm}

\begin{definition} A right quasigroup is a set $X$ with a binary operation 
$\circ$ such that for each $a, b \in X$ there exists a unique $x \in X$ such that 
$a \, \circ \, x \, = \, b$. We write the solution of this equation 
$x \, = \, a \, \bullet \, b$. 

 An idempotent right quasigroup (irq) is a  right quasigroup $(X,\circ)$ such that 
 for any $x \in X$ $x \, \circ \, x \, = \, x$. Equivalently, it can be seen as a 
  set $X$ endowed with two  operations $\circ$ and $\bullet$, which satisfy the following axioms: for any $x , y \in X$  
\begin{enumerate}
\item[(R1)] \hspace{2.cm} $\displaystyle x \, \circ \, x \, = \, x \, \bullet \, x \,  = \,  x$
\item[(R2)] \hspace{2.cm} $\displaystyle x \, \circ \, \left( x\, \bullet \,  y \right) \, = \, x \, \bullet \, \left( x\, \circ \,  y \right) \, = \, y$
\end{enumerate}
\label{defquasigroup}
\end{definition}

The Reidemeister III move concerns the sliding of an arc (indifferent of 
orientation) under a crossing. In the next figure it is shown only one 
possible sliding movement. 

\vspace{.5cm}

\centerline{\includegraphics[angle=270, width=0.7\textwidth]{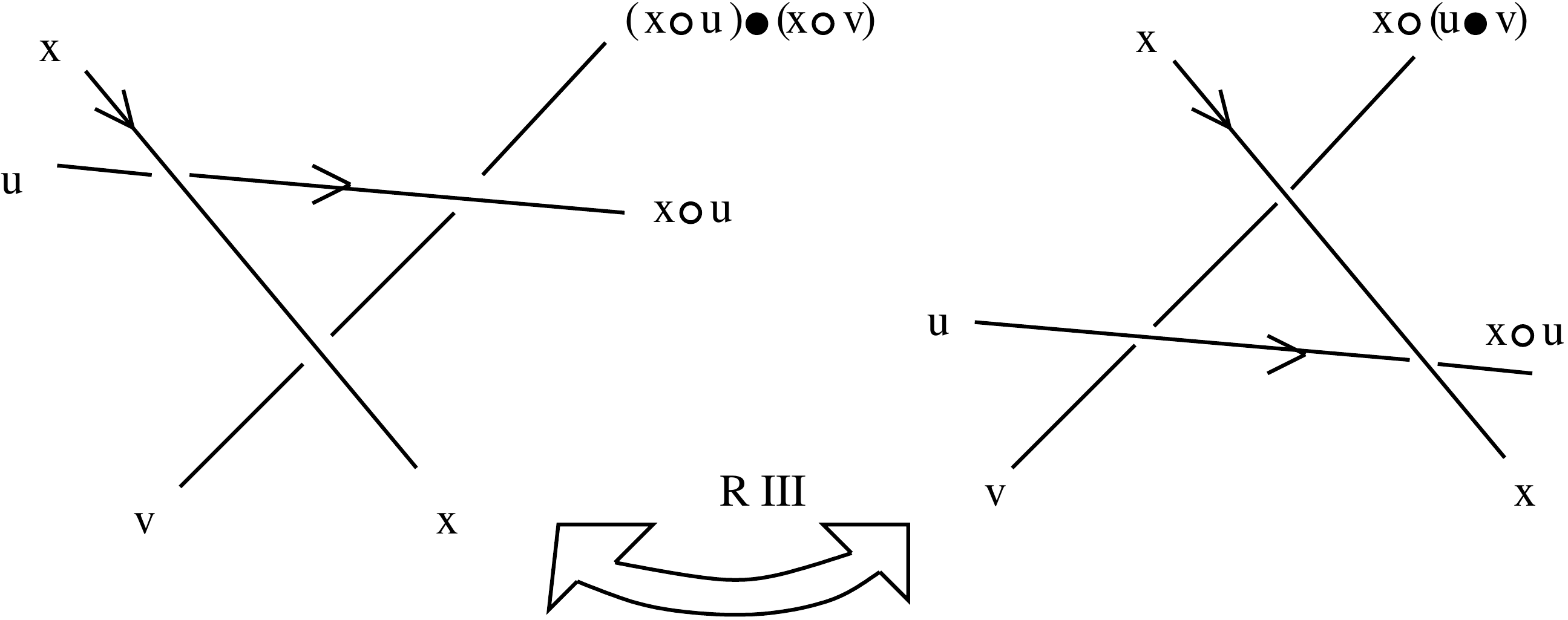}}

\vspace{.5cm}

Such a sliding move is possible, without modifying the coloring, if and 
only if  the operation 
 $\circ$ is left distributive with respect to the operation 
 $\bullet$. (For the other possible choices of crossings, the "sliding" movement
 corresponding to the Reidemeister III move is possible if and only if 
 $\bullet$ is distributive with respect to $\circ$ and also the operations 
 $\circ$ and $\bullet$ are self distributive). 

With this self-distributivity property, $(X,\circ,\bullet)$ is called a quandle. A well known quandle (therefore also an irq) 
is the Alexander quandle: consider $X = \mathbb{Z}[ \varepsilon, \varepsilon^{-1}]$ with the operations 
$$ x \circ y \, = \, x + \varepsilon \left( -x + y\right) \quad , \quad x \bullet y \, = \, x + \varepsilon^{-1} \left( -x + y\right)$$
The operations in the Alexander quandle are therefore dilations in euclidean spaces.

\paragraph{Important remark.} Further  I shall NOT see oriented
tangle diagrams  as objects associated to a tangle in three-dimensional space. 
That is because I am going to renounce to the Reidemeister III move. 
This interpretation, of being projections of tangles in space, is 
only for keeping a visually based vocabulary, like "over", "under", "sliding an 
arc under another" and so on.

\subsection{Emergent algebras and tangles with decorated crossings}

I shall adapt the tangle diagram coloring, presented in the previous 
section, for better understanding of the formalism of dilation structures. 
In fact we shall arrive to a more algebraic concept, more basic in some sense
that the one of dilation structures, named "emergent algebra". 

The first step towards this goal is to consider richer decorations as
previously. We could decorate not only the connected components of 
tangle diagrams  but also the crossings. I use for crossing decorations the
scale parameter. Formally the scale parameter belongs to a commutative group 
$\Gamma$. In this paper is comfortable to think that $\Gamma = (0,+\infty)$ with
the operation of multiplication or real numbers.  

Here is the rule of decoration of tangle diagrams, by using a dilation
structure: 

\vspace{.5cm}

\centerline{\includegraphics[angle=270, width=0.7\textwidth]{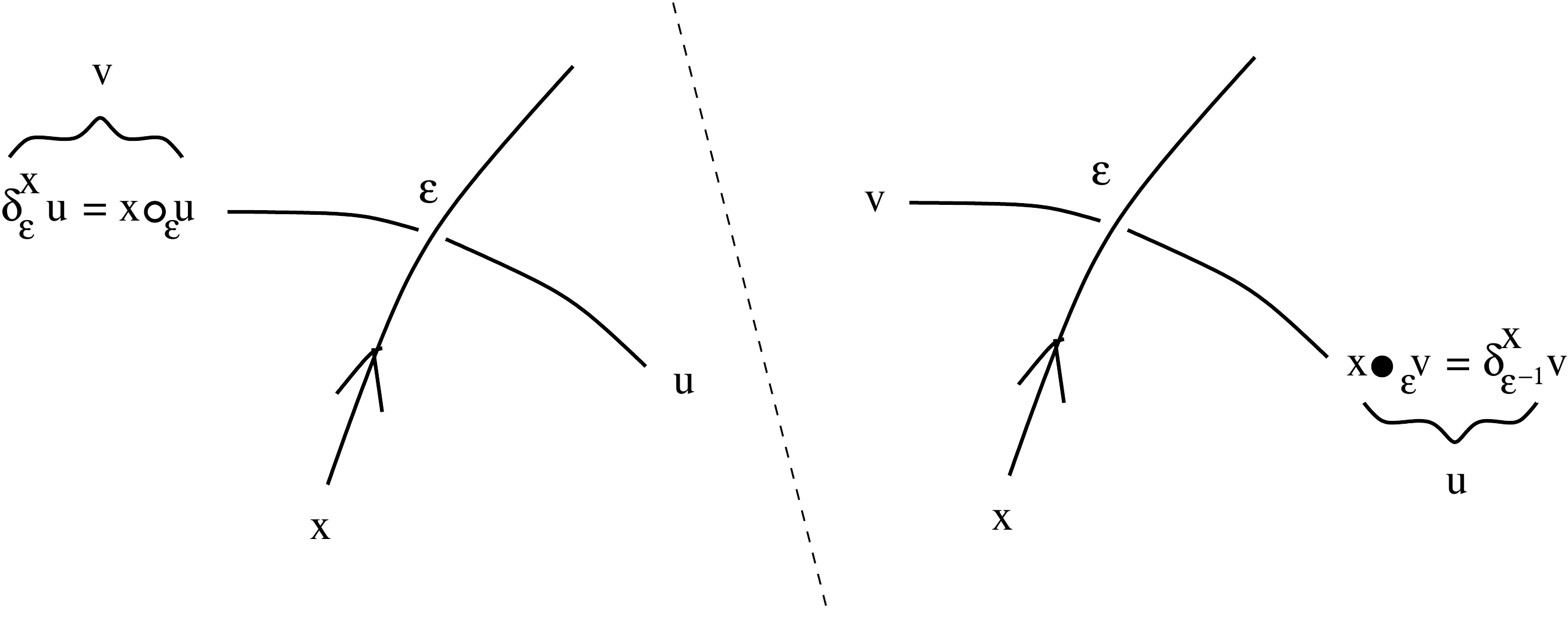}}

\vspace{.5cm}

In terms of idempotent right quasigroups, instead of one $(X,\circ, \bullet)$,
we have a family  
$\displaystyle (X, \circ_{\varepsilon}, \bullet_{\varepsilon})$, for all 
$\varepsilon \in \Gamma$. In terms of dilation structures, the operations are: 
$$x \circ_{\varepsilon} u \, = \, \delta^{x}_{\varepsilon} u \quad , \quad 
x \bullet_{\varepsilon} u \, = \, \delta^{x}_{\varepsilon^{-1}} u$$

This implies that virtual crossings are allowed. A virtual crossing is just a crossing where
nothing happens, a crossing with decoration $\varepsilon = 1$.

\vspace{.5cm}

\centerline{\includegraphics[angle=270, width=0.7\textwidth]{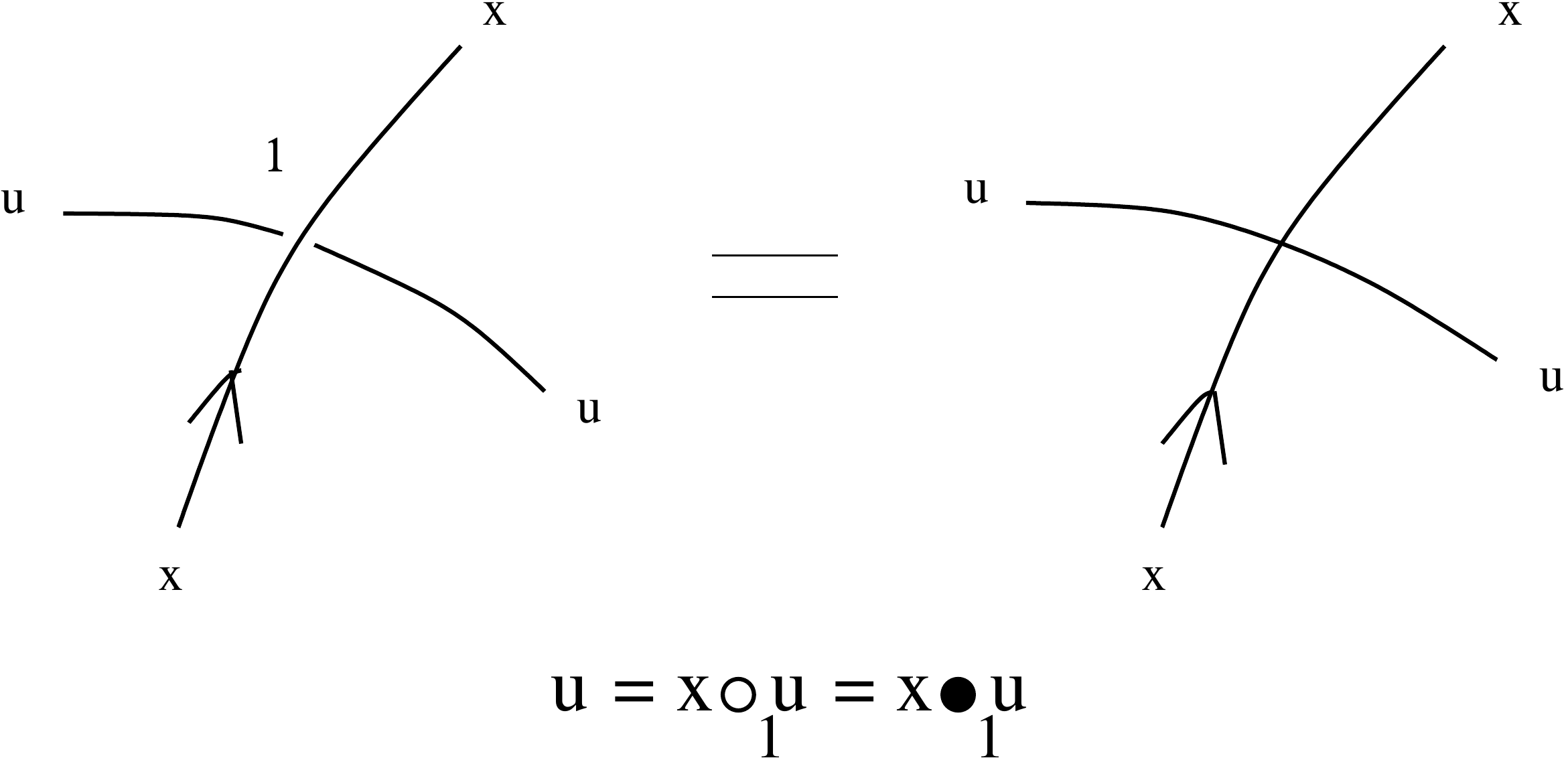}}

\vspace{.5cm}

Equivalent with the first two axioms of dilation structures, 
is that for all $\varepsilon \in \Gamma$  the 
triples $\displaystyle (X, \circ_{\varepsilon}, \bullet_{\varepsilon})$ 
 are idempotent right quasigroups (irqs), moreover we want that for any 
 $x \in X$ the mapping 
$$ \varepsilon \in \Gamma \,  \mapsto \, x \circ_{\varepsilon} \, ( \cdot )$$
to be an action of $\Gamma$ on $X$. This reflects into the following rules for  
combinations of decorated crossings.

\vspace{.5cm}

\centerline{\includegraphics[angle=270, width=0.7\textwidth]{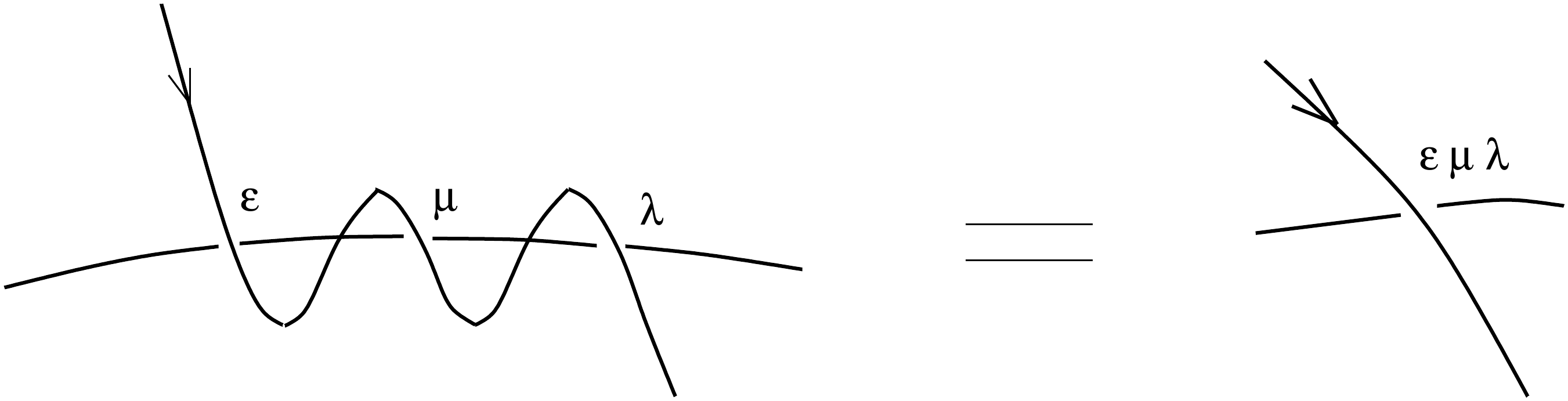}}

\vspace{.5cm}

The equality sign means that we can replace one tangle diagram by the other.

In particular, we get an interpretation for the crossing decorated by a scale 
parameter. Look first to this equality of tangle diagrams.

\vspace{.5cm}

\centerline{\includegraphics[angle=270, width=0.7\textwidth]{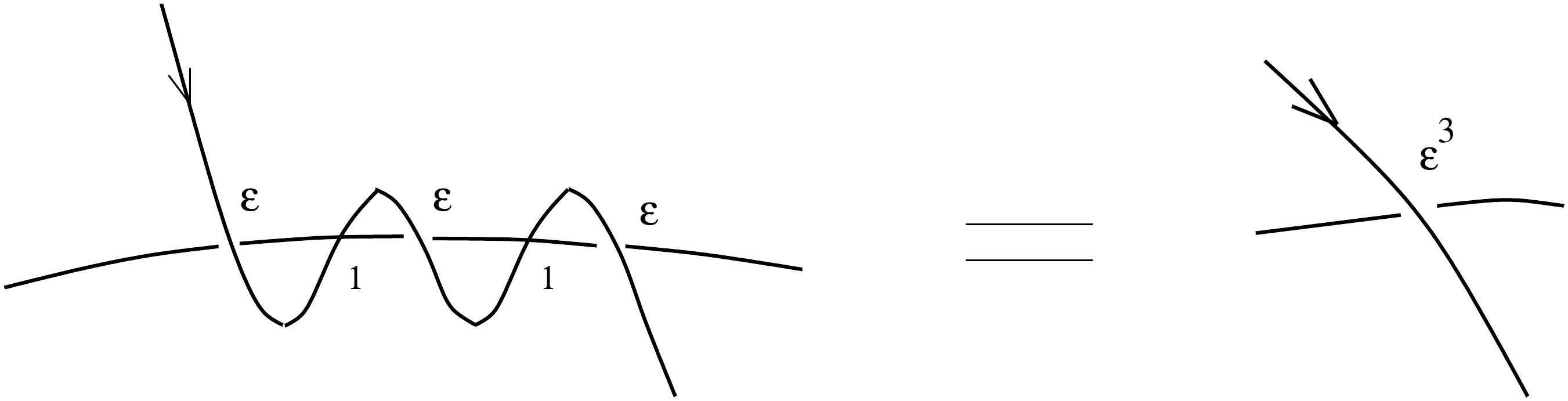}}

\vspace{.5cm}

If we fix the $\varepsilon$, take for example $\varepsilon = 1/2$, then 
any crossing decorated by a power of this $\varepsilon$ is equivalent with 
a chain of crossings decorated with $\varepsilon$, with virtual crossings
inserted in between. 

The usual interpretation of virtual crossings is that these are crossings which 
are not really there. Alternatively, but only to get an intuitive image, we 
may imagine that a crossing decorated with  $\displaystyle \varepsilon^{n}$ is 
equivalent with the projection of a helix arc with $n$ turns around an 
imaginary cylinder.
 
\vspace{.5cm}

\centerline{\includegraphics[angle=270, width=0.5\textwidth]{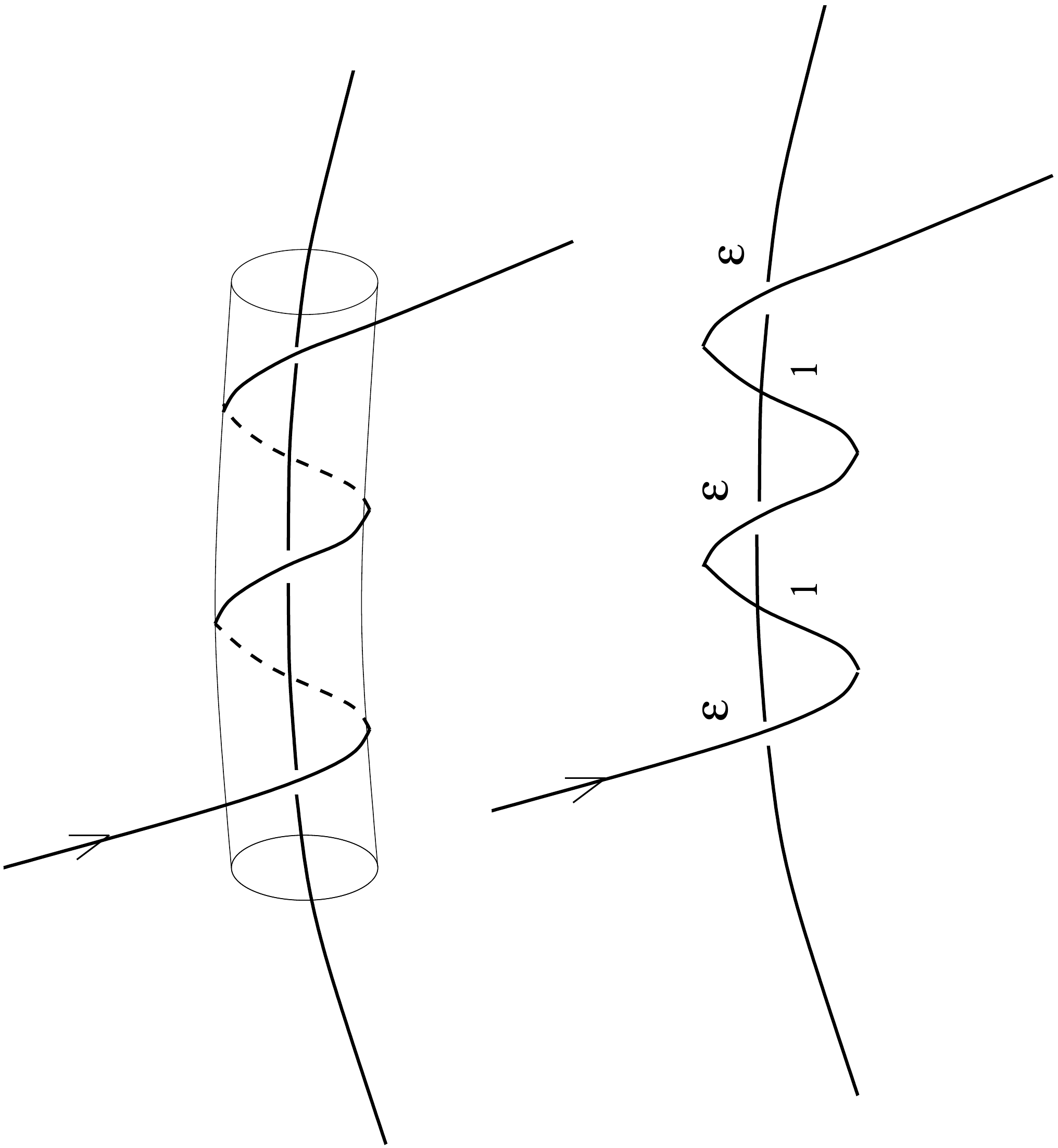}}

\vspace{.5cm}

Thus, for example if we take $\varepsilon = 1/2$ as a "basis", then a crossing 
decorated with the scale parameter $\mu$ could be imagined as the projection 
of a helix arc with $\displaystyle - \log_{2} \mu$ turns.

The sequence of irqs is the same as the algebraic object called a 
$\Gamma$-irq. The definition is given further (remember that for the needs of 
this paper $\Gamma = (0,+\infty)$).

\begin{definition}
Let $\Gamma$ be a commutative group. A $\Gamma$-idempotent right quasigroup 
is a set $X$ with a function $\displaystyle \varepsilon \in \Gamma \mapsto 
\circ_{\varepsilon}$ such that for any $\varepsilon \in \Gamma$ the pair $\displaystyle (X, \circ_{\varepsilon})$ is a irq
and moreover for any $\varepsilon, \mu \in \Gamma$ and any $x, y \in X$ we have 
$$x \, \circ_{\varepsilon} \, \left( x \, \circ_{\mu} \, y \right) \, = \, 
x \, \circ_{\varepsilon \mu} \, y$$
\label{defgammairq}
\end{definition}

\paragraph{Rules concerning wires.}

(W1) We may join two wires decorated by the same element of the $\Gamma$-irq and with the same
orientation.

\vspace{.5cm}

\centerline{\includegraphics[angle=270, width=0.6\textwidth]{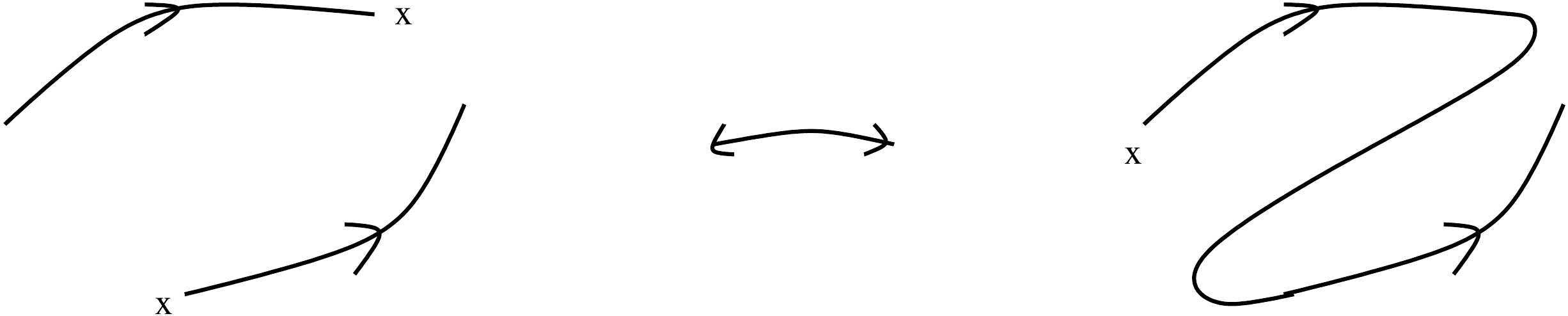}}
%\caption{}
%\label{join}

\vspace{.5cm}

(W2) We may change the orientation in a wire which passes over others, but
we must invert (power "-1") the decoration of each crossing. 
\vspace{.5cm}

\centerline{\includegraphics[angle=270, width=0.5\textwidth]{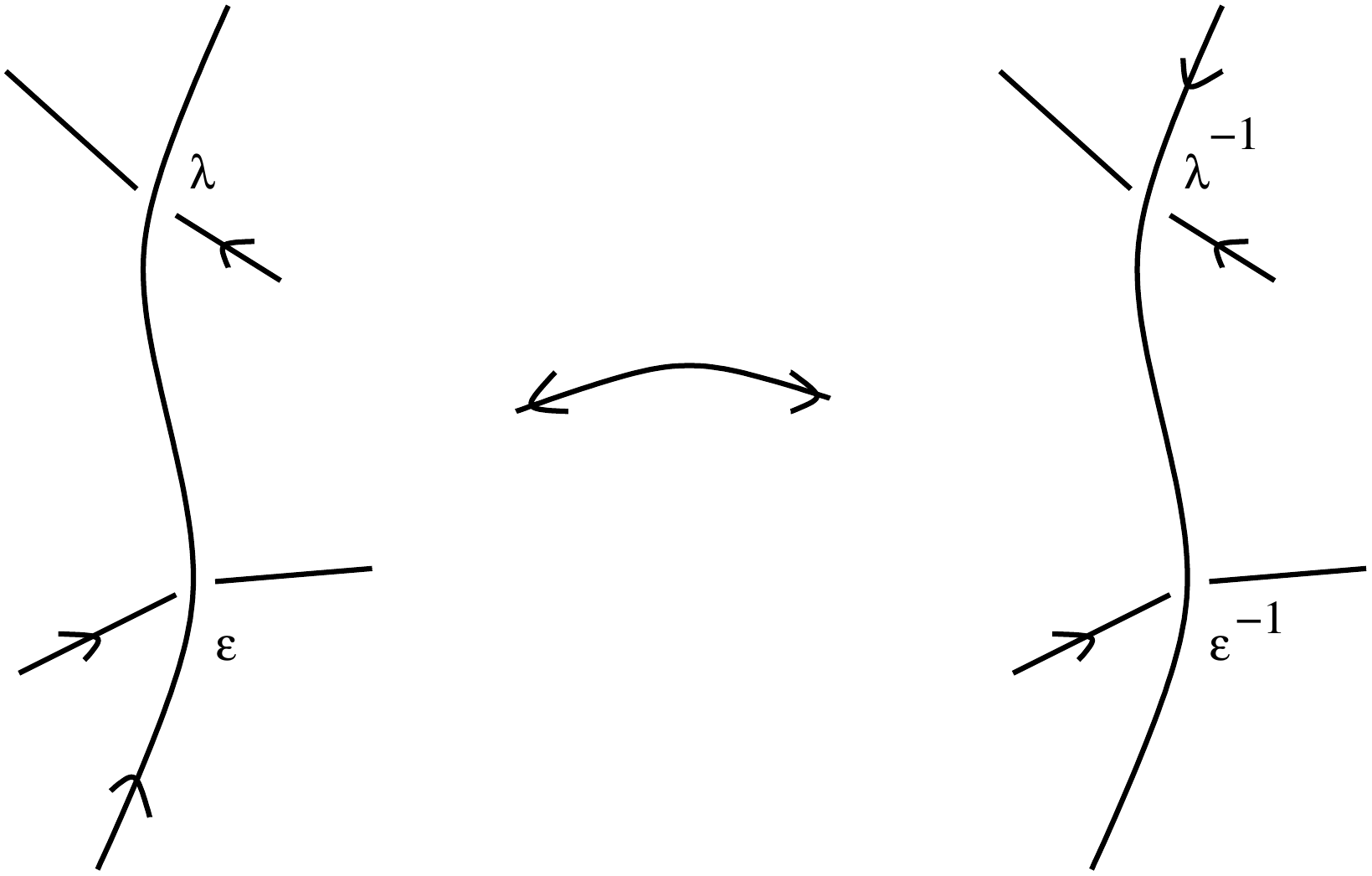}}
%\caption{}
%\label{changeorientation}

\vspace{.5cm}

\subsection{Decorated binary trees}

Here I use the second image of a oriented tangle diagram in order to 
understand the rules of decoration and movements described in the previous
section.  

In this interpretation an oriented tangle diagram is  an  
oriented planar graph with 3-valent and 1-valent nodes (input or exit nodes), 
 connected by wires which could cross  (crossings of wires in this graph has
  no meaning). 

In fact we consider trivalent oriented  planar graphs (together with 1-valent nodes representing
inputs and outputs), with wires (and input, output nodes) decorated by elements
of a $\Gamma$-irq $\displaystyle \varepsilon \in \Gamma \mapsto 
(X, \circ_{\varepsilon}, \bullet_{\varepsilon})$. The nodes are 
 either undecorated (corresponding to FAN-OUT gates) or
decorated by pairs $(\circ, \varepsilon)$ or $(\bullet, \varepsilon)$, with 
$\varepsilon \in \Gamma$. The rules of decoration are the following. 

\vspace{.5cm}

\centerline{\includegraphics[angle=270, width=0.7\textwidth]{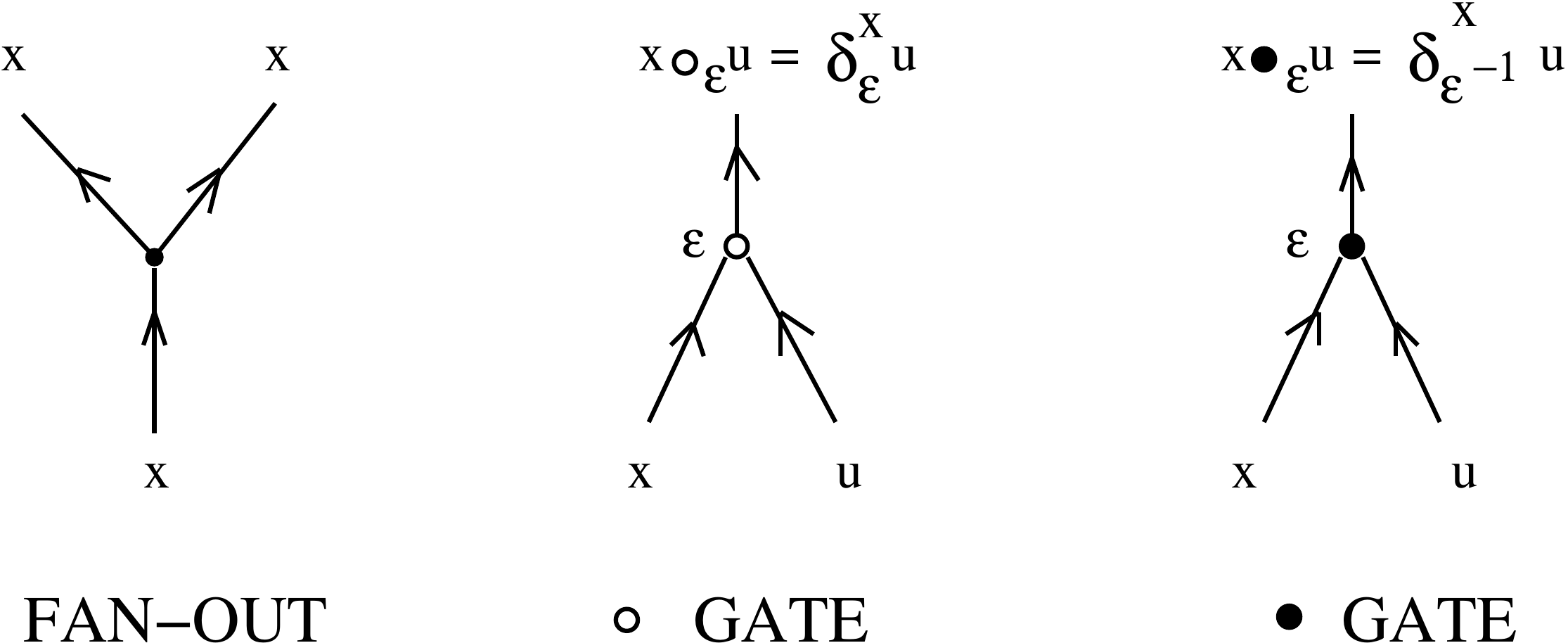}}
%\caption{}
%\label{join}

\vspace{.5cm}

The  trivalent graph is obtained from the tangle diagram by the procedure 
of replacing crossings with pairs of gates consisting of one FAN-OUT and 
one of the $\circ$ or $\bullet$ gates, explained in the next figure.

\vspace{.5cm}

\centerline{\includegraphics[angle=270, width=0.7\textwidth]{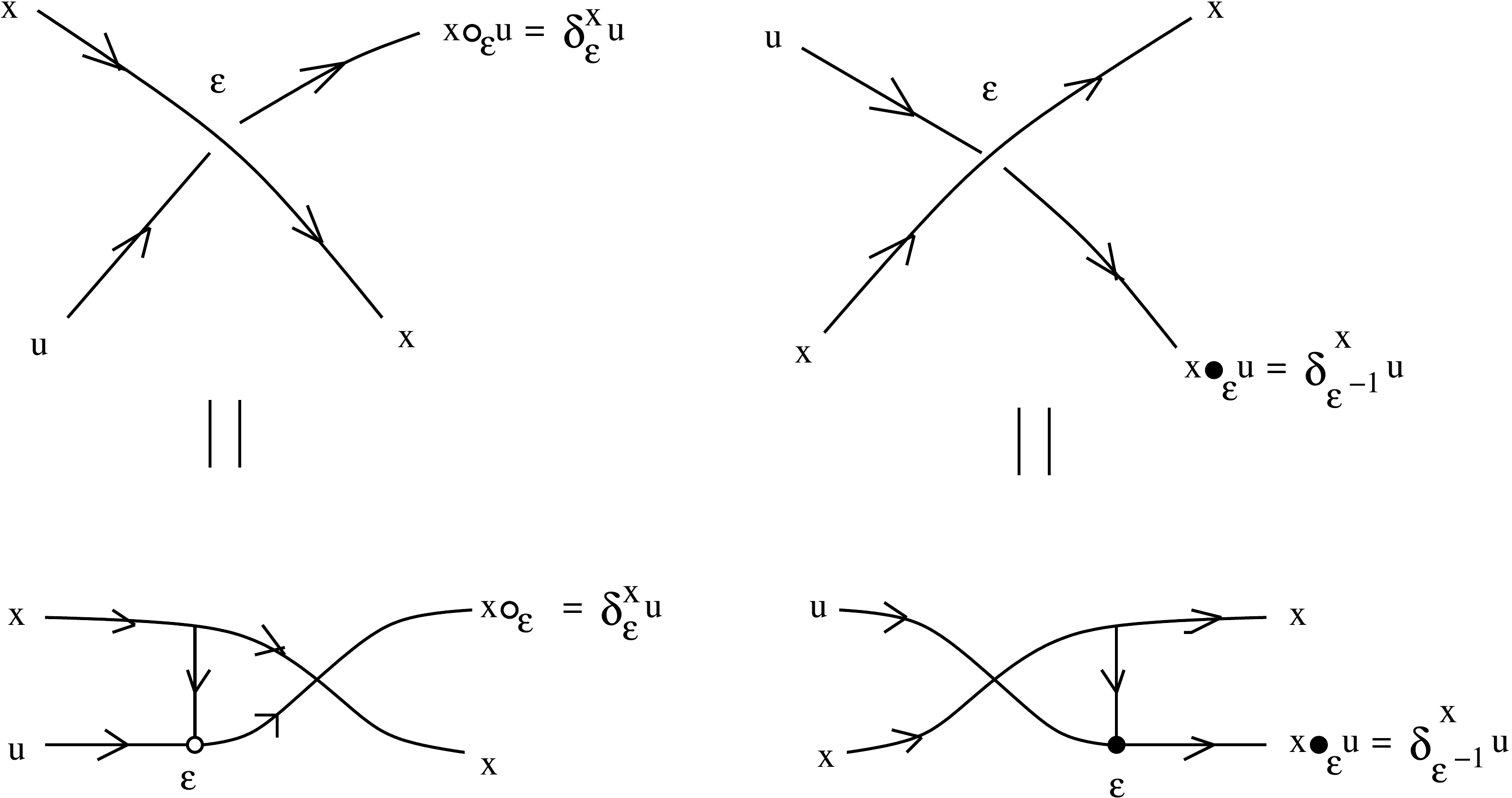}}
%\caption{}
%\label{join}

\vspace{.5cm}

In terms of trivalent graphs, the condition (R1) from definition 
\ref{defquasigroup} applied for the irq 
$\displaystyle (X,\circ_{\varepsilon}, \bullet_{\varepsilon})$ is graphically 
translated into the following identity (passing from one term of the identity
to another is a "Reidemeister I move"). 

\vspace{.5cm}

\centerline{\includegraphics[angle=270, width=0.4\textwidth]{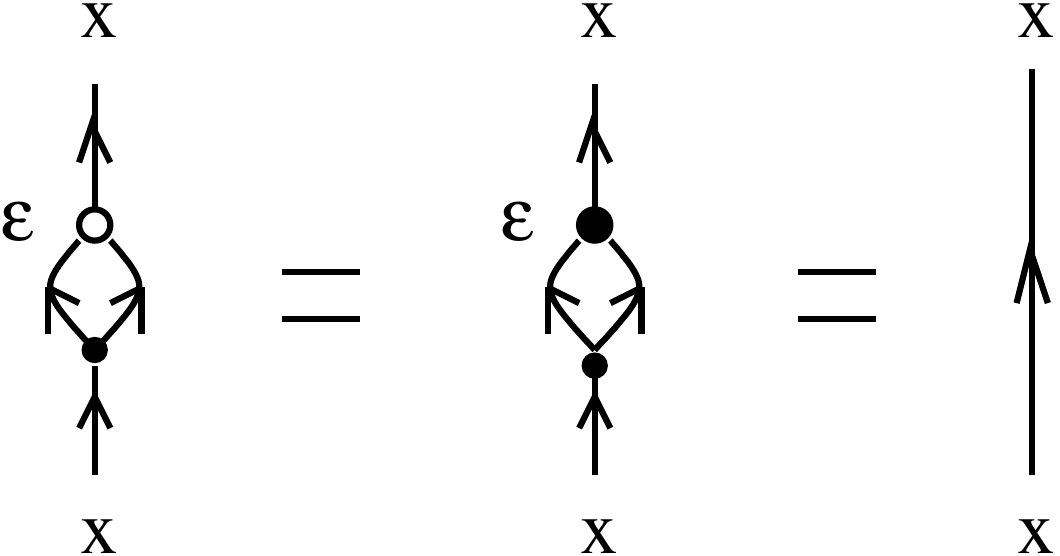}}
%\caption{}
%\label{join}

\vspace{.5cm}

There are two more groups of identities (or moves), which describe the
mechanisms of coloring trivalent graphs  $\Gamma$-irqs. The first group
consists in "triangle moves". This corresponds to Reidemeister move II and 
to the condition from the end of definition \ref{defgammairq}.

\vspace{.5cm}

\centerline{\includegraphics[angle=270, width=0.7\textwidth]{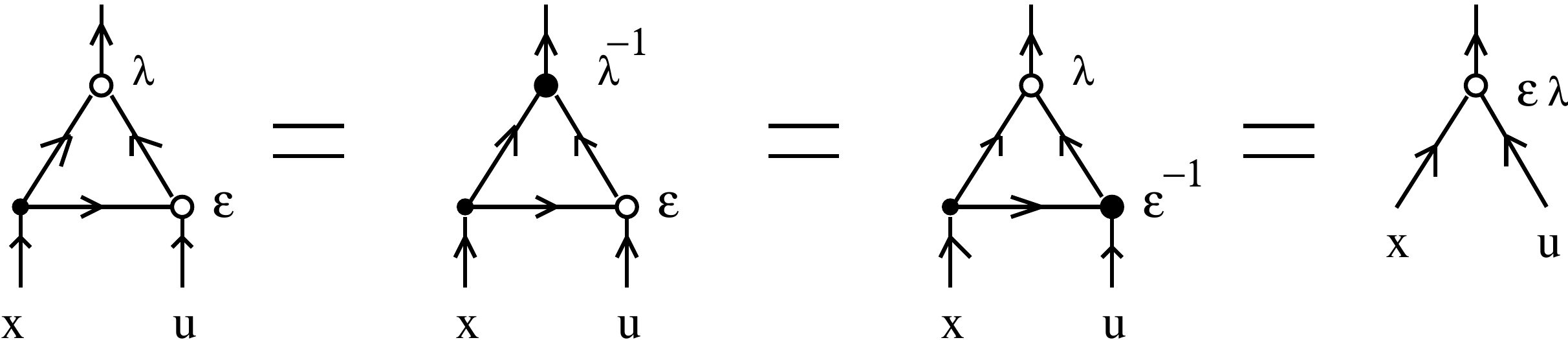}}
%\caption{}
%\label{join}

\vspace{.5cm}

The second groups is equivalent with the re-wiring move (W1) and the relation 
$\displaystyle x \circ_{1} u = x \bullet_{1} u = u$.

\vspace{.5cm}

\centerline{\includegraphics[angle=270, width=0.6\textwidth]{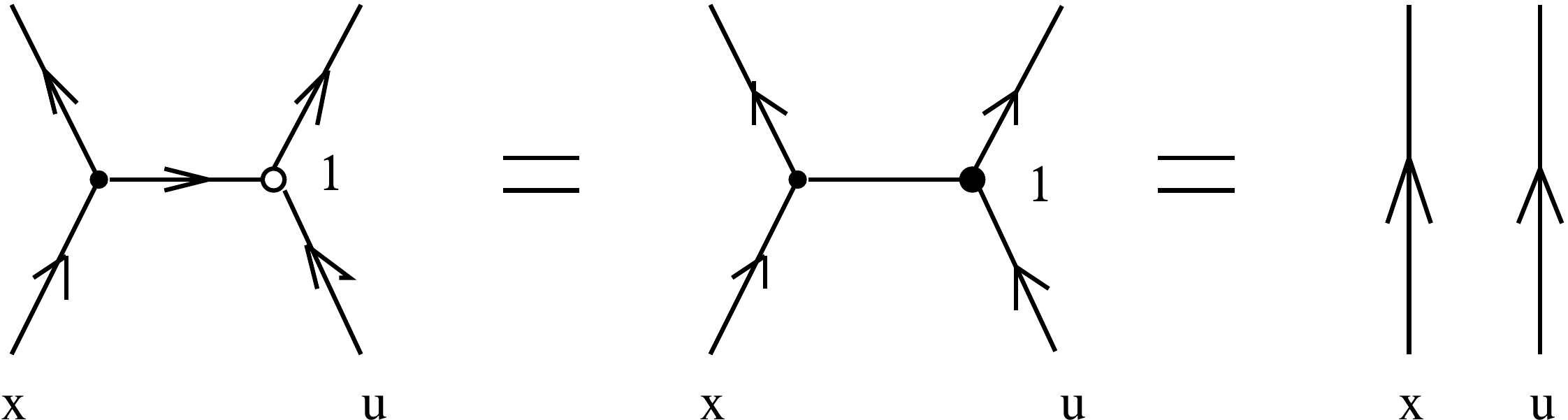}}
%\caption{}
%\label{join}

\vspace{.5cm}

\subsection{Linearity, self-similarity, Reidemeister III move}

Let $f: Y \rightarrow X$ be an invertible function. We can use the tangle
formalism for picturing the function $f$. To $f$ is associated a special curved
segment, figured by a double line. The crossings passing under this double line 
are colored following the rules explained in this figure.

\vspace{.5cm}

\centerline{\includegraphics[angle=270, width=0.6\textwidth]{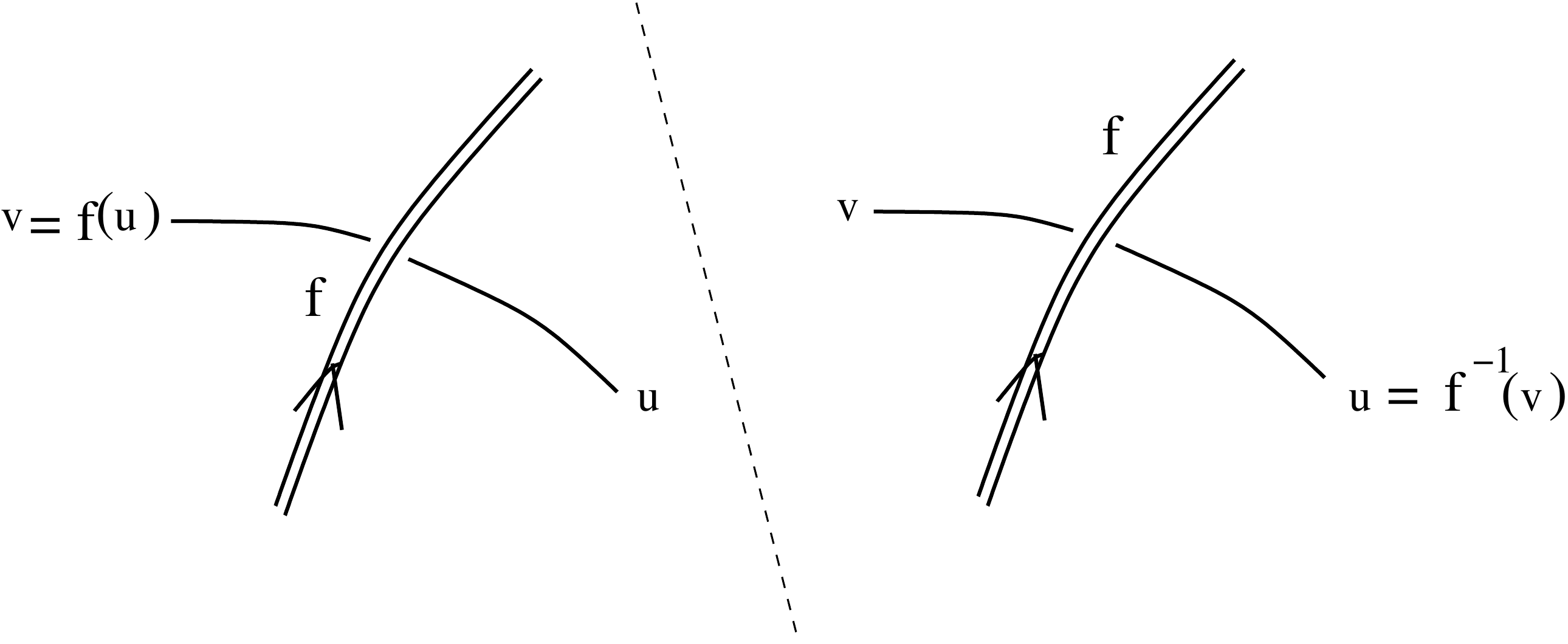}}
%\caption{}
%\label{reidemeistermoves1and2}

\vspace{.5cm}

Suppose that $X$ and $Y$ are endowed with a $\Gamma$-irq structure (in
particular, we may suppose that they are endowed with dilation structures).
Consider the following sliding movement. 

\vspace{.5cm}

\centerline{\includegraphics[angle=270, width=0.6\textwidth]{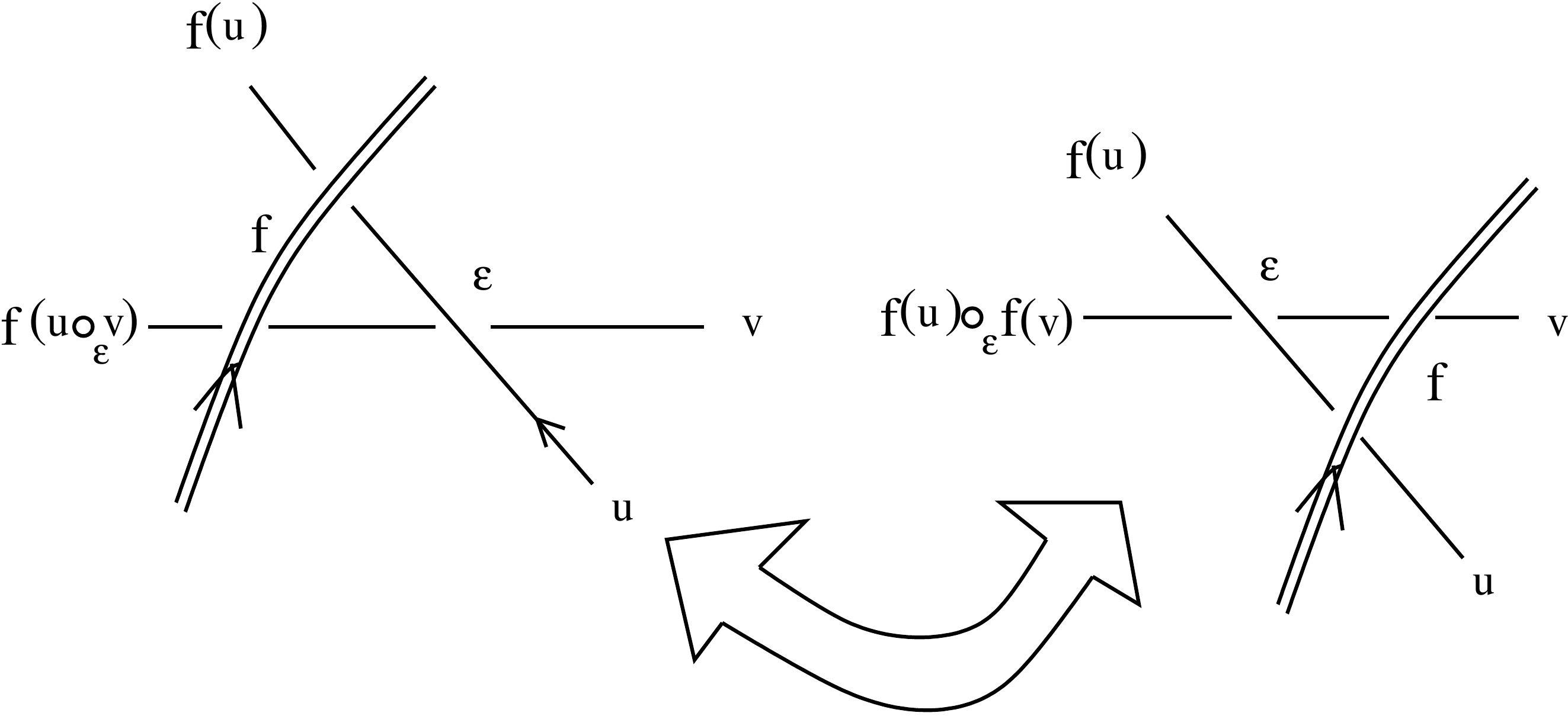}}
%\caption{}
%\label{reidemeistermoves1and2}

\vspace{.5cm}

The crossing decorated with $\varepsilon$ from the left hand side diagram is 
in $X$ (as well as the rules of decoration with the $\Gamma$-irq of $X$).
Similarly, the crossing decorated with $\varepsilon$ from the right hand 
side diagram is in $Y$. Therefore these two diagrams are equal (or we may pass 
from one to another by a sliding movement) if and only if $f$ transforms an 
operation into another, equivalently if $f$ is a morphism of $\Gamma$-irqs. 

This sliding movement becomes the Reidemeister III move in the case of 
$X = Y$ and $f$ equal to a dilation of $X$, 
$\displaystyle f = \delta^{x}_{\mu}$.

\begin{definition}
A function $f: X \rightarrow Y$ is linear if and only if it is a morphism 
of $\Gamma$-irqs (of $X$ and $Y$ respectively). Moreover, if $X$ and $Y$ are 
endowed with dilation structures then $f$ is linear if it is a morphism, written
in terms of dilations notation as: for any $u,v \in X$ and any $\varepsilon \in
\Gamma$  
$$f(\delta^{u}_{\varepsilon} v) \, = \, \delta^{f(u)}_{\varepsilon} f(v)$$ 
which is also a Lipschitz map from $X$ to $Y$ as metric spaces. 

A dilation structure $(X,d,\delta)$ is $(x,\mu)$ self-similar (for a 
$x \in X$ and $\mu\in \Gamma$, different from $1$, the neutral element
of $\Gamma = (0,+\infty)$) if the dilation $\displaystyle f = \delta^{x}_{\mu}$ is 
linear from $(X,d,\delta)$ to itself and moreover for any $u,v \in X$ we have 
$$d(\delta^{x}_{\mu} u, \delta^{x}_{\mu} v) \, = \, \mu \, d(u,v)$$

A dilation structure is linear if it is self-similar with respect to any $x \in
X$ and $\mu \in \Gamma$. 
\label{dline}
\end{definition}

Thus the Reidemeister III move is compatible with the tangle coloring by a 
dilation structure if and only if the dilation structure is linear. 

Conical groups are groups endowed with a one-parameter family of 
dilation morphisms. From the viewpoint of $\Gamma$-irqs, they are equivalent
with linear dilation structures (theorem 6.1 \cite{buligairq}, see also the
Appendix). 

A real vector space is a particular case of Carnot group. It is a commutative 
(hence nilpotent)  group with the addition of vectors operation and it has 
a one-parameter family of dilations defined by the multiplication of vectors 
by positive scalars. 

Carnot groups which are not commutative provide therefore a generalization of a 
vector space. Noncommutative Carnot groups are aplenty, in particular the 
simplest noncommutative Carnot groups are the Heisenberg groups, that is 
the simply connected Lie groups with the Lie algebra defined by the 
Heisenberg noncommutativity relations. 

For me Carnot groups, or conical groups, are just linear objects. (By extension,
manifolds, which are assemblies of open subsets of vector spaces, are locally 
linear objects as well. Moreover, they are "commutative", because the model of
the tangent space at a point is a commutative Carnot group.)

It is also easy to explain graphically the transport of a dilation structure, 
or of a $\Gamma$-irq from $X$ to $Y$, by using $\displaystyle f^{-1}$. 

\vspace{.5cm}

\centerline{\includegraphics[angle=270, width=0.6\textwidth]{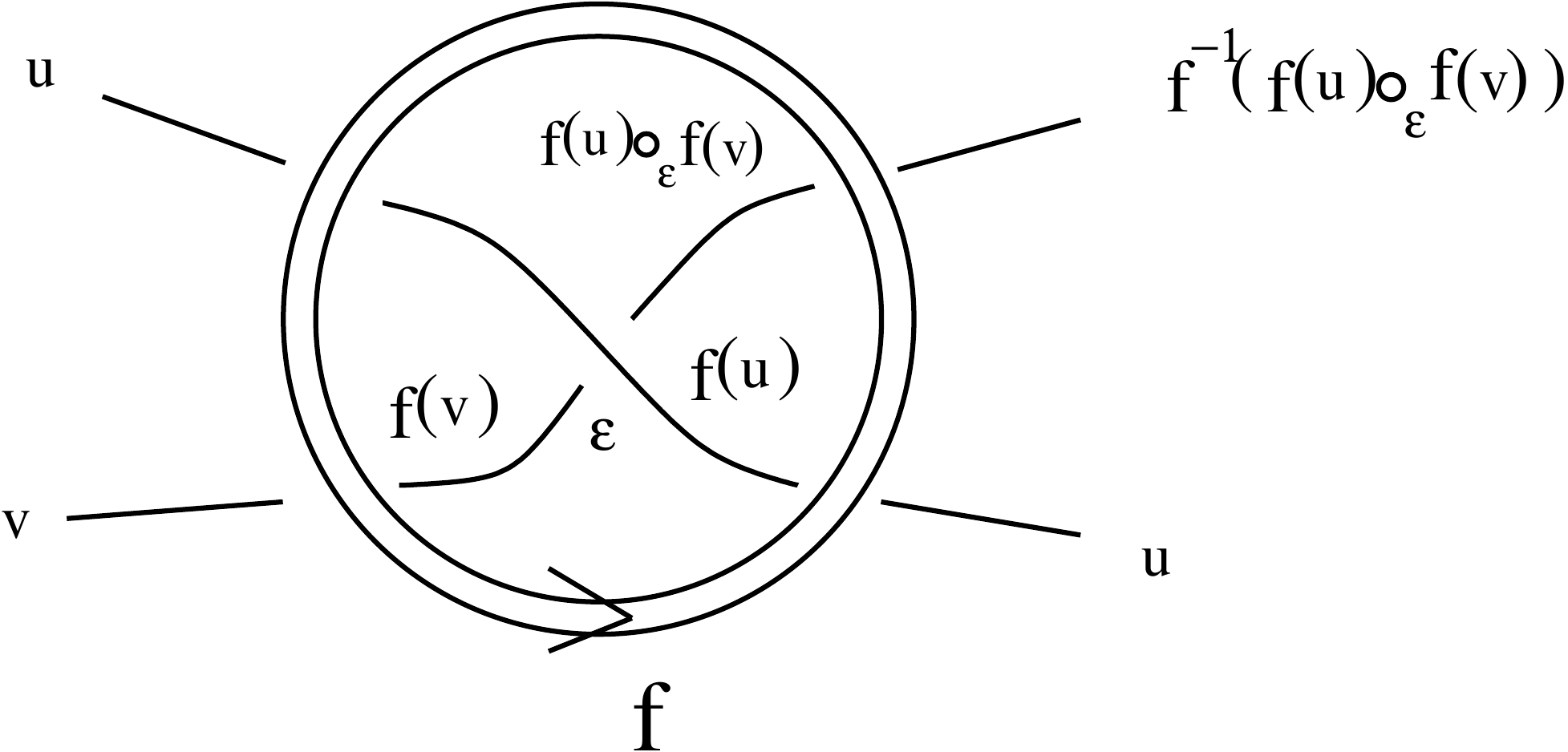}}
%\caption{}
%\label{reidemeistermoves1and2}

\vspace{.5cm}

The transport operation amounts to adding a double circle decorated by $f$,
which overcrosses the whole diagram (in this case a diagram containing only one 
crossing). If we use the map-territory distinction, then inside the circle we
are in $Y$, outside we are in $X$. 

It is obvious that $f$ is linear if and only if the transported dilation
structure on $X$ (by $\displaystyle f^{-1}$) coincides with the dilation
structure on $X$. Shortly said, encirclings by linear functions can be removed 
from the diagram. 

\vspace{.5cm}

\centerline{\includegraphics[angle=270, width=0.6\textwidth]{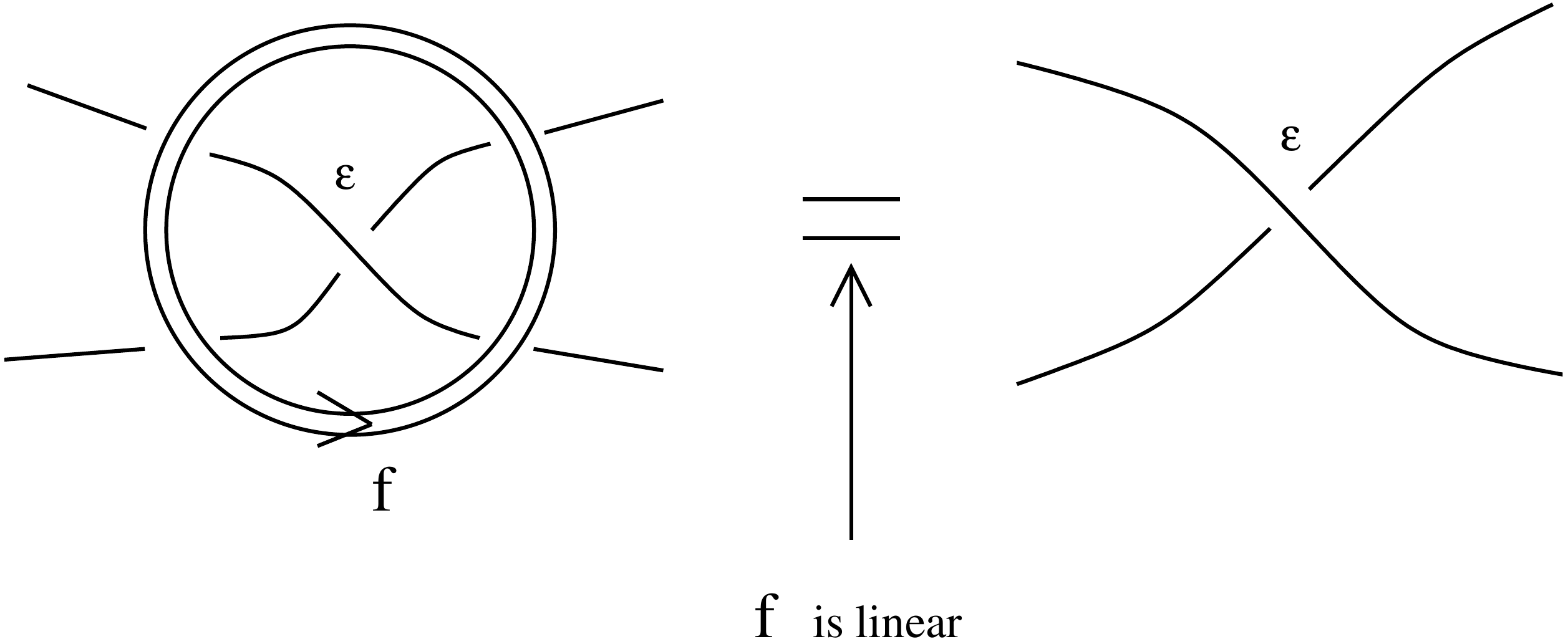}}
%\caption{}
%\label{reidemeistermoves1and2}

\vspace{.5cm}

 In this tangle decoration formalism we have no 
reason to suppose that the dilations structures which we use are linear. This
would be an unnecessary limitation of the dilation structure 
(or emergent algebra) formalism. That is why the Reidemeister III move is not
an acceptable move in this formalism.

 The last axiom (A4) of dilation structures can be translated into 
an algebraic statement which will imply a weak form of the Reidemeister 
III move, namely that this move can be done IN THE LIMIT. 

\subsection{Acceptable tangle diagrams}

Consider  a tangle diagram  with decorated crossings, but with undecorated 
segments.  

A  notation for such a diagram is $T[\varepsilon,\mu,\eta,...]$, where 
$\varepsilon, \mu, \eta, ...$ are decorations of the 3-valent nodes (in the 
second image) or decorations of the crossings (in the first image of a tangle
diagram).

A tangle diagram $T[\varepsilon,\mu,\eta,...]$ is "acceptable" if there exists at
least a decoration of the input segments such that all the segments of the 
diagram can be decorated according to the rules specified previously, maybe
non-uniquely.

A set of parameters of an acceptable  tangle diagram is any coloring of a 
part of the segments of the diagram such that  any coloring of the 
input segments which are already not colored by parameters, can be completed in a way 
which is unique for the output segments (which are not already colored by
parameters).   

Given an acceptable tangle diagram which admits a set of parameters, we shall
see it as a function from the colorings of the input to the colorings of the
output, with parameters from the set of parameters and with "scale parameters" 
the decorations of the crossings. 

Given one acceptable tangle diagram which admits a set of parameters, given a 
set of parameters for it, we can choose one or more crossings and 
their decorations (scale parameters) as "scale variables". 
This is equivalent to considering a  sequence of acceptable tangle diagrams, 
indexed by a multi-index of scale variables (i.e. taking values in some 
cartesian power of $\Gamma$). Each member of the sequence has the same 
tangle diagram, with the same set of parameters, with the same decorations 
of crossings which are not variables; the only difference is in the decoration 
of the crossings chosen as variables.  

Associated to such a sequence  is the sequence of input-output functions 
of these diagrams. We shall consider uniform convergence of these functions 
with respect to compact sets of inputs. 

All this is needed to formulate the emergent algebra correspondent of axiom A4 
of dilation structures.

\subsection{Going to the limit: emergent algebras}

Basically, I see a decorated tangle diagram as an expression dependent on the 
decorations of the crossings. More precisely, I shall reserve the letter 
$\varepsilon$ for an element of $\Gamma$ which will be conceived as going 
to zero. This is the same kind of reasoning as for the zoom sequences in the
section dedicated to maps.

Why is such a thing interesting? Let me give some examples. 

\paragraph{Finite differences.} We use the convention of adding to the tangle
diagram supplementary arcs decorated by homeomorphisms. Let $f: X \rightarrow Y$
such a homeomorphism. I want to be able to differentiate the homeomorphism $f$, 
in the sense of dilation structures. 

For this I need a notion of finite differences. These appear as the following diagram. 

\vspace{.5cm}

\centerline{\includegraphics[angle=270, width=0.6\textwidth]{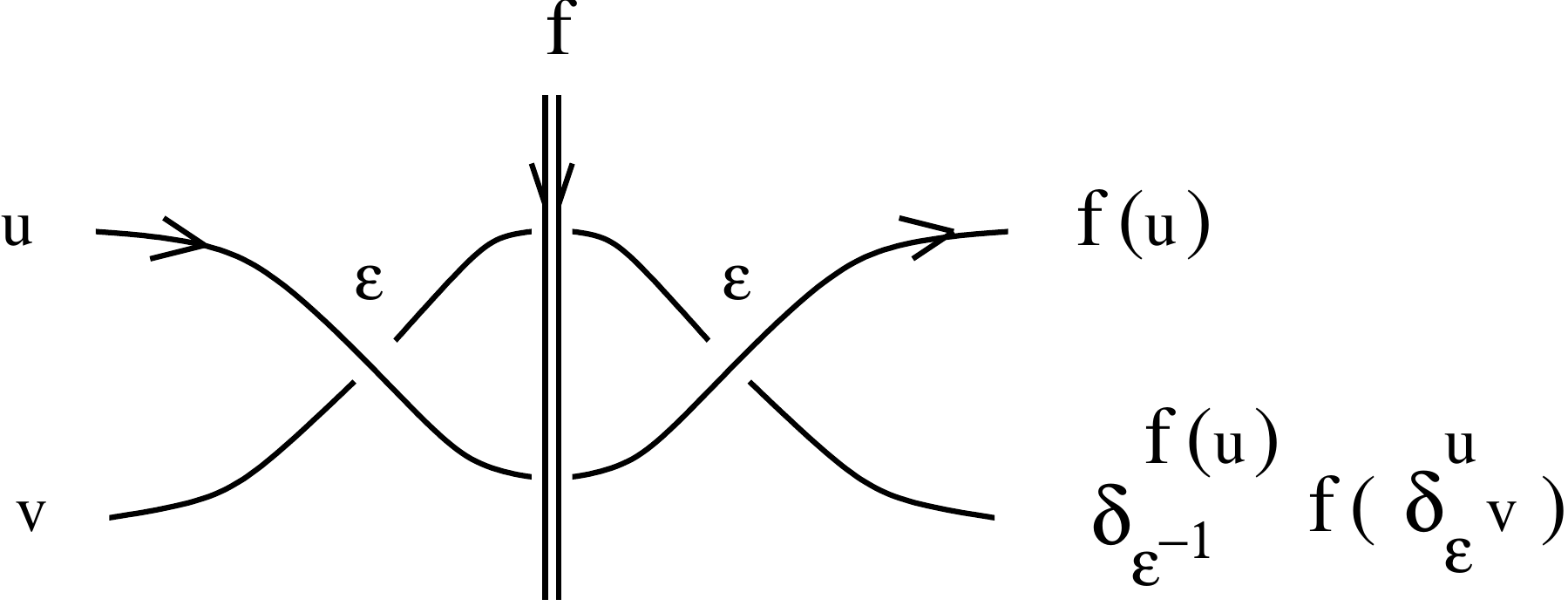}}
%\caption{}
%\label{reidemeistermoves1and2}

\vspace{.5cm}

Indeed, suppose for simplicity that $X$ and $Y$ are finite dimensional normed
vector spaces, with distance given by the norm and dilations 
$$\delta^{u}_{\varepsilon} u \, = \, u + \varepsilon (-u+v)$$
Then we have: 
$$\delta^{f(u)}_{\varepsilon^{-1}} f \left( \delta^{u}_{\varepsilon} v \right)
\, = \, f(u) \, + \, \frac{1}{\varepsilon} \left( f( u + \varepsilon (-u+v)) -
f(u)\right)$$ 

Pansu \cite{pansu} generalized this definition of finite differences from 
real vector spaces to Carnot groups, which are nilpotent graduated simply
connected Lie groups, a particular example of conical groups.  

It is true that the diagram which encodes finite differences is not, technically
speaking, of the type explain previously, because it has a segment (the one 
decorated by $f$), which is different from the other segments. But it is easy to
see that all the mathematical formalism can be modified easily in order to
accommodate such edges decorated with homeomorphism. 

The notion of differentiability of $f$ is obtained by asking that the sequence
of input-output functions associated to the  "finite difference" diagram, 
 with parameter "$x$" and variable "$\varepsilon$", converges uniformly on 
 compact sets. See definition \ref{defdiffer}.

\paragraph{Difference gates.} 
For any $\varepsilon \in \Gamma$, the $\varepsilon$-difference gate, 
is described by the next figure.

\vspace{.5cm}

\centerline{\includegraphics[angle=270, width=0.6\textwidth]{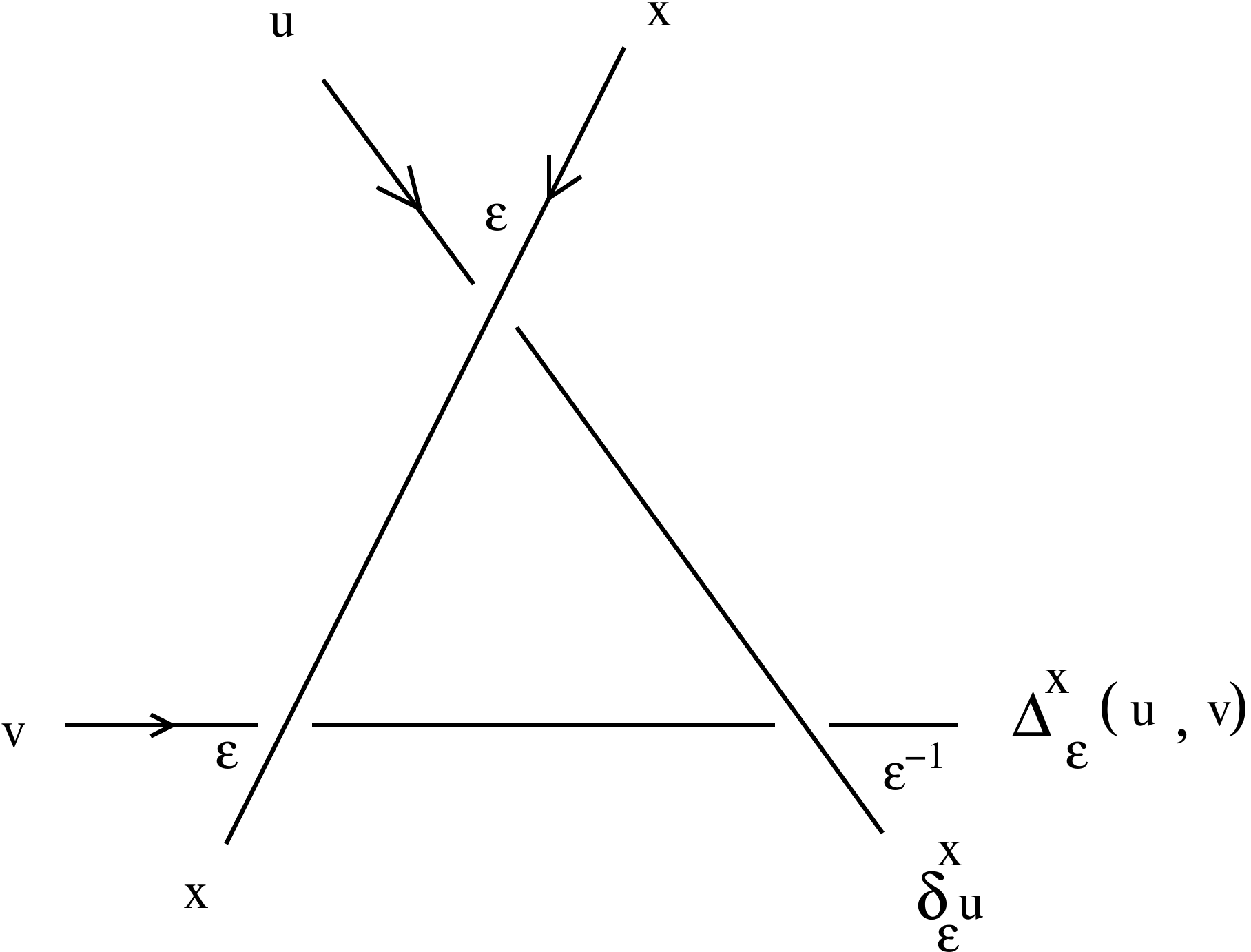}}
%\caption{}
%\label{reidemeistermoves1and2}

\vspace{.5cm}

 Here   $\displaystyle  \Delta^{x}_{\varepsilon} (u,v)$ is a construct made 
 from operations $\displaystyle \circ_{\varepsilon}$, 
 $\displaystyle \bullet_{\varepsilon}$. It corresponds to the difference coming
 from changing the viewpoint, in the map-territory frame. In terms of 
 dilation structures, is the approximate difference which appears in axiom A4. 
 In terms of notations of a $\Gamma$-irq,  from the figure we 
can compute $\displaystyle  \Delta^{x}_{\varepsilon} (u,v)$ as 
$$\displaystyle  \Delta^{x}_{\varepsilon} (u,v) \, = \, \left( x \circ_{\varepsilon} u \right) \, 
\bullet_{\varepsilon} \, \left( x \circ_{\varepsilon} v \right)$$

The geometric meaning of $\displaystyle \Delta^{x}_{\varepsilon}(u,v)$ is that 
it is indeed a kind of approximate difference between the vectors 
$\displaystyle \vec{xu}$ and $\vec{xv}$, by means of a generalization of 
the parallelogram law of vector addition.  This is  shown in the following 
figure, where straight lines have been replaced by slightly curved ones in 
order to suggest that this construction has meaning in settings far more 
general than euclidean spaces, like in Carnot-Caratheodory or sub-riemannian geometry, as shown in \cite{buligadil1}, or generalized (noncommutative) affine geometry 
 \cite{buligadil2}, for length metric spaces with dilations \cite{buligadil3} or even for normed groupoids 
 \cite{buligagr}.

\vspace{.5cm}

\centerline{\includegraphics[angle=270, width=0.6\textwidth]{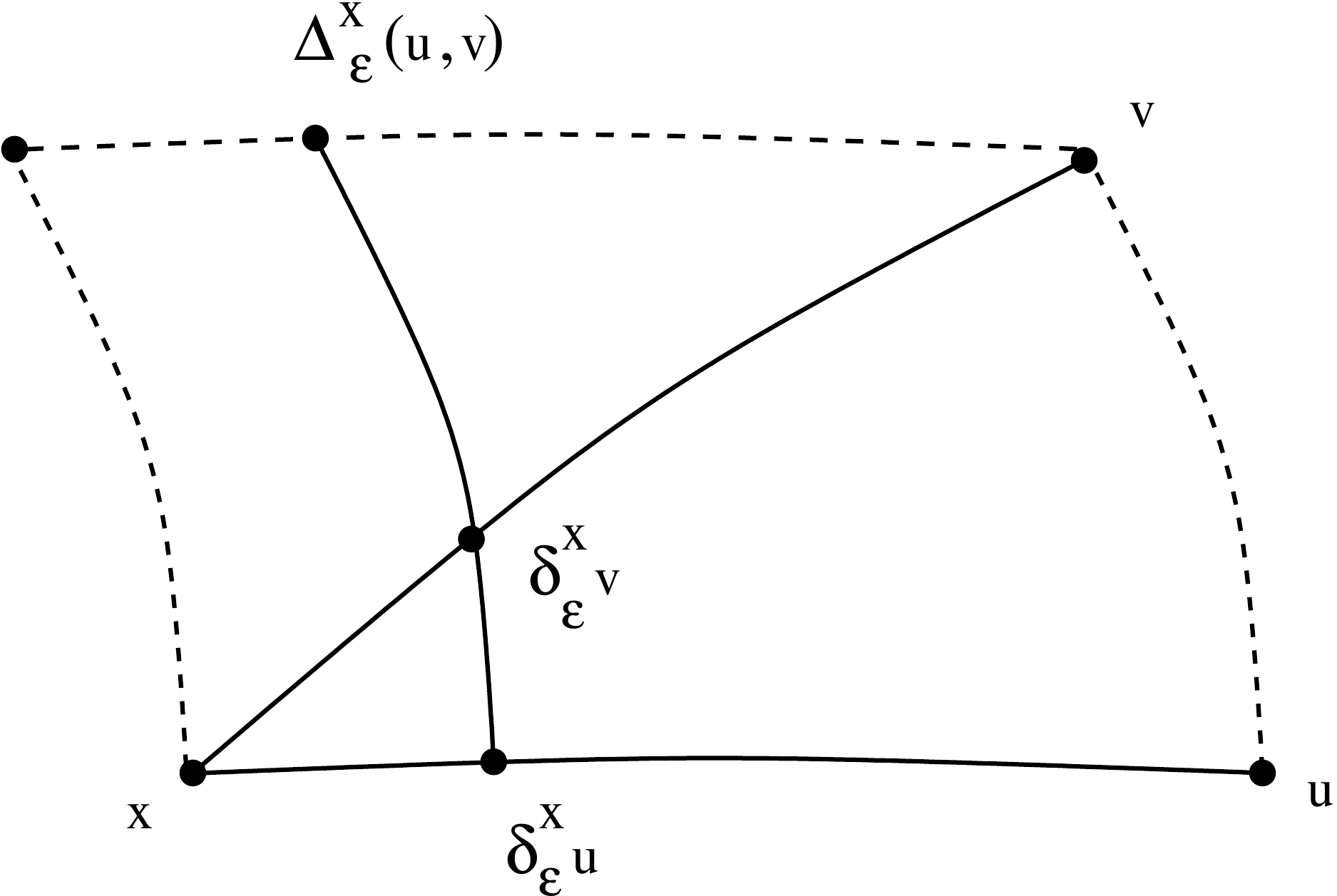}}
%\caption{}
%\label{reidemeistermoves1and2}

\vspace{.5cm}

The $\varepsilon$-sum gate is described in the next figure. 

\vspace{.5cm}

\centerline{\includegraphics[angle=270, width=0.6\textwidth]{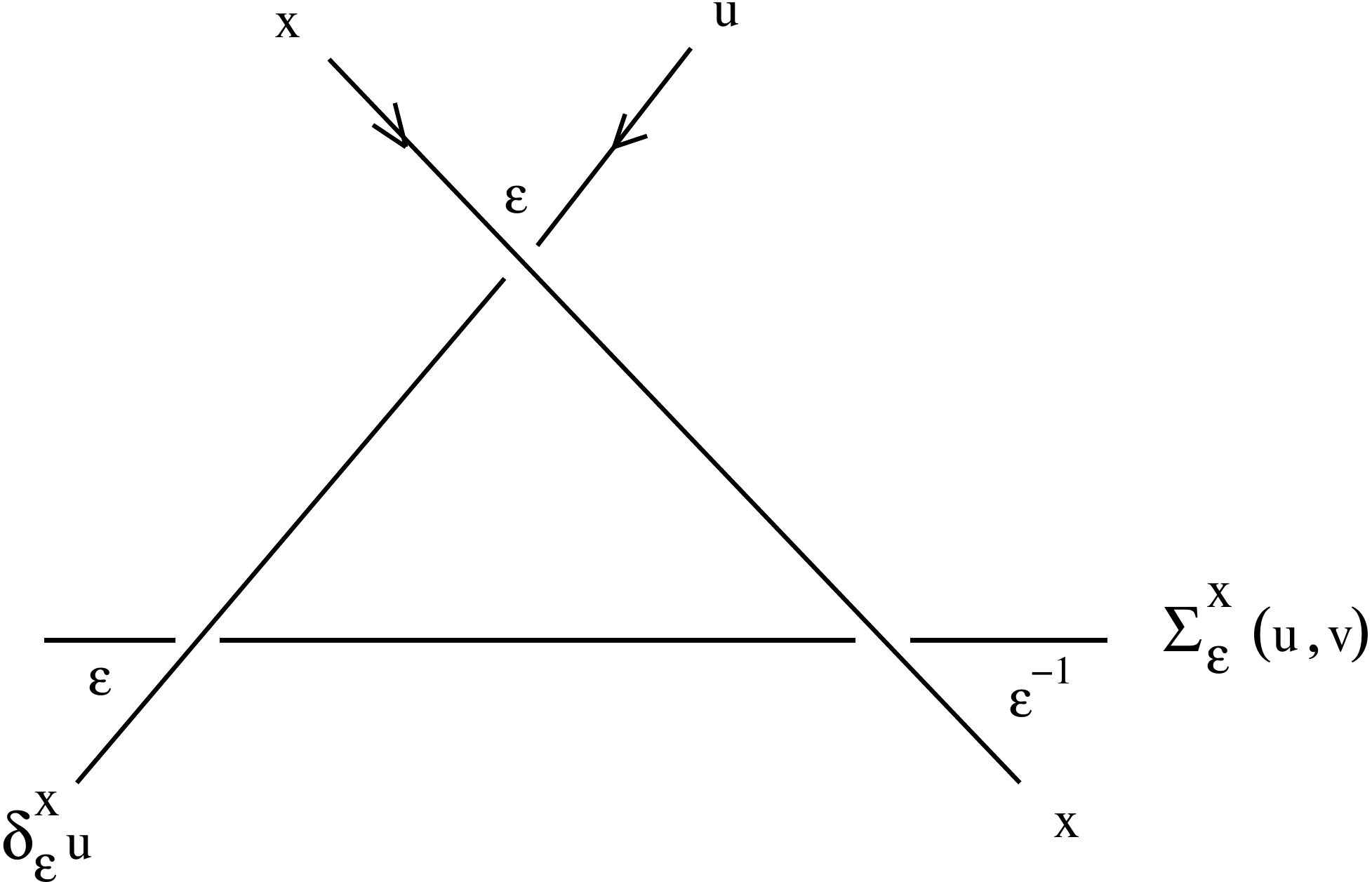}}
%\caption{}
%\label{reidemeistermoves1and2}

\vspace{.5cm}

Similar comments can be made, concerning the sum gate. It is the approximate sum
appearing in the axiom A4 of dilation structures. 

Finally, there is another important tangle diagram, called $\varepsilon$-inverse
gate. It is, at closer look, a particular case of a difference gate (take $x =
u$ in the difference diagram). 

\vspace{.5cm}

\centerline{\includegraphics[angle=270, width=0.6\textwidth]{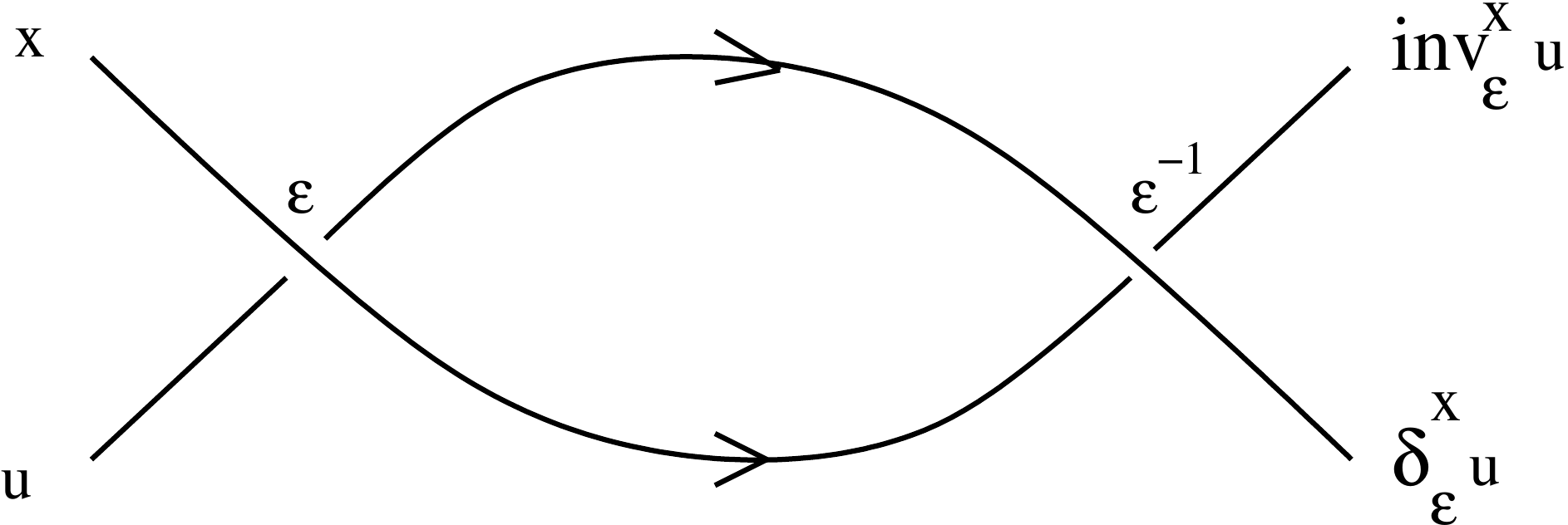}}
%\caption{}
%\label{reidemeistermoves1and2}

\vspace{.5cm}

The relevant outputs of the previously introduced gates, namely the approximate difference, sum and 
inverse functions, are described in the next definition, in terms of decorated
binary trees (trivalent graphs). I am going to ignore the trees constructed from
 FAN-OUT  gates, replacing them by patterns of decorations (of leaves of the
 binary trees). In the following all tree nodes are
 decorated with the same label $\varepsilon$ and edges are oriented upwards. 

\begin{definition}
We define the difference, sum and inverse trees  given by:

\centerline{\includegraphics[width=120mm]{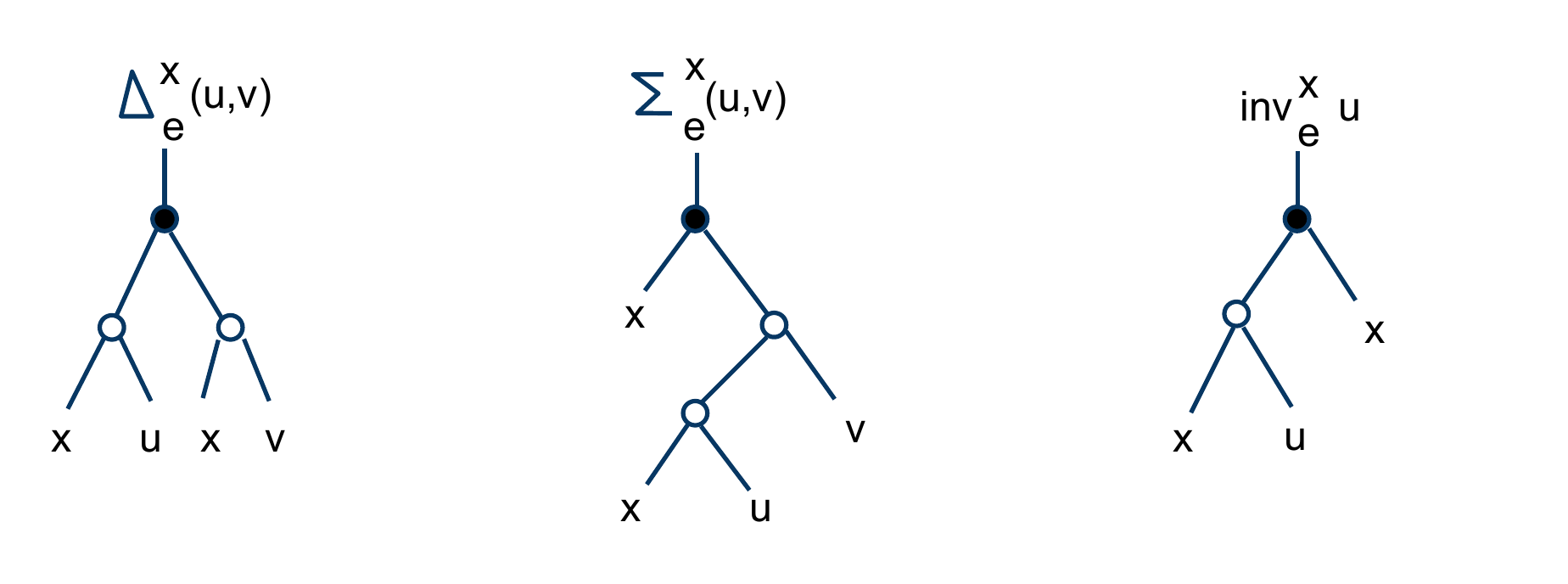}}

\label{defsumdif}
\end{definition}

The following proposition contains the main relations between the difference,
sum and inverse gates. They can all be proved by this tangle diagram formalism. 
In \cite{buligadil1} I explained these relations as appearing from the 
equivalent formalism using binary decorated trees.

\begin{proposition}
Let $\displaystyle (X,\circ_{\varepsilon})_{\varepsilon \in \Gamma}$ be
a $\Gamma$-irq. Then we have the relations: 
\begin{enumerate}
\item[(a)] $\displaystyle \Delta^{x}_{\varepsilon}(u,
\Sigma^{x}_{\varepsilon}(u,v)) \, = \, v$  (difference is the inverse of sum) 
\item[(b)] $\displaystyle \Sigma^{x}_{\varepsilon}(u,
\Delta^{x_{\varepsilon}}(u,v)) \, = \, v$  (sum is the inverse of difference) 
\item[(c)]  $\displaystyle \Delta^{x}_{\varepsilon}(u, v) \, = \, \Sigma^{x \circ_{\varepsilon} u}_{\varepsilon} 
(inv_{\varepsilon}^{x} u , v)$ (difference  approximately equals  the sum of the inverse) 
\item[(d)]  $\displaystyle inv_{\varepsilon}^{x\circ u} \, inv_{\varepsilon}^{x}
\, u  \, = \, u $ (inverse operation is approximatively an involution) 
\item[(e)] $\displaystyle \Sigma^{x}_{\varepsilon}(u, \Sigma^{x\circ_{\varepsilon}u}_{\varepsilon}
(v , w)) \, = \, \Sigma^{x}_{\varepsilon}(\Sigma^{x}_{\varepsilon}(u,v), w) $ (approximate associativity of the sum)  
\item[(f)] $\displaystyle  inv^{x}_{\varepsilon} \, u \, = \,  \Delta^{x}_{\varepsilon}( u , x)$
\item[(g)]  $\displaystyle  \Sigma^{x}_{\varepsilon} (x, u) \, = \,  u $
(neutral element at right).
\end{enumerate}
\label{pplay}
\end{proposition}

We shall use the tree formalism to prove  some 
of these relations.   For complete proofs see \cite{buligadil1}. 

For example, in order to prove  (b) we do the following calculus:

\centerline{\includegraphics[width=120mm]{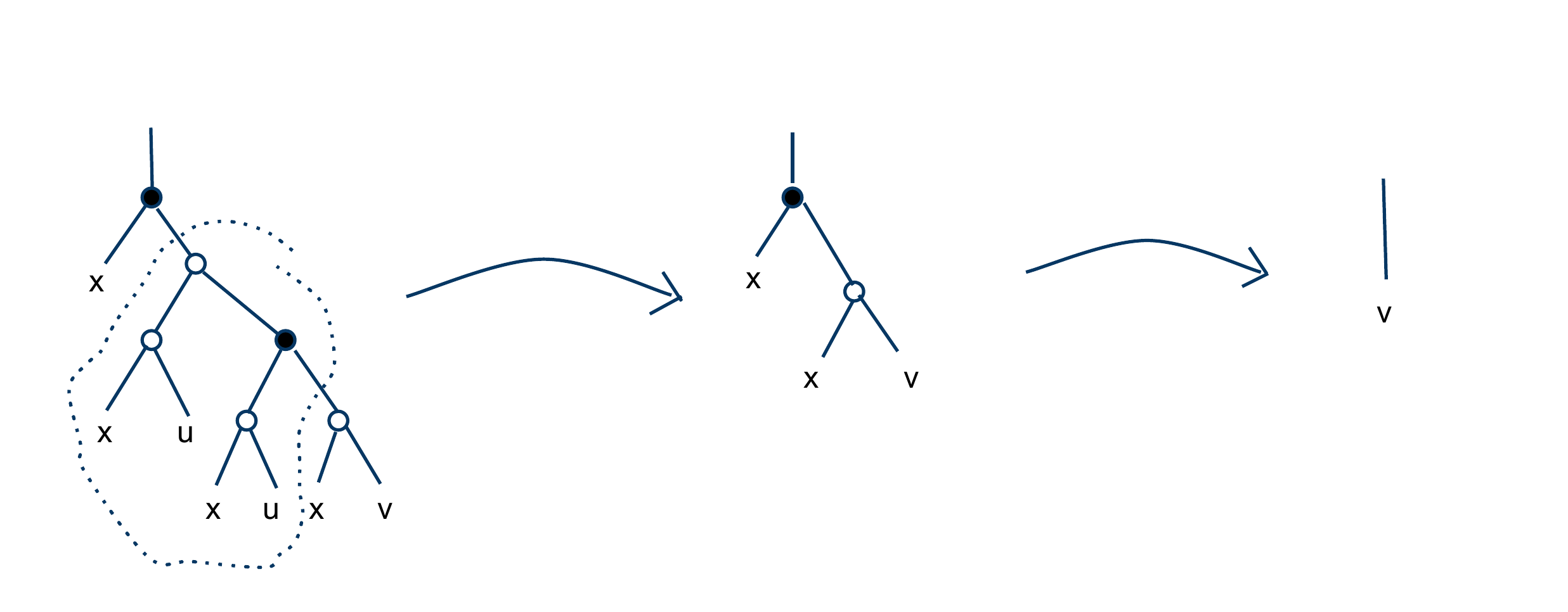}}

The relation (c) is obtained from: 

\centerline{\includegraphics[width=100mm]{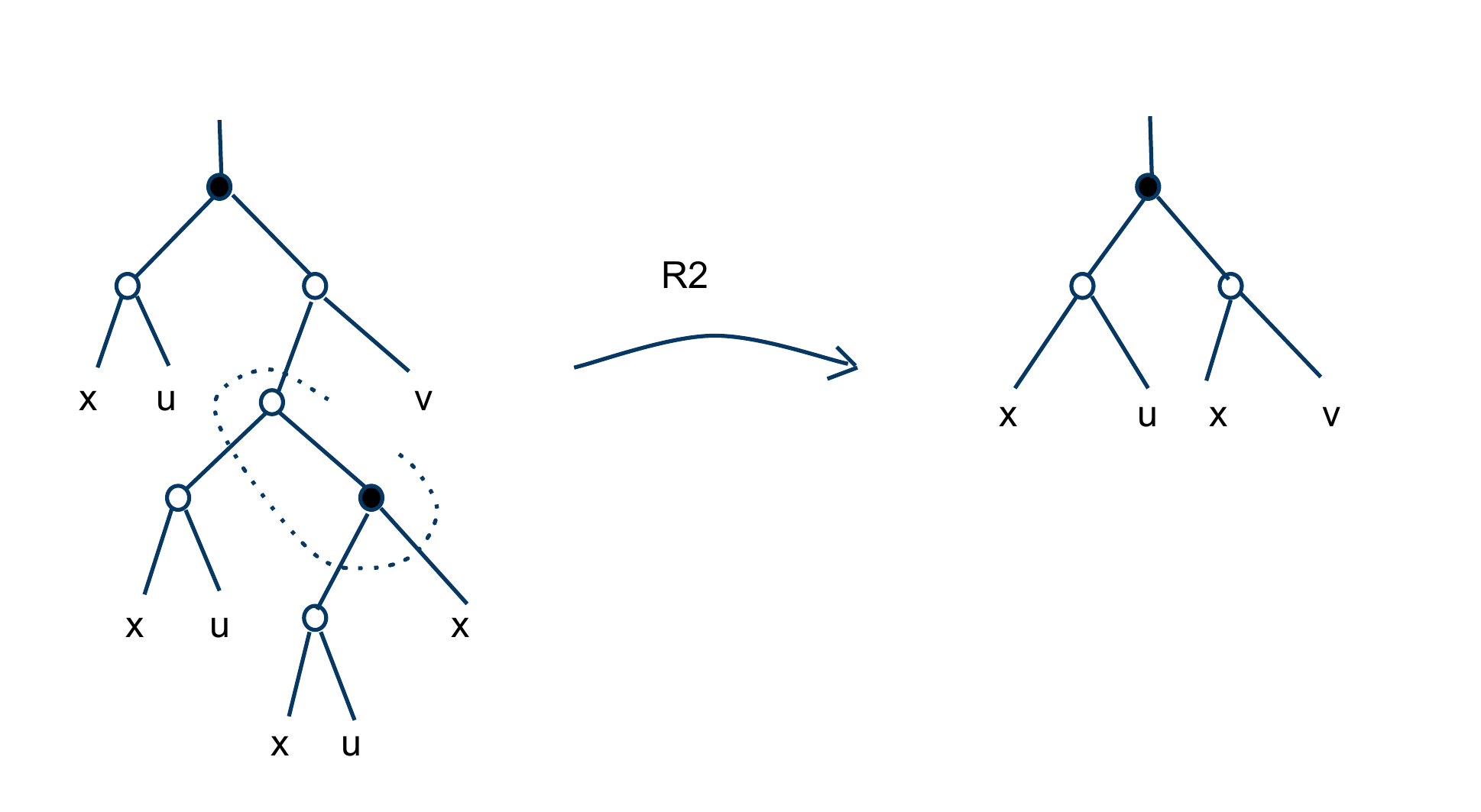}}

Relation (e) (which is  a kind of associativity relation) is obtained from: 

\centerline{\includegraphics[width=120mm]{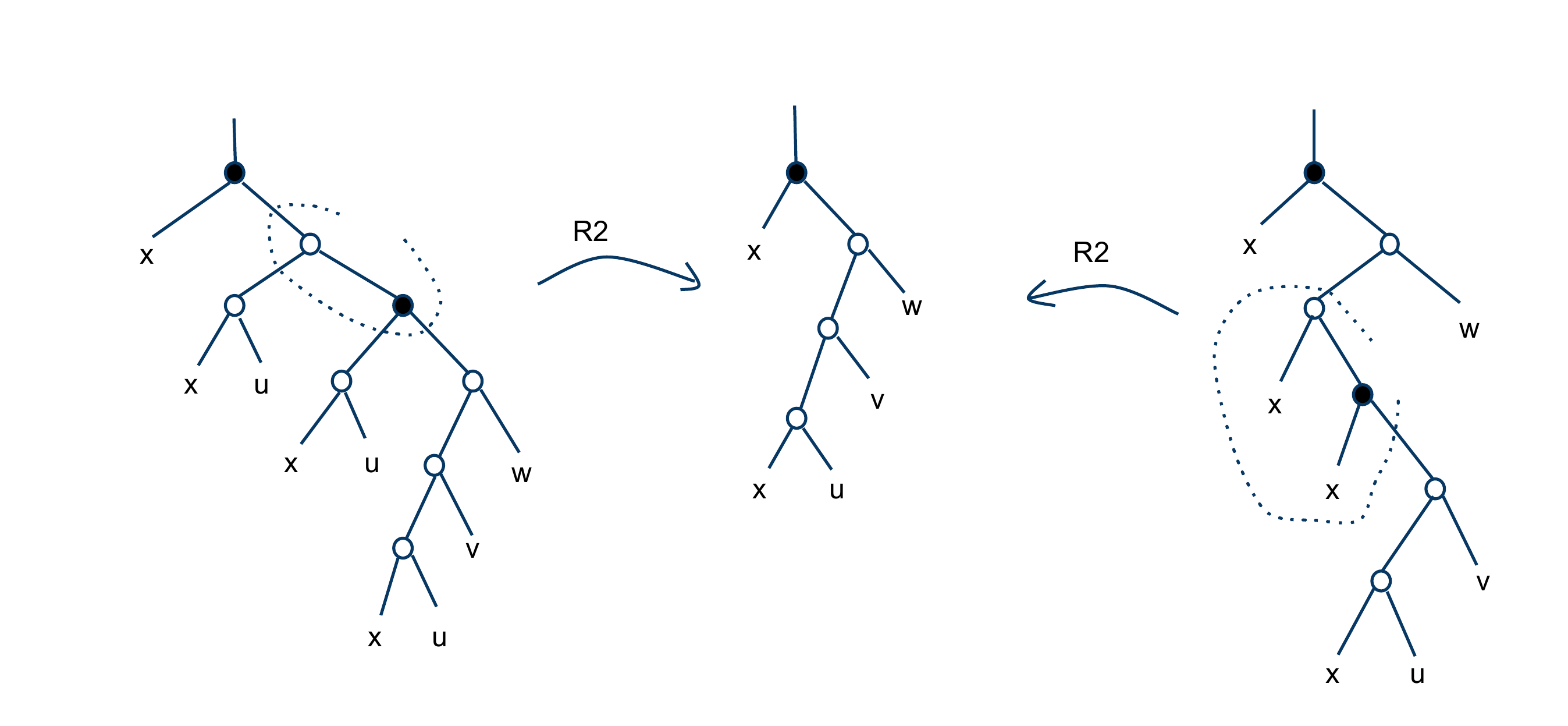}} 

Finally, for proving  relation (g) we use also the rule (R1). 

\centerline{\includegraphics[width=80mm]{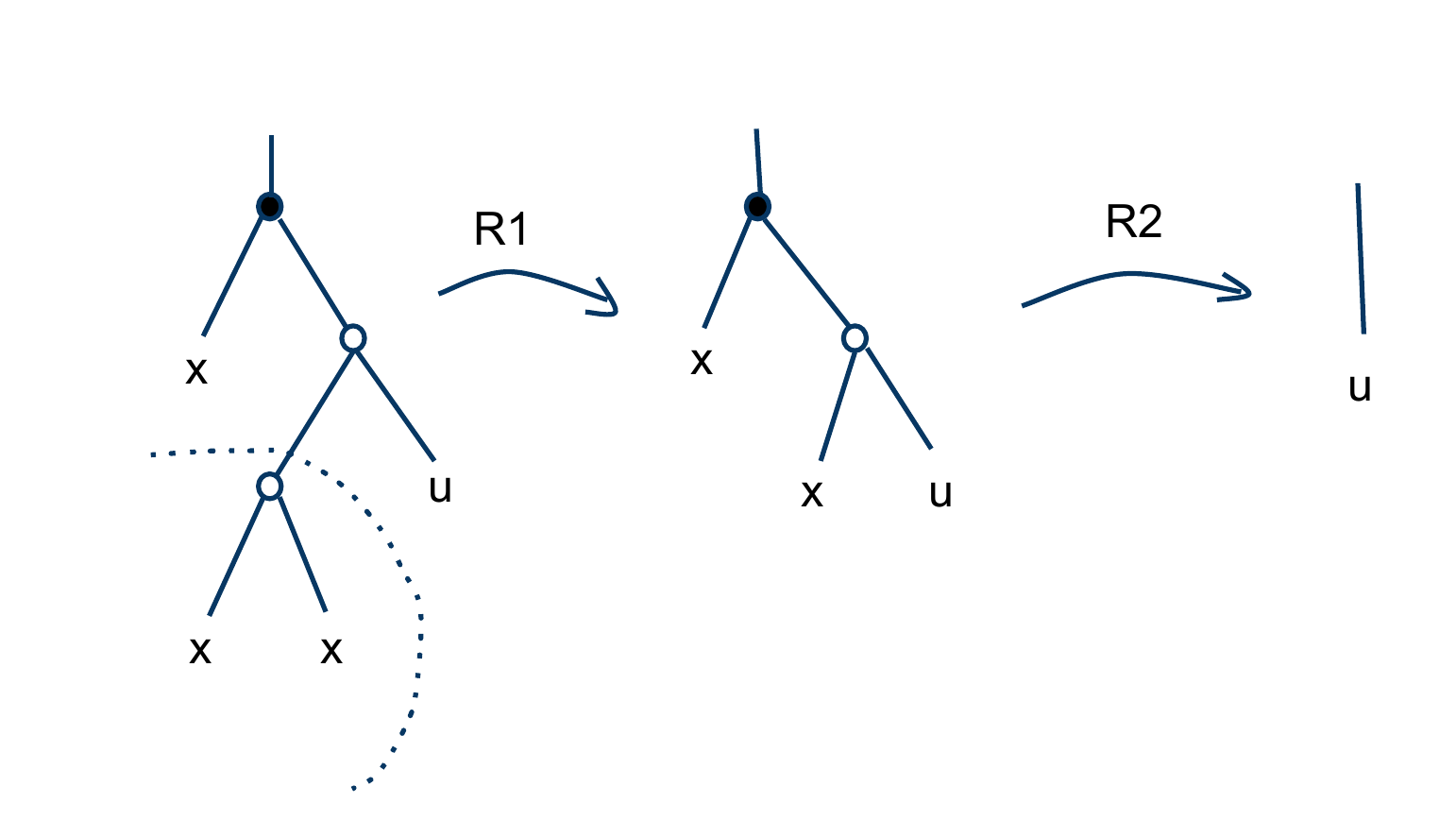}}

\paragraph{ Emergent algebras.}
See \cite{buligabraided} for all details. 

\begin{definition}
A $\Gamma$-uniform irq, or emergent algebra  $(X, \circ, \bullet)$ is a separable uniform  
space $X$ which is also a $\Gamma$-irq, with continuous operations,  such that: 
\begin{enumerate}
\item[(C)] the operation $\circ$ is compactly contractive: for each compact set 
$K \subset X$ and open set $U \subset X$, with $x \in U$, there is an open set 
$\displaystyle A(K,U) \subset \Gamma$ with $\mu(A) = 1$ for any $\mu \in 
Abs(\Gamma)$ and for any $u \in K$ and
$\varepsilon \in A(K,U)$, we have $\displaystyle x \circ_{\varepsilon} u \in U$; 
\item[(D)] the following limits exist for any $\mu \in Abs(\Gamma)$ 
$$ \lim_{\varepsilon \rightarrow \mu} \Delta_{\varepsilon}^{x}(u,v)  \, = \, 
\Delta^{x}(u,v) \quad , \quad \lim_{\varepsilon \rightarrow \mu} \Sigma_{\varepsilon}^{x}(u,v) 
 \,  = \, \Sigma^{x}(u,v) $$
and are uniform with respect to $x, u, v$ in a compact set. 
\end{enumerate}
\label{deftop}
\end{definition}

Dilation structures are also emergent algebras. In fact, emergent algebras are 
 generalizations of dilation structures, where the distance is no
longer needed. 

The main property of a uniform irq is the following. It is a consequence of 
relations from proposition \ref{pplay}. 

\begin{theorem}
Let $(X, \circ, \bullet)$ be a uniform irq. Then for any $x \in X$  the operation 
$\displaystyle (u,v) \mapsto  \Sigma^{x}(u,v) $ gives $X$ the structure of
a conical group with the dilation $\displaystyle u \mapsto x \circ_{\varepsilon} u$.
\label{mainthm}
\end{theorem}

\paragraph{Proof.} Pass to the limit in the relations from proposition
 \ref{pplay}.  We can do this exactly because of the uniformity 
 assumptions. We therefore have a series of algebraic relations which can 
 be used to get the conclusion. \hfill $\square$
 
%\newpage

\section{The difference as a universal gate}

Because Reidemeister 3 moves are forbidden, the following configurations, called 
difference gates, become important. Each of them is related to either the
$\varepsilon$-sum or $\varepsilon$-difference functions.

\vspace{.5cm}

\centerline{\includegraphics[angle=270, width=0.7\textwidth]{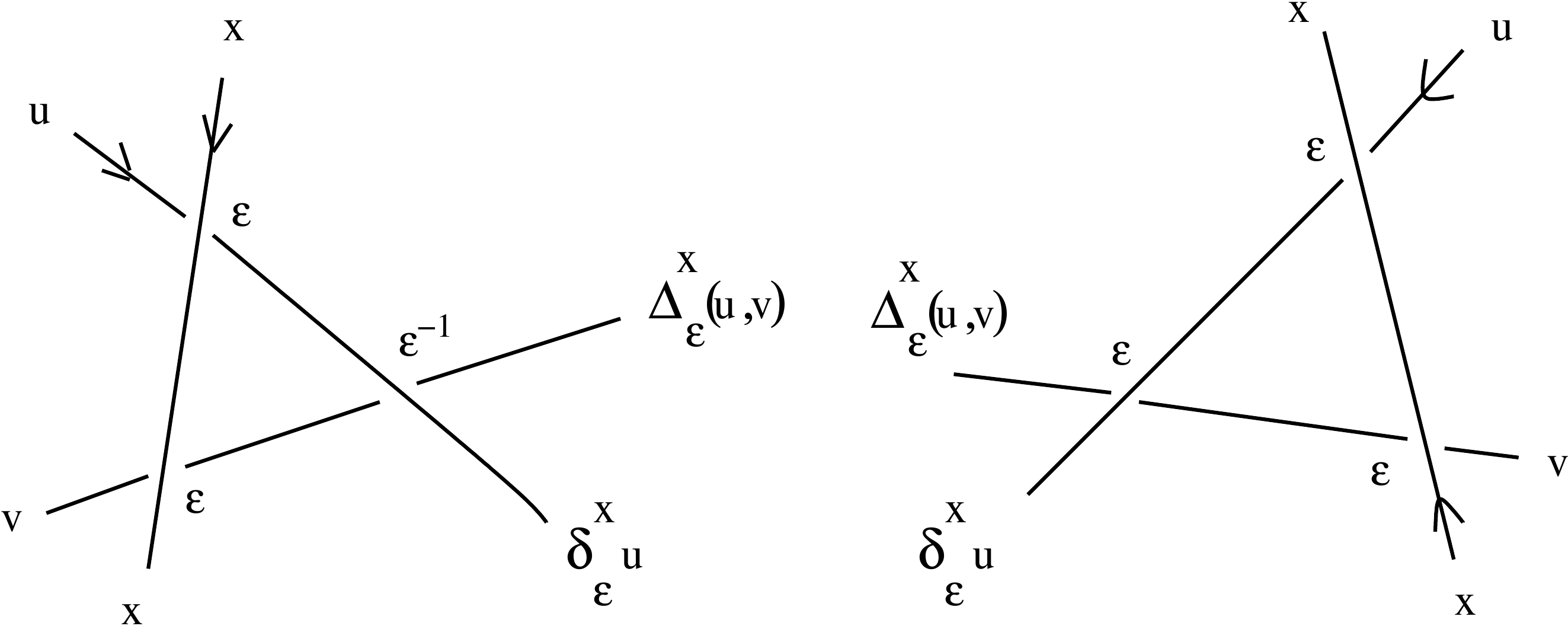}}
%\caption{}
%\label{carrier1}

\vspace{.5cm}

\vspace{.5cm}

\centerline{\includegraphics[angle=270, width=0.7\textwidth]{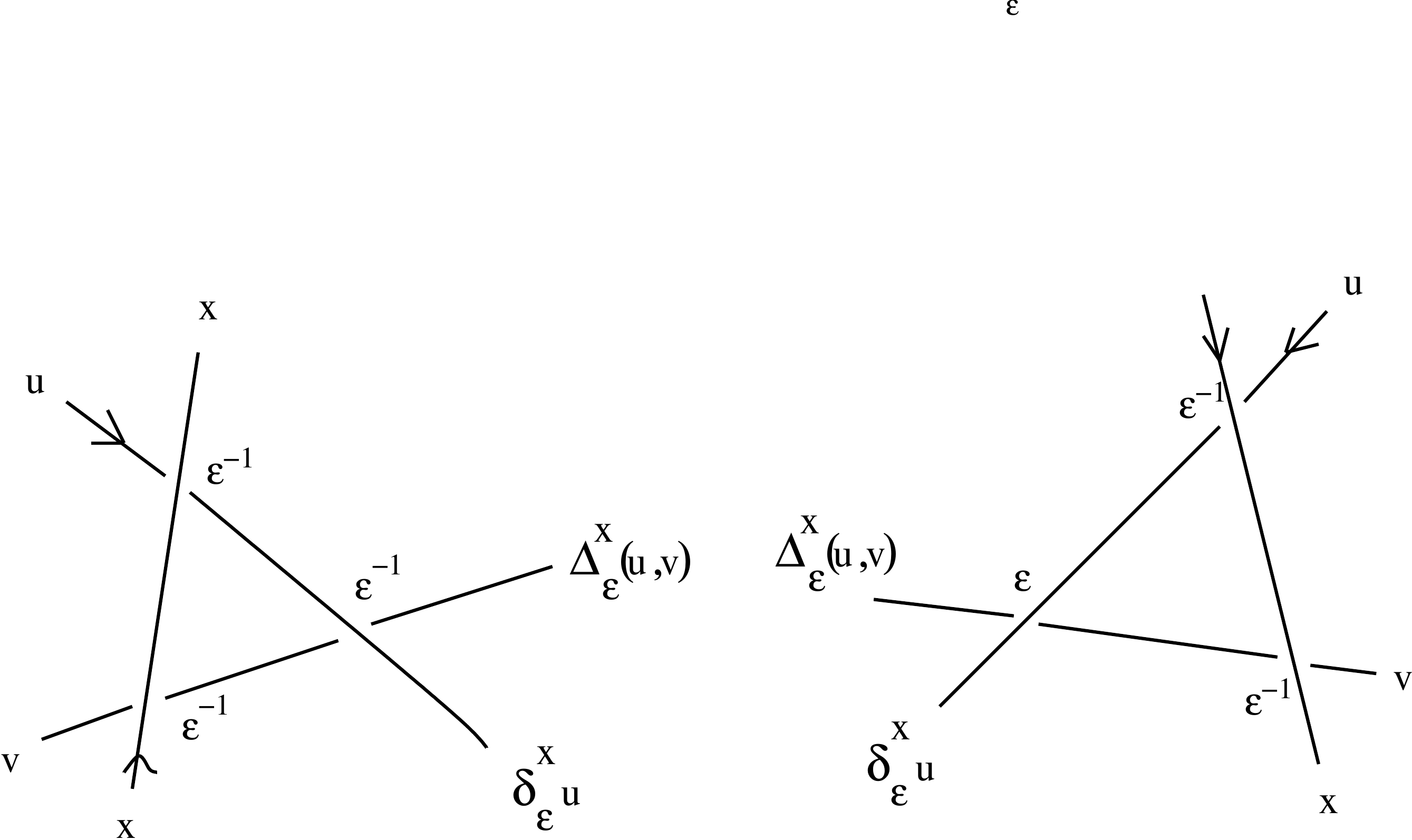}}
%\caption{}
%\label{carrier1}

\vspace{.5cm}

There are 8 difference gates. Remark that in these 4 drawings the line which
passes under the other two is not oriented. This line can be oriented in two
possible ways, therefore to each of the 4 figures correspond two figures with
all lines oriented. Also, depending on the choice of the input-output,
differences are converted into sums.

We may take this diagrams as primitives, as universal gates. It is enough to
notice that a decorated crossing can be expressed as a difference diagram, or
gate. 

\vspace{.5cm}

\centerline{\includegraphics[angle=270, width=0.7\textwidth]{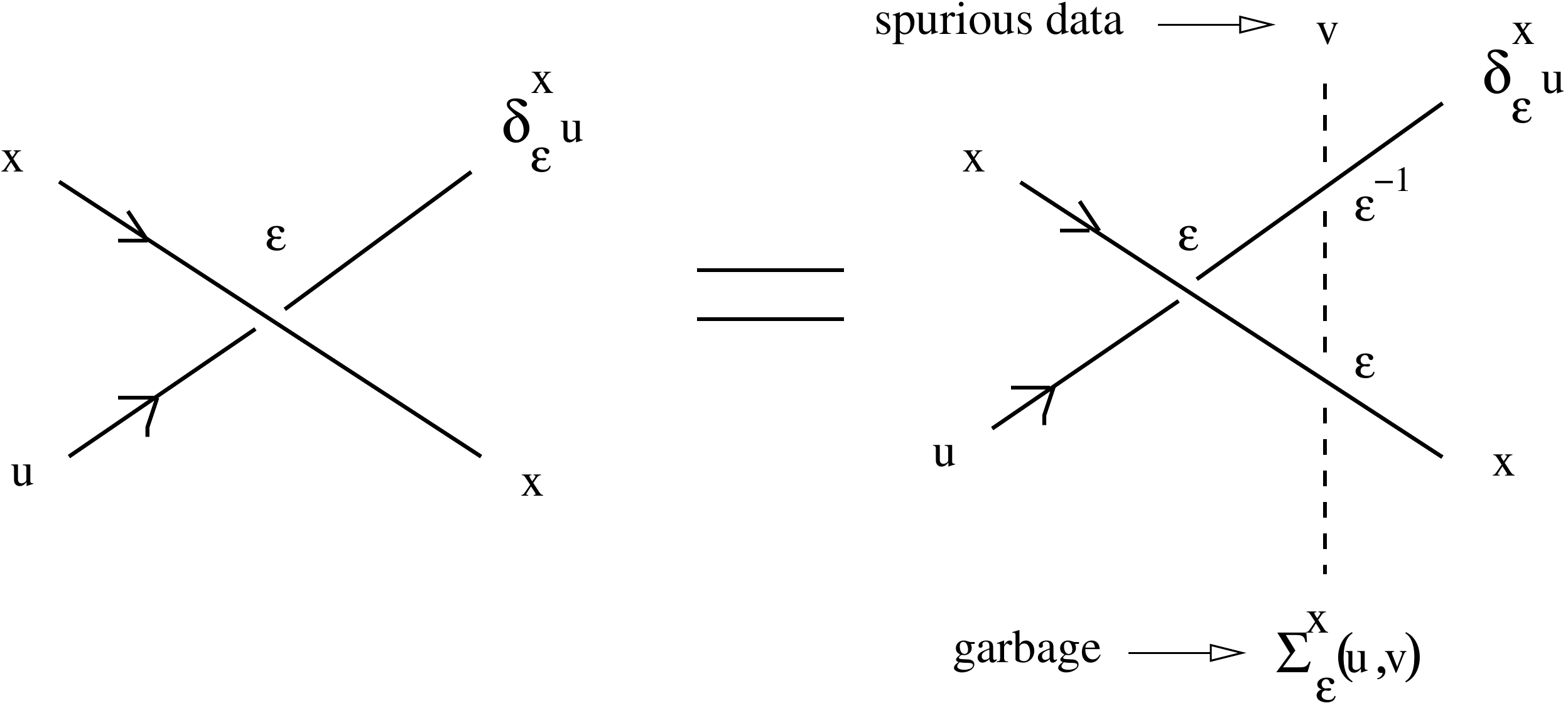}}
%\caption{}
%\label{carrier1}

\vspace{.5cm}

In the spirit of Conservative Logic of Fredkin and Toffoli \cite{toffoli}, we can
see a decorated tangle diagram as a circuit build  from several gates, all
derived from the universal difference gates. From this point of view, the
following gate could be seen as a kind of  approximate fan-out gate, 
called $\varepsilon$-fan-out.

\vspace{.5cm}

\centerline{\includegraphics[angle=270, width=0.7\textwidth]{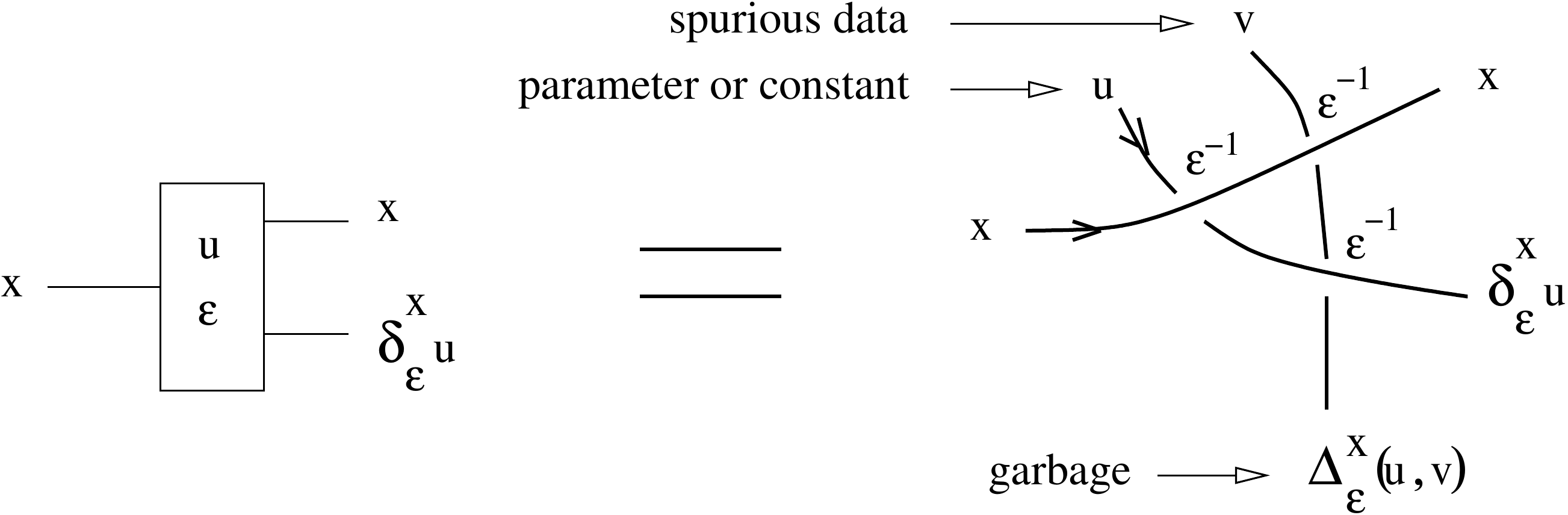}}
%\caption{}
%\label{carrier1}

\vspace{.5cm}

The $\varepsilon$-inverse gate can be seen as a NOT gate, in particular because 
it is involutive.

\vspace{.5cm}

\centerline{\includegraphics[angle=270, width=0.7\textwidth]{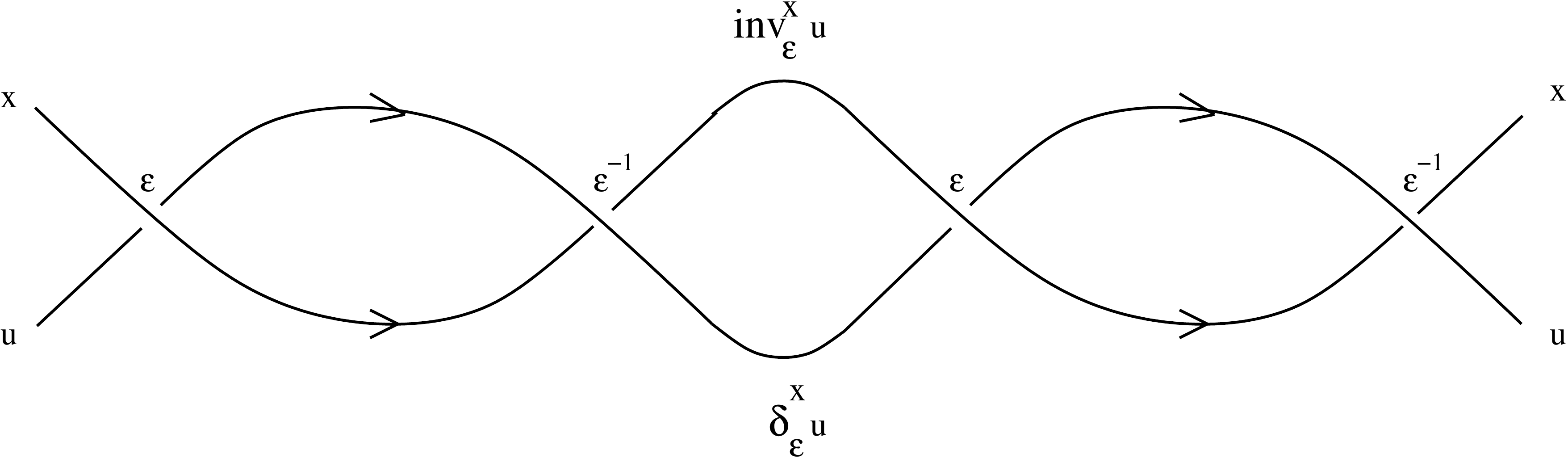}}
%\caption{}
%\label{carrier1}

\vspace{.5cm}

A difference is constructed from three crossings. Each crossing in turn may be 
expressed by a difference gate. By recycling  data and garbage, we see that a difference is self-similar: 
it decomposes  into three differences, as shown in the next figure. 

\vspace{.5cm}

\centerline{\includegraphics[angle=270, width=0.6\textwidth]{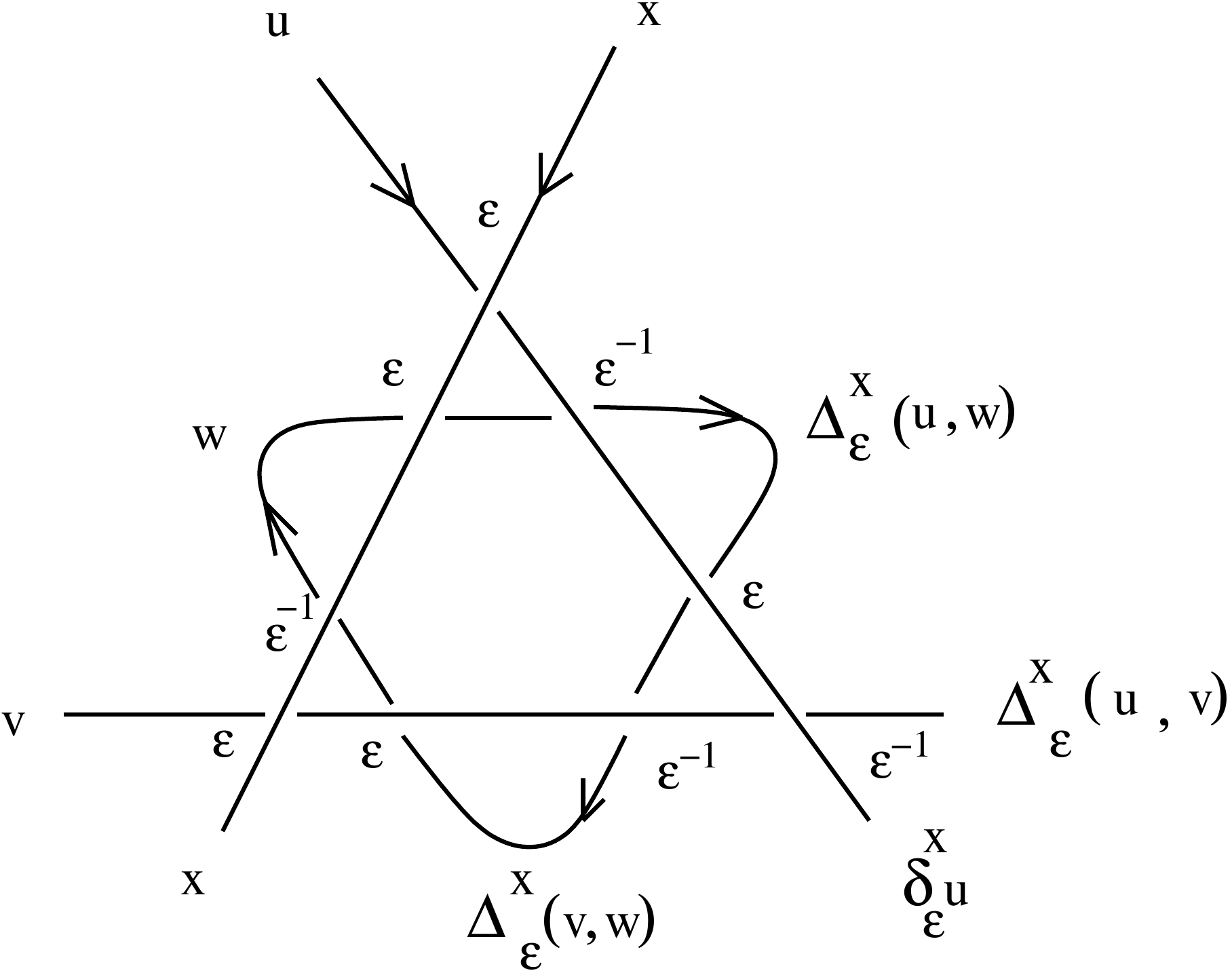}}
%\caption{A carrier is self-similar: it decomposes into three carriers}
%\label{carrierselfsimilar}

\vspace{.5cm}

%\newpage

\section{The chora}

\subsection{Definition of chora and nested tangle diagrams}

The origin of the name "chora" is in Plato' dialogue Timaeus, see relevant
quotations form this dialogue in the section dedicated to it. 
The name "chora" means "place" in ancient greek. In the light of the introduction, indeed a chora 
decorated with $x \in X$ and scale parameter $\varepsilon$ seems appropriate 
for defining a part of the territory $X$ around $x$, seen at scale 
$\varepsilon$. 
A chora is a oriented tangle diagram with decorated crossings which can be
transformed after a finite number of moves into a diagram like the one figured
next. 

%\vspace{.5cm}

\centerline{\includegraphics[angle=270, width=0.5\textwidth]{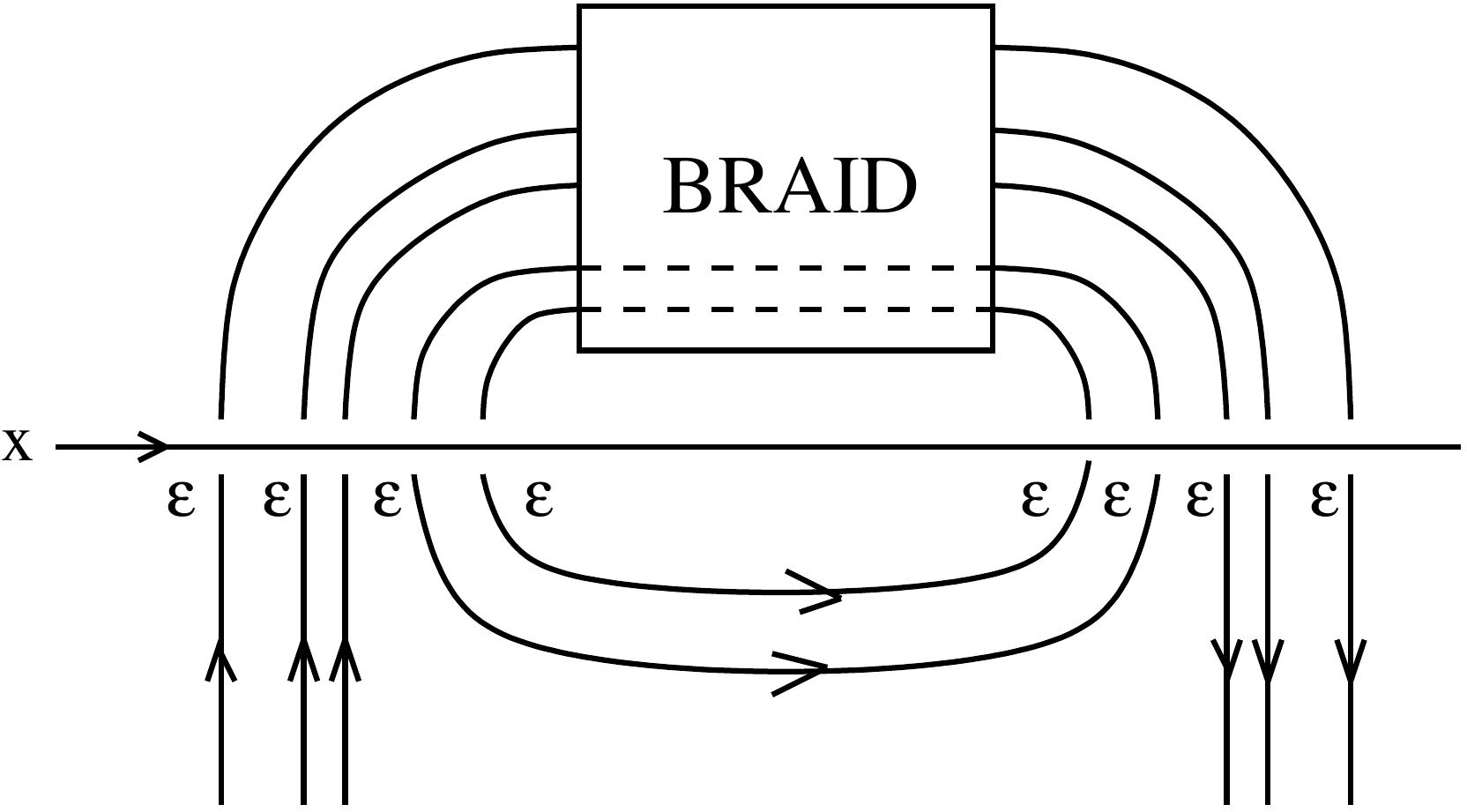}}
%\caption{}
%\label{notchora1}

\vspace{.5cm}

The chora has a "boundary", which is decorated in the figure by $x \in X$ 
and an "interior", which is a braid, such that several conditions are
fulfilled. 

\paragraph{Boundary of a chora.} As is seen in the preceding figure, the
boundary of a chora is equivalent by re-wiring to a simple closed curve oriented in 
the counter-clock sense, decorated by an element $x \in X$, such that it passes
over in all crossings and it is not involved in any virtual crossing.
Moreover, all crossings where the boundary is involved are decorated with the
same scale variable ($\varepsilon \in \Gamma$ in the figure). 

After re-wiring, a chora (and its boundary) looks like this: 

\vspace{.5cm}

\centerline{\includegraphics[angle=270, width=0.7\textwidth]{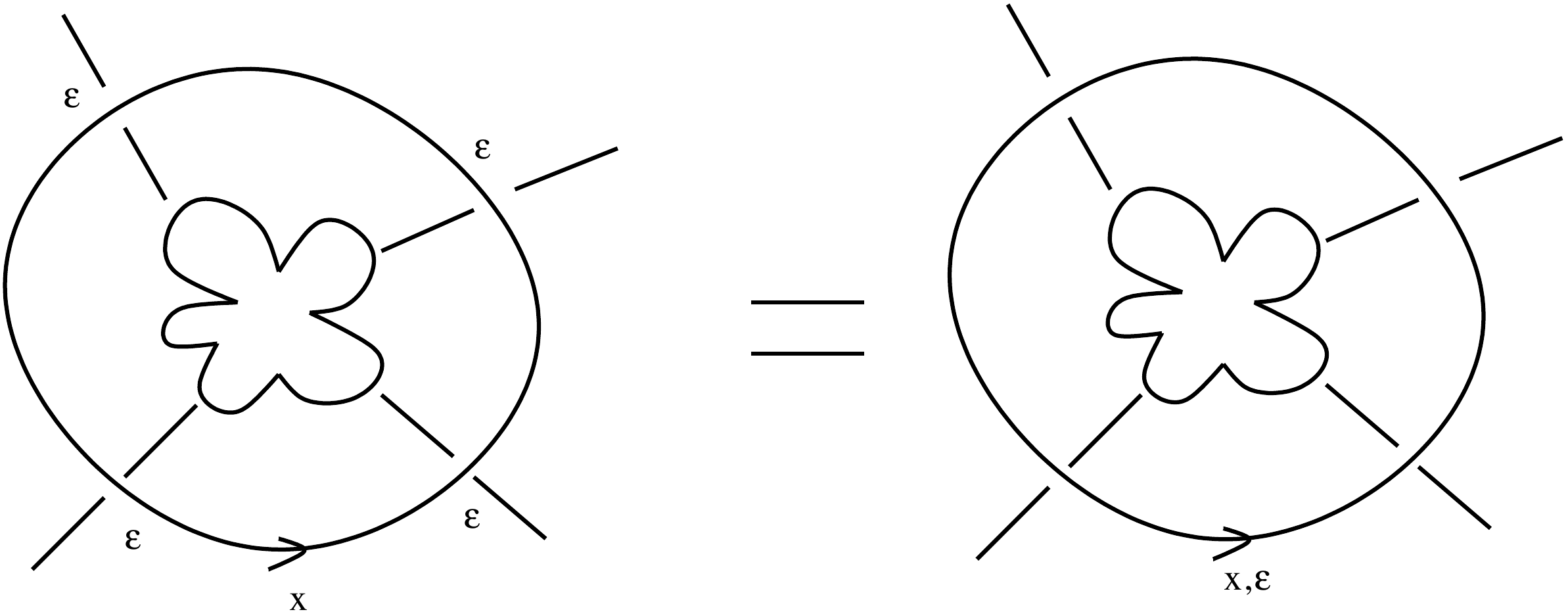}}
%\caption{}
%\label{notchora1}

\vspace{.5cm}

As explained in the first section, we adopt the following notation for a chora: 
we choose not to decorate the  crossings with the simple closed curve obtained
by joining segments decorated with $x$, because 
they all have the same decoration $\varepsilon$, but instead decorate the 
whole curve with $x,\varepsilon$. 

\paragraph{Interior of a chora.} The interior of a chora is the 
diagram located in the bounded region of the plane encircled by the close
curve which defines the chora. The diagram is such that it can be transformed by 
Reidemeister II moves into a braid, which has an input (segments entering 
into the chora), output (segments going out from the chora) and parameters 
which decorate input segments connected with output segments (the arcs figured
in the first drawing of a chora, which connect inputs with outputs). 

\paragraph{Nested choroi.} A nested configuration of choroi ("choroi" is the 
greek plural of "chora") is a collection of simple closed curves
satisfying the conditions of a boundary of a chora, such that for any two such
curves either the interior of one is inside the other, or they have disjoint
interiors. 

Here is an example of a diagram containing three choroi, in a nested
configuration. (The irrelevant decorations of the segments were ignored.)

%\vspace{.5cm}

\centerline{\includegraphics[angle=270, width=0.7\textwidth]{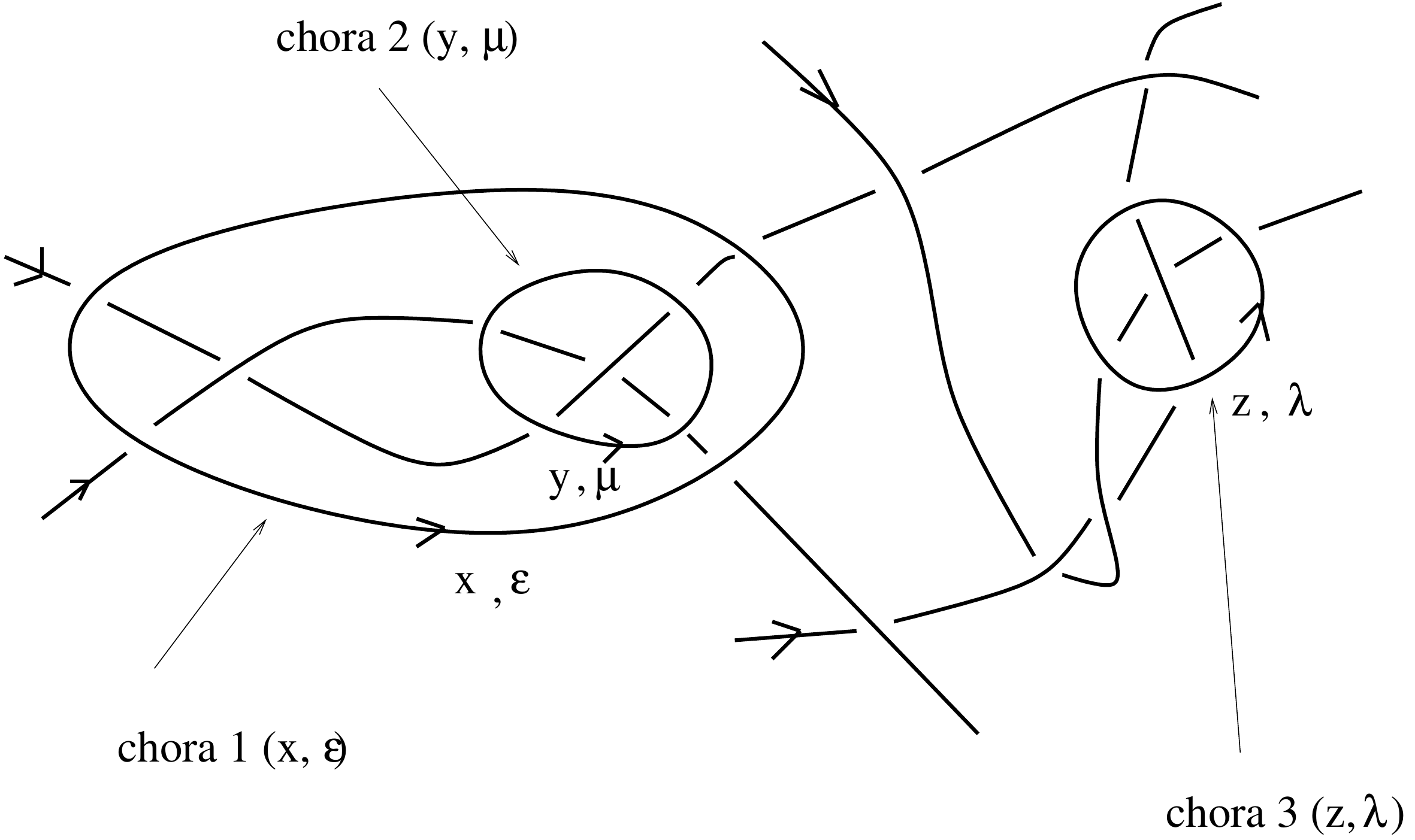}}
%\caption{}
%\label{notchora1}

%\vspace{.5cm}

\paragraph{On parameters.} By definition, a simple chora is one without
parameters, i.e. without arcs connecting inputs with outputs. We want to be
able to have nested configurations of choroi. For this we need parameters. 

Indeed, look at the preceding figure, at the chora $1(x,\varepsilon)$. Let us
transform it into a chora diagram as the one figured at the beginning of this
section. At the beginning the diagram looks like this. "1", "2" are the inputs,
"3", "4" are the outputs. 

\vspace{.5cm}

\centerline{\includegraphics[angle=270, width=0.5\textwidth]{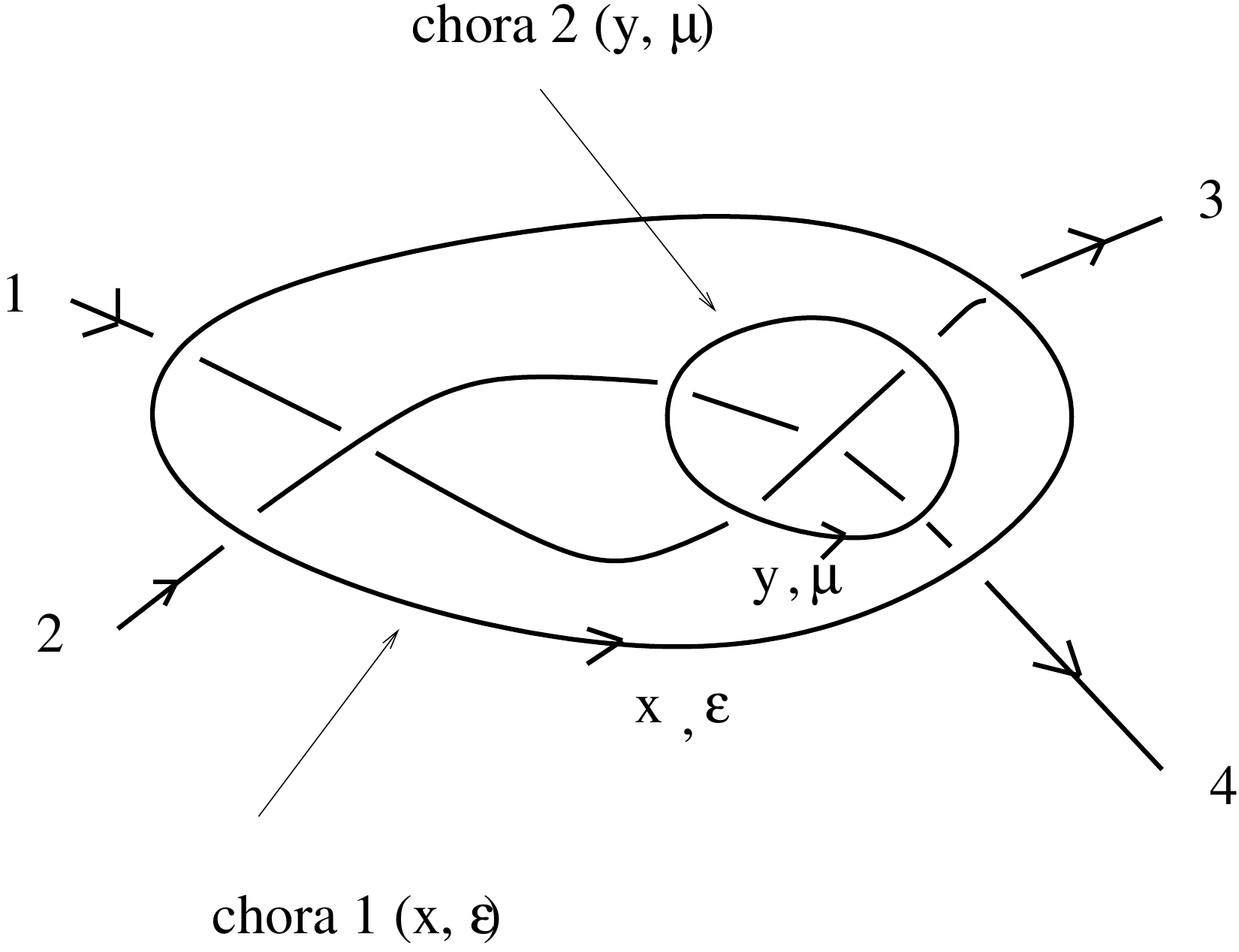}}
%\caption{}
%\label{notchora1}

\vspace{.5cm}

After re-wiring, the diagram starts to look like a chora, only that the
interior is not a braid. 

\vspace{.5cm}

\centerline{\includegraphics[angle=270, width=0.5\textwidth]{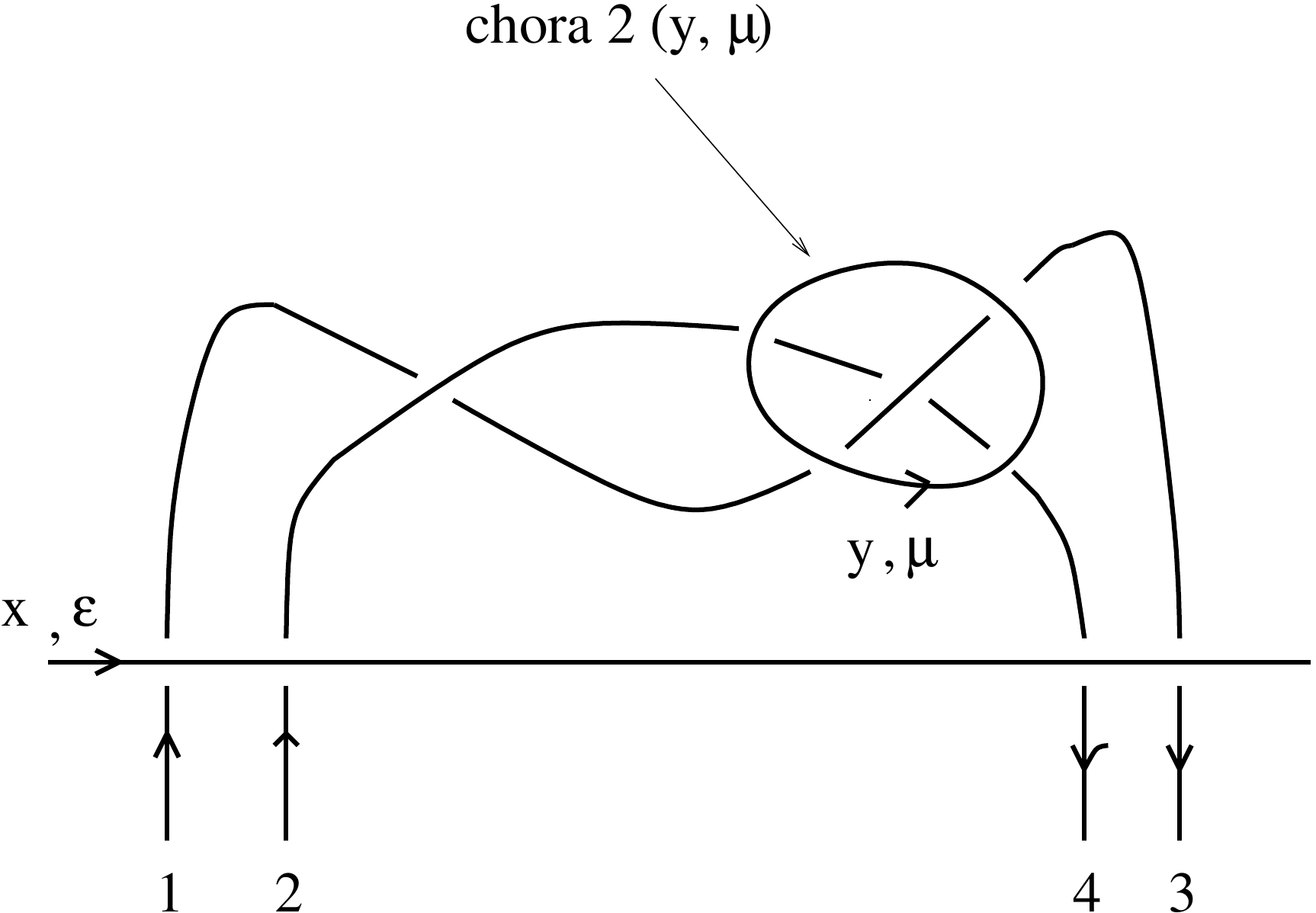}}
%\caption{}
%\label{notchora1}

\vspace{.5cm}

We pull the boundary of chora $2(y,\mu)$ down, slide it under the boundary 
of chora $1(x,\varepsilon)$ and we finally obtain a figure like the first one
in this section.

\vspace{.5cm}

\centerline{\includegraphics[angle=270, width=0.5\textwidth]{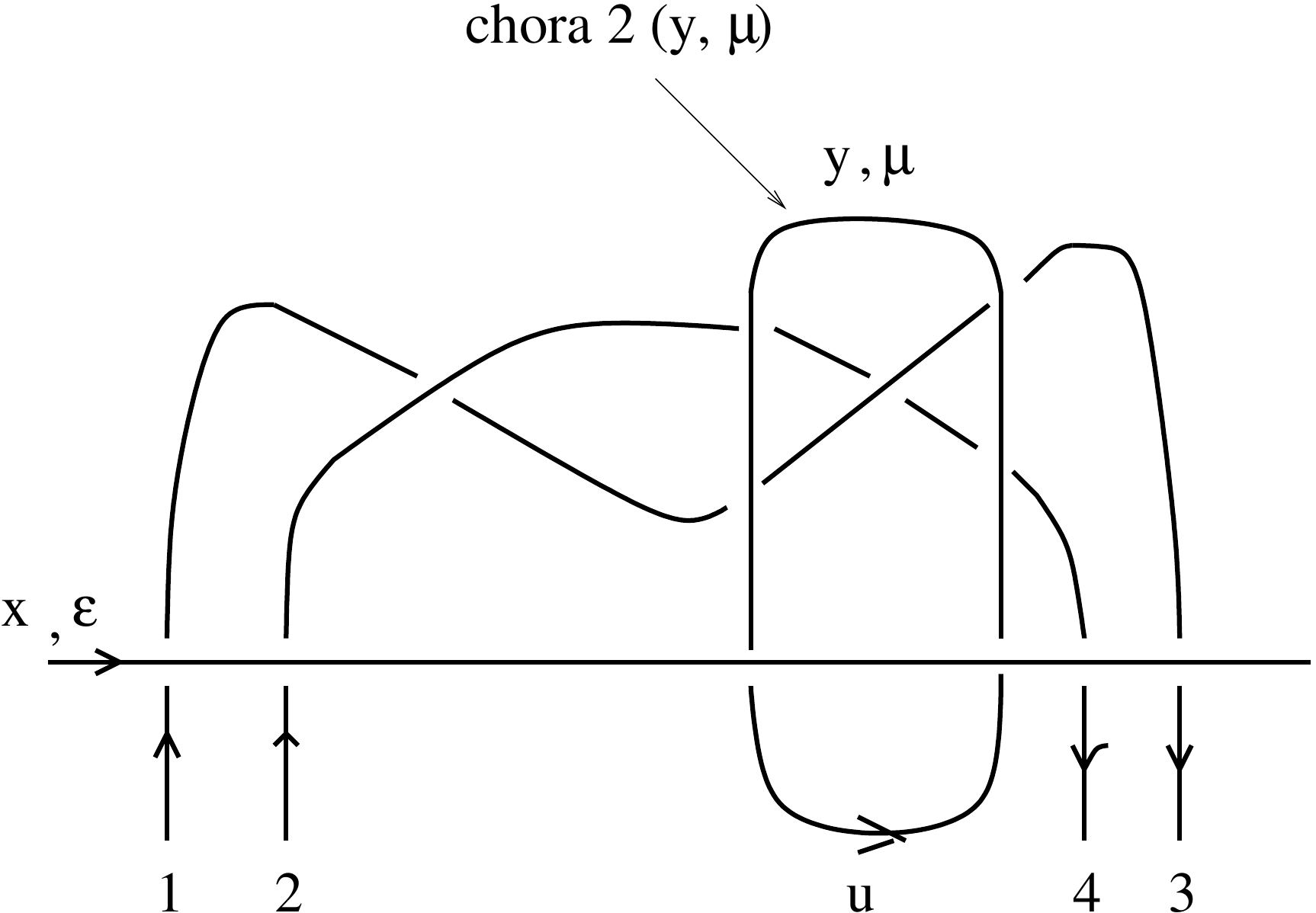}}
%\caption{}
%\label{notchora1}

\vspace{.5cm}

The new segment which appeared is decorated by  the parameter $u$. By the rules
of decoration, we have a relation between $x, u, y, \varepsilon$: 
$$ y = \delta^{x}_{\varepsilon} u $$
so the chora $2(y,\mu)$ is in fact decorated by $\displaystyle
(\delta^{x}_{\varepsilon} u, \mu)$.

Note that, however, a "chora" is a notation on the tangle diagram. We choose 
to join segments in a simple closed curve, but this  could be done in several 
ways, starting from the same decorated tangle diagram. In this sense, as 
Plato writes, a chora "is apprehended without the help of sense, by a kind of 
spurious reason, and is hardly real; [...]  Of these and other things of the 
same kind, relating to the true and waking reality of nature, 
we have only this dreamlike sense, and we are unable to cast off sleep and 
determine the truth about them."

For example, look at the first  diagram from the next figure. 
This diagram is  called an "elementary chora". 
Such a diagram could be obtained from several other diagrams, among them 
the following two (second and third from the next figure), called states of  
the elementary chora.

\vspace{.5cm}

\centerline{\includegraphics[angle=270, width=0.7\textwidth]{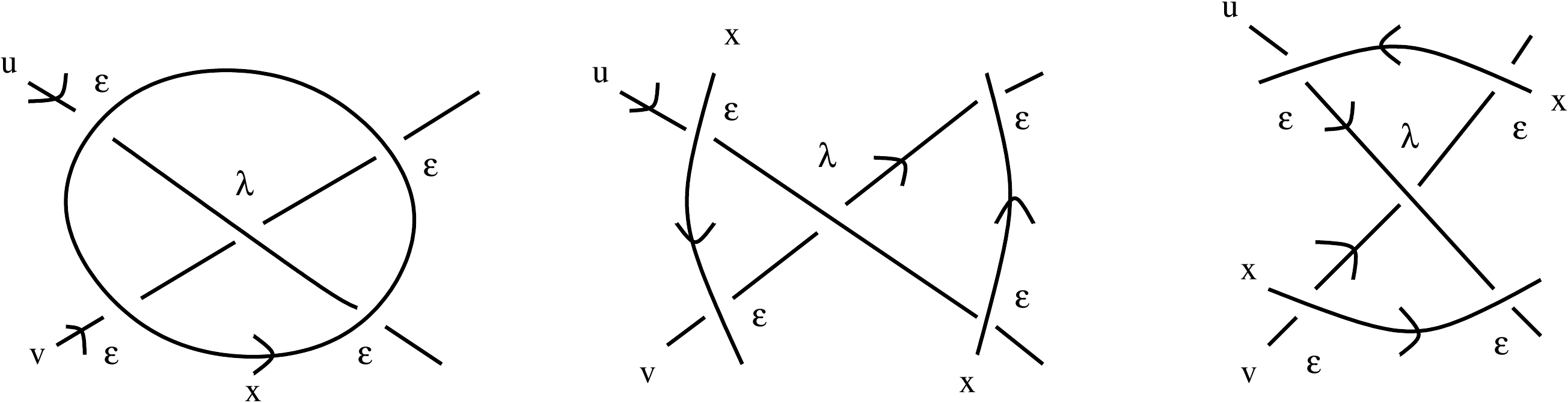}}
%\caption{An elementary chora and its two states}
%\label{elemchora2states}

\vspace{.5cm}

By making Reidemeister 2 moves, we may split a big chora diagram 
into smaller choroi. 
Conversely, we may join two choroi in a bigger chora. 
Moreover, starting from a big chora, we can split it into  
 into elementary choroi. An example is figured further.

\vspace{.5cm}

\centerline{\includegraphics[angle=270, width=0.7\textwidth]{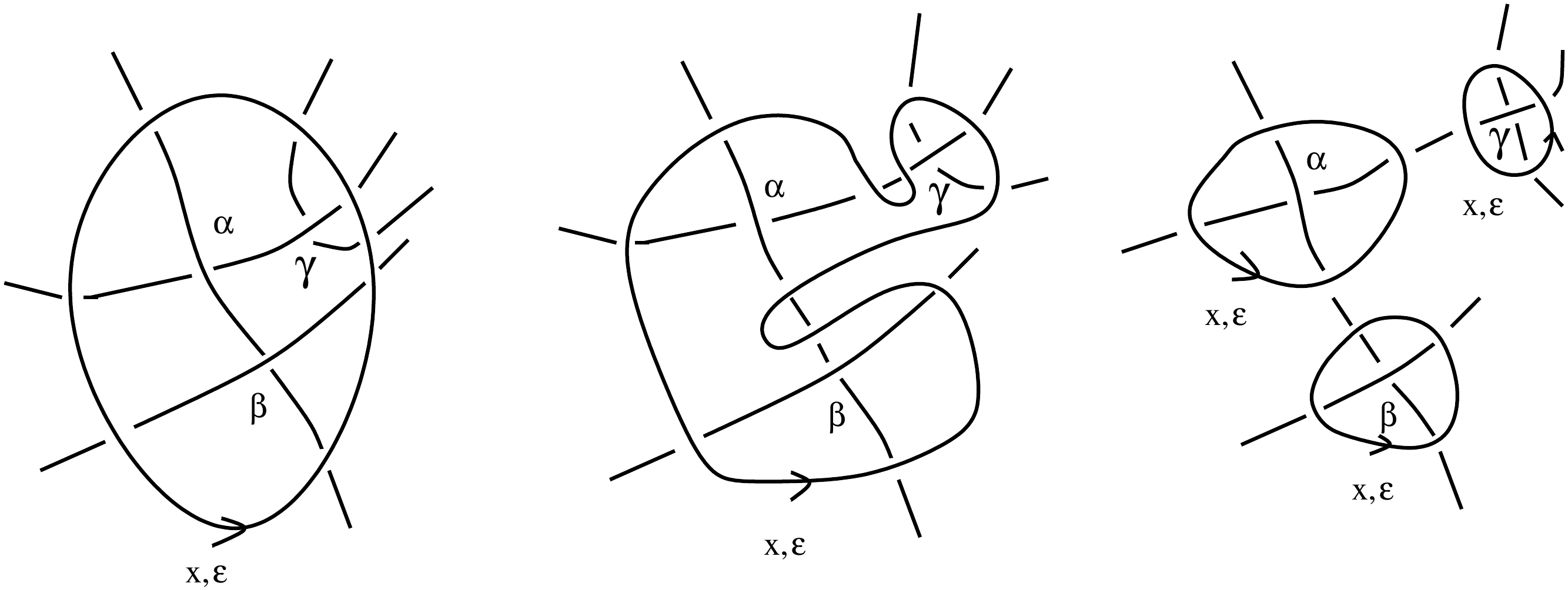}}
%\caption{}
%\label{notchora2}

\vspace{.5cm}

\subsection{Decompositions into choroi and differences}

An elementary chora can be decomposed in two ways into 
two differences and a crossing. Therefore an elementary chora is a particular
combination of three universal difference gates.

\vspace{.5cm}

\centerline{\includegraphics[angle=270, width=0.7\textwidth]{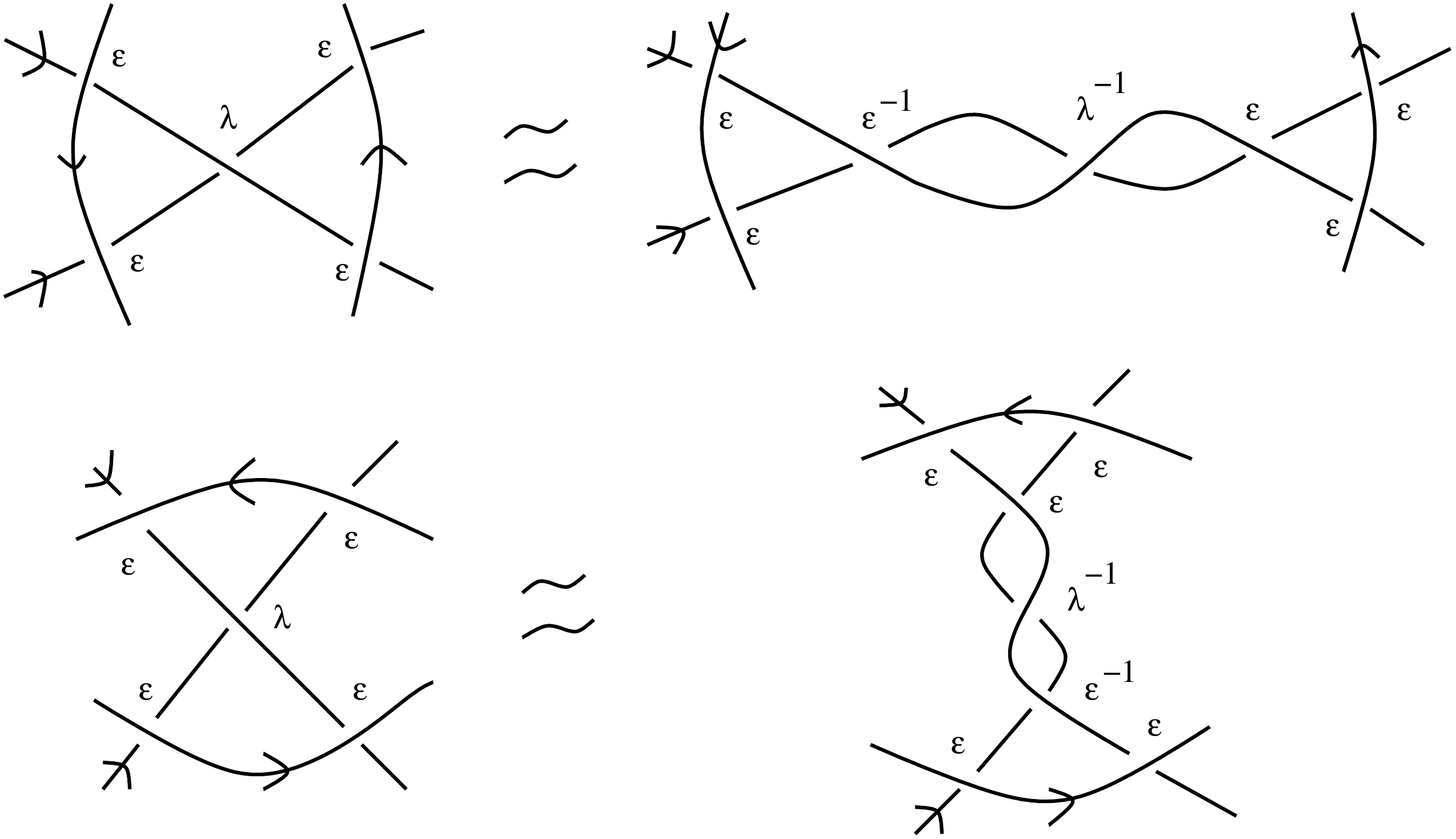}}
%\caption{}
%\label{choradecomposition}

\vspace{.5cm}

Likewise, the chora inside a chora diagram  decomposes into a chora and four 
difference gates.

\vspace{.5cm}

\centerline{\includegraphics[angle=270, width=0.8\textwidth]{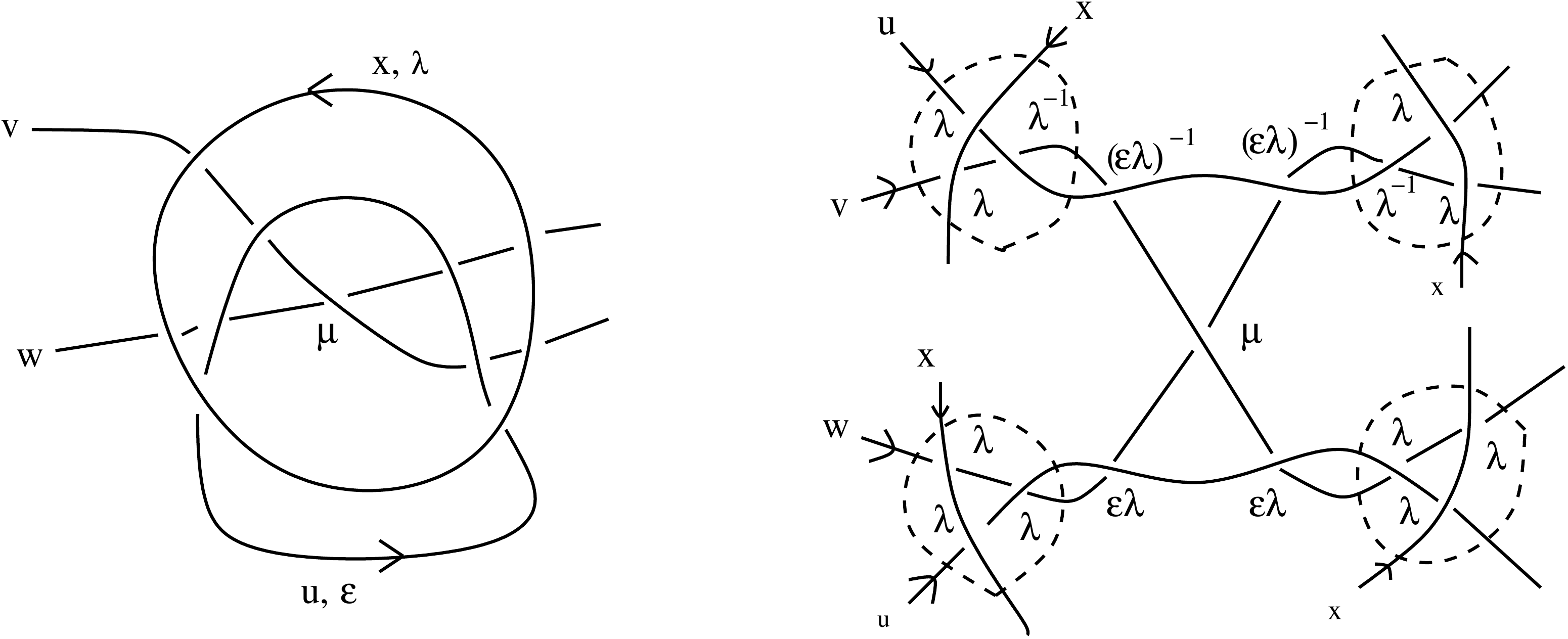}}
%\caption{}
%\label{chorainchora3}

\vspace{.5cm}

Here I applied repeatedly a trick based on the Reidemeister II move, look at the
small parts of the diagram encircled by dashed curves.

\vspace{.5cm}

\centerline{\includegraphics[angle=270, width=0.5\textwidth]{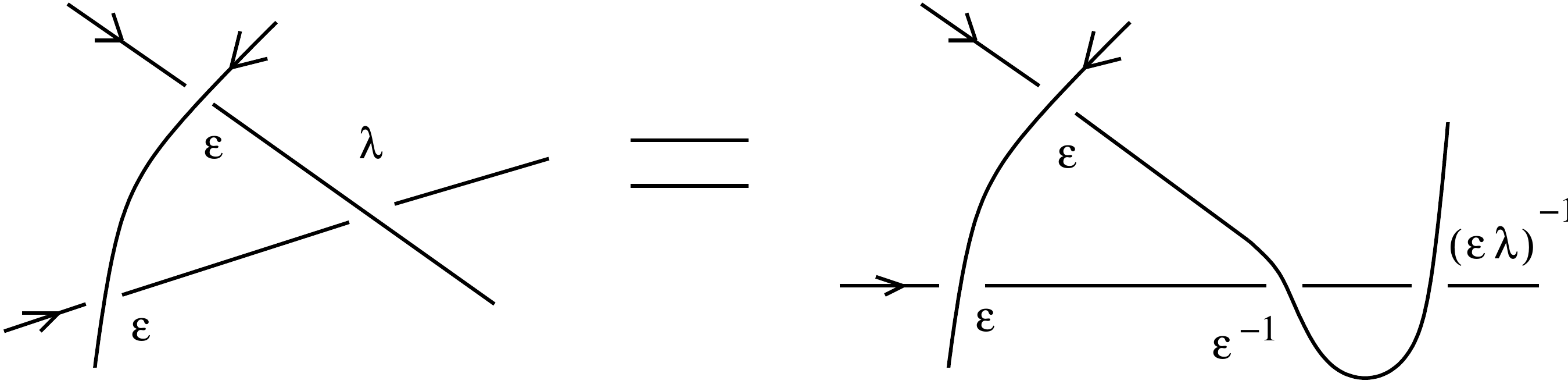}}
%\caption{}
%\label{chorainchora3}

\vspace{.5cm}

Imagine now that we contemplate a diagram presenting a nested collection of
choroi. In order to simplify it as much as possible, we can apply the 
chora-inside-a-chora decomposition, but this operation could leave us with 
difference-inside-a-chora diagrams. 

\vspace{.5cm}

\centerline{\includegraphics[angle=270, width=0.45\textwidth]{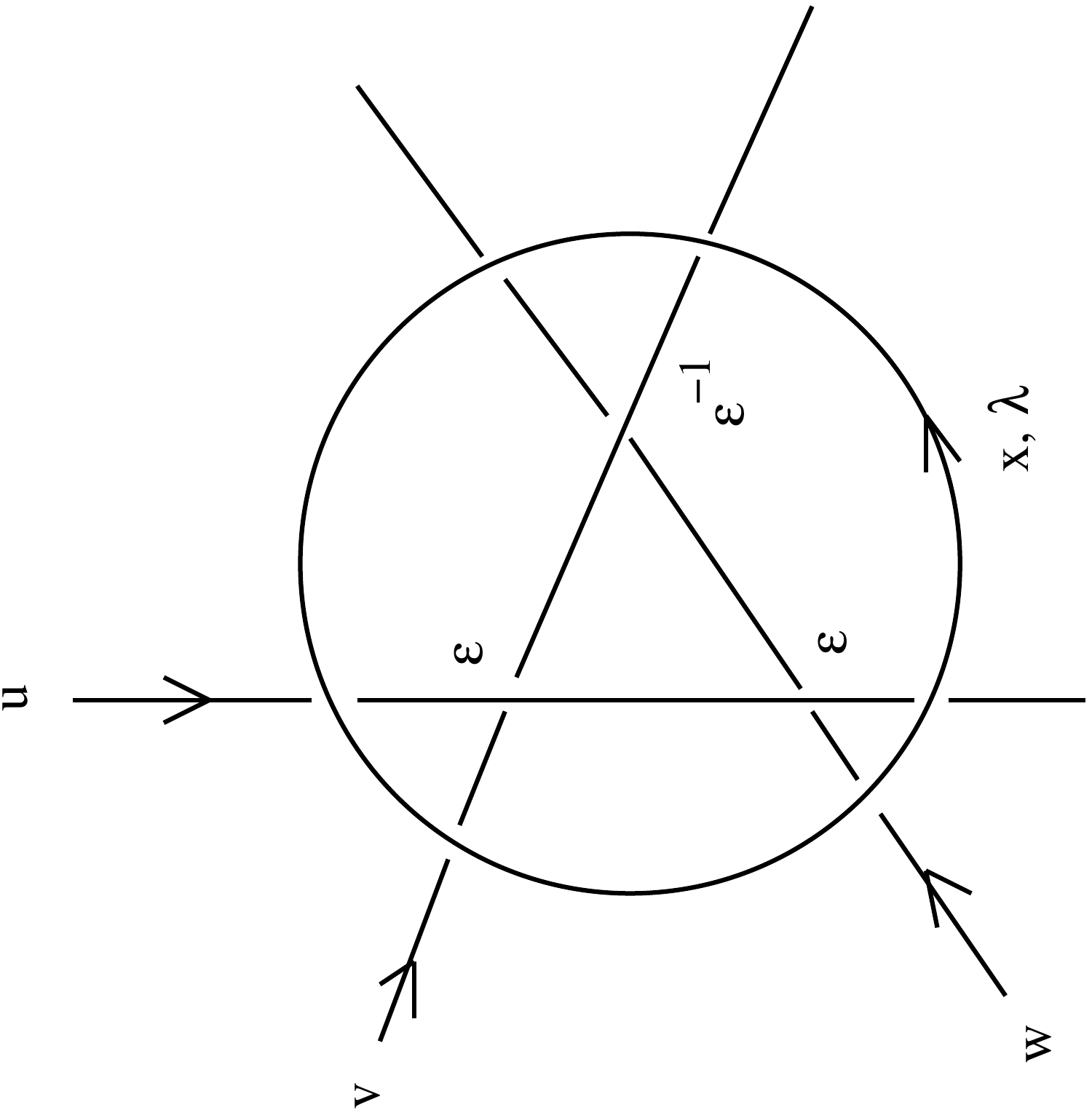}}
%\caption{}
%\label{chorainchora3}

\vspace{.5cm}

 By applying a trick involving Reidemeister II moves, resembling with the
 preceding one, we find out that the difference-inside-a-chora decomposes as four differences.

%\vspace{.5cm}

\centerline{\includegraphics[angle=270, width=0.6\textwidth]{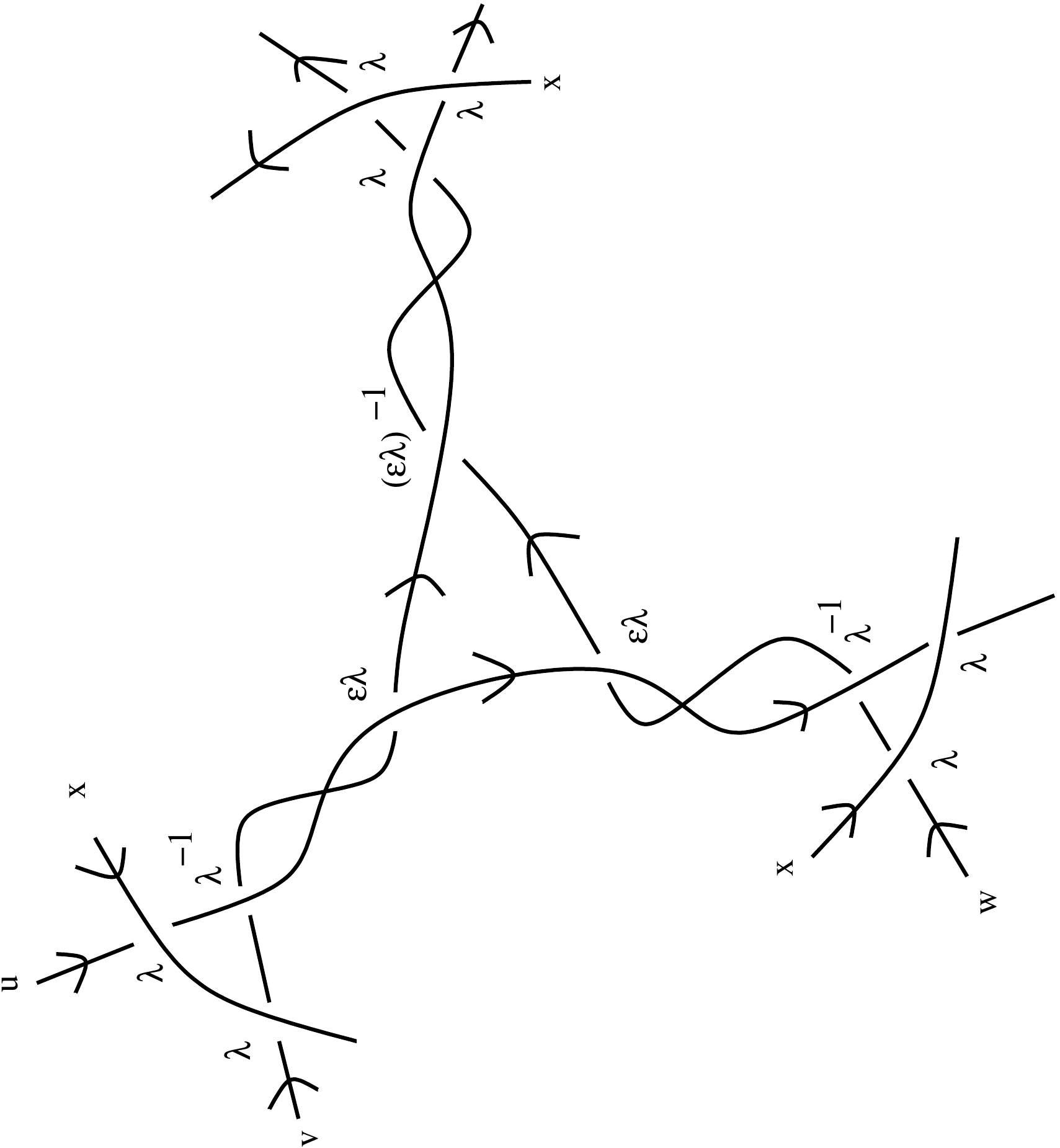}}
%\caption{Carrier inside chora. It decomposes as four carriers.}
%\label{carrierinchora}

%\vspace{.5cm}

\begin{theorem}
Let $T$ be an acceptable oriented tangle diagram formed by a collection of
nested choroi, colored with a uniform 
$\Gamma$-irq, such that any crossing is either a crossing with a 
chora or is inside a chora. For this diagram take as  scale variables and 
parameters the decorations of the choroi.    Then $T$ is equivalent with a 
diagram made only with difference gates with crossings decorated by one of the scale variables and 
with elementary choroi decorated with products of the scale variables and with 
elements of $X$ which are expressed as functions of scale variables, 
parameters, and input decorations. 
\label{thmdecomp}
\end{theorem}

\paragraph{Proof.} The hypothesis is such that we can apply repeatedly
chora-in-chora or difference-in-chora decompositions, until we transform all 
choroi into elementary choroi. We see that the difference gates which appear
will be decorated by scale variables and that the decorations from $X$ of the 
intermediate and final choroi are all expressible as outputs of diagrams  
decorated with scale variables, parameters,  input and output decorations. But 
the output decorations are themselves functions of parameters, scale parameters
and inputs. \hfill $\square$

\section{How to perform Reidemeister III moves}

In this formalism we cannot perform the Reidemeister III move. However, we can 
do this move inside a chora, in a approximate sense. 

Indeed, I shall use for this a void chora, equivalently, a particular decorated
braid which is equivalent by Reidemeister 2 moves (in the language of tangles) 
with the identity braid. In the next figure we see this braid at the left and 
the equivalent void chora at the right. 

\vspace{.5cm}

\centerline{\includegraphics[angle=270, width=0.7\textwidth]{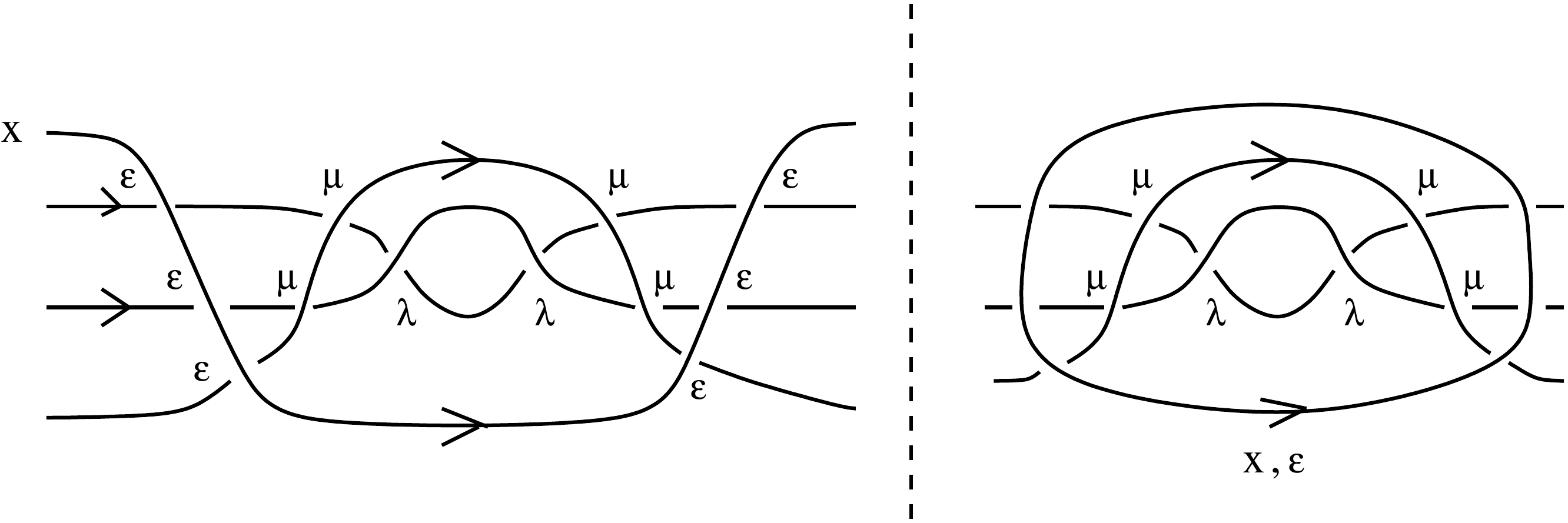}}
%\caption{Carrier inside chora. It decomposes as four carriers.}
%\label{carrierinchora}

\vspace{.5cm}

Let us put this void chora  at the right of the diagram which is the subject of 
a Reidemeister III move (inside a chora). The diagrams interact by re-wiring and Reidemeister II moves.

\vspace{.5cm}

\centerline{\includegraphics[angle=270, width=0.7\textwidth]{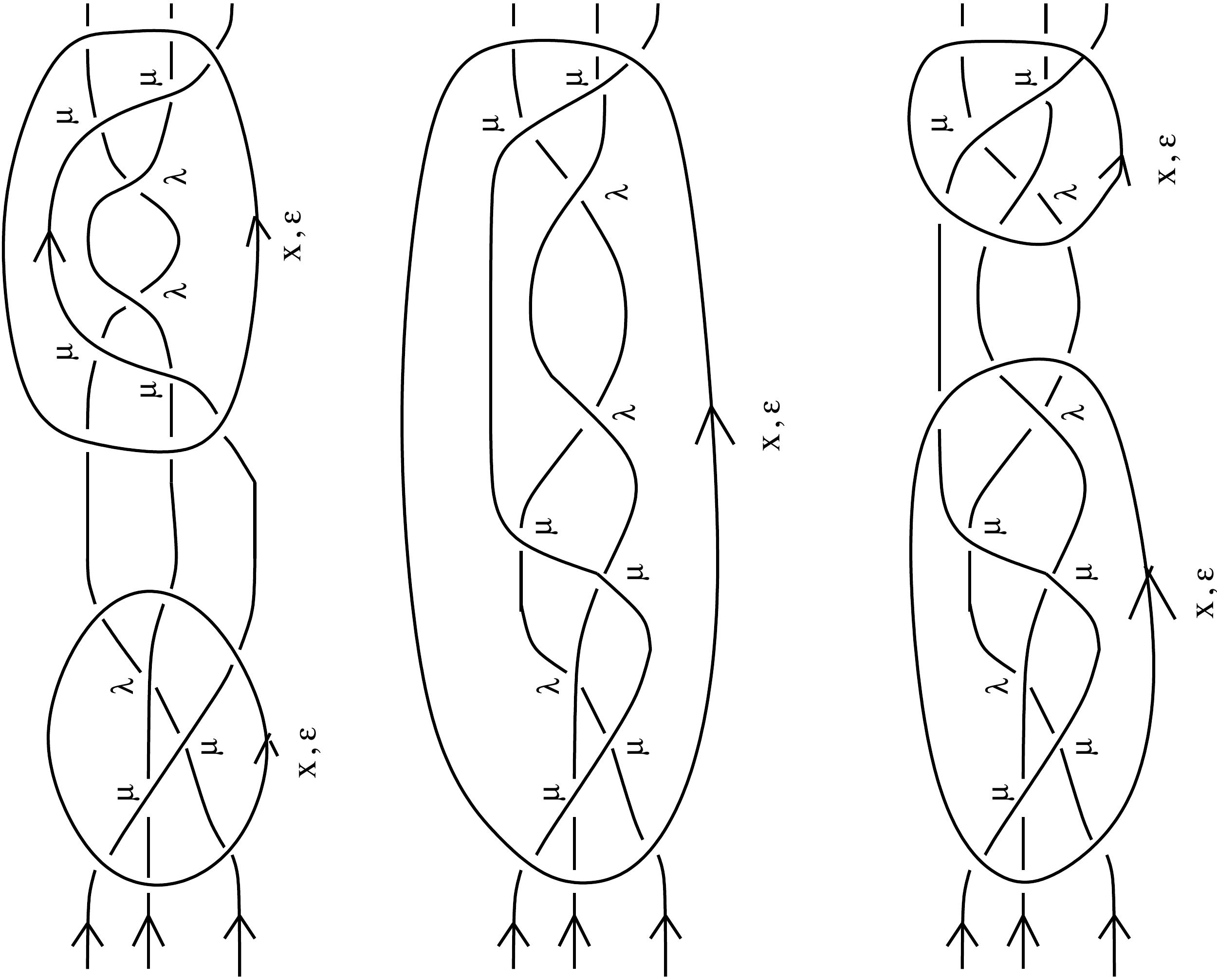}}
%\caption{Carrier inside chora. It decomposes as four carriers.}
%\label{carrierinchora}

\vspace{.5cm}

Now we see that the Reidemeister III move was performed (the diagram 
from the right of the last line), only that there is a residual diagram, 
the one from the left of the last line. I reproduce this residue, in two
equivalent forms, in the next figure. 

\vspace{.5cm}

\centerline{\includegraphics[angle=270, width=0.7\textwidth]{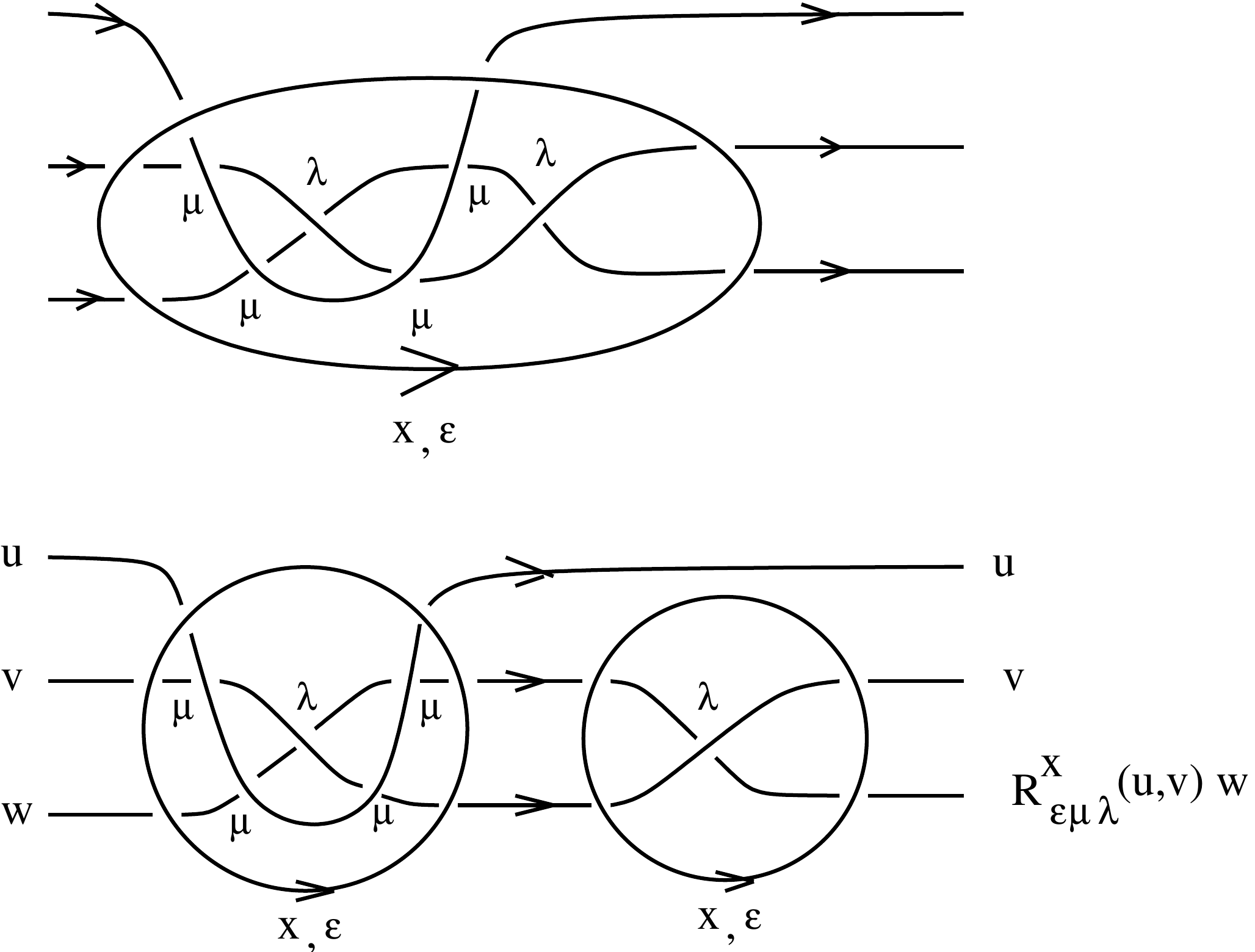}}
%\caption{Carrier inside chora. It decomposes as four carriers.}
%\label{carrierinchora}

\vspace{.5cm}

\begin{theorem}
With the notation of the last figure, we have: 
$$\lim_{\varepsilon \rightarrow 0} R^{x}_{\varepsilon \mu \lambda} (u,v) w \, =
\, w$$
uniformly with respect to $u,v,w$ in compact sets. 
\end{theorem}

\paragraph{Proof.} Let us look at a decorated elementary chora. 

\vspace{.5cm}

\centerline{\includegraphics[angle=270, width=0.5\textwidth]{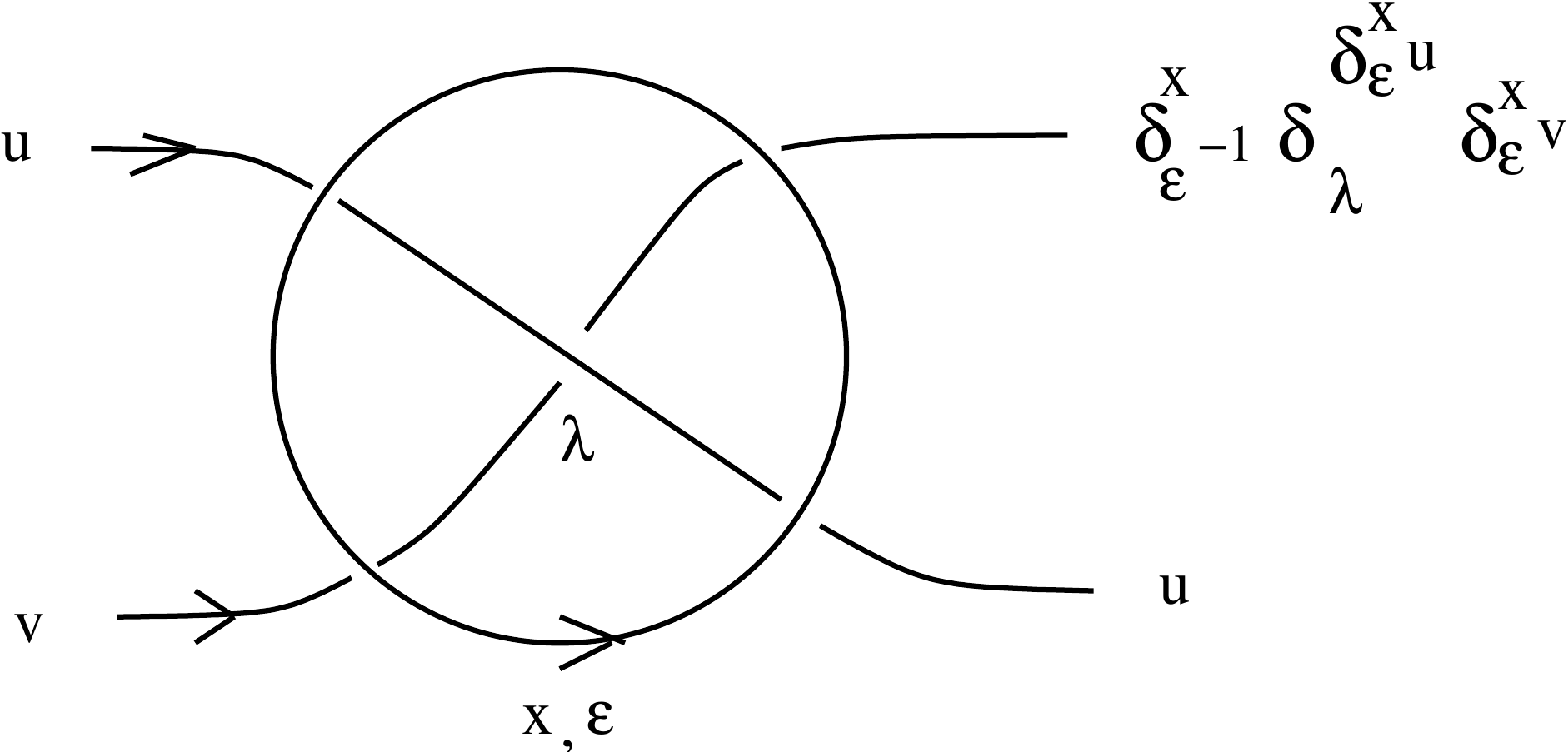}}
%\caption{Carrier inside chora. It decomposes as four carriers.}
%\label{carrierinchora}

\vspace{.5cm}

Because it decomposes into difference and sum gates which converge as
$\varepsilon$ goes to zero, in follows that the (input-output function
associated to the)  elementary chora diagram converges,  uniformly with 
respect to $x,u,v$ in compact sets, to a (function associated to) the gate: 

\vspace{.5cm}

\centerline{\includegraphics[angle=270, width=0.5\textwidth]{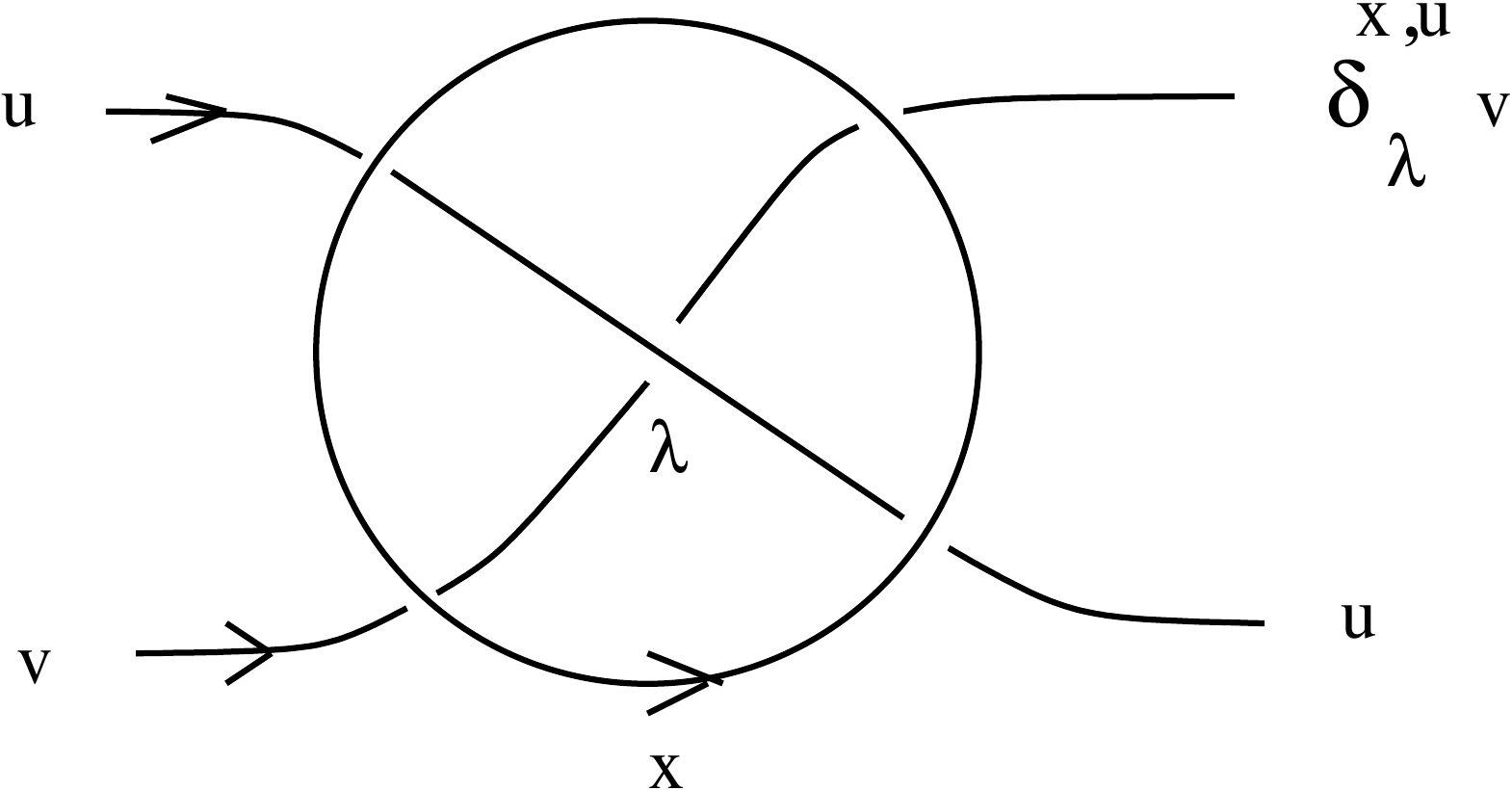}}
%\caption{Carrier inside chora. It decomposes as four carriers.}
%\label{carrierinchora}

\vspace{.5cm}

Here the function $\displaystyle (\lambda,u,v) \mapsto \delta^{x,u}_{\lambda}v$ is 
the dilation structure in the tangent space at $x \in X$. Equivalently, is the 
limit, as $\varepsilon$ goes to zero, of the transport of the original dilation
structure by the "map" $\displaystyle \delta^{x}_{\varepsilon^{-1}}$. 
This is a linear
dilation structure, fact which "explains" why Reidemeister  III move can be done 
in the limit (that is, for tangles decorated with linear dilation structures). 

I shall exploit the fact that the convergence is uniform. A written proof with 
lots of symbols can be easily written, but instead I draw  the residue diagram 
in the following equivalent form (using the chora-inside-chora decomposition 
and some Reidemeister II moves):

\vspace{.5cm}

\centerline{\includegraphics[angle=270, width=0.7\textwidth]{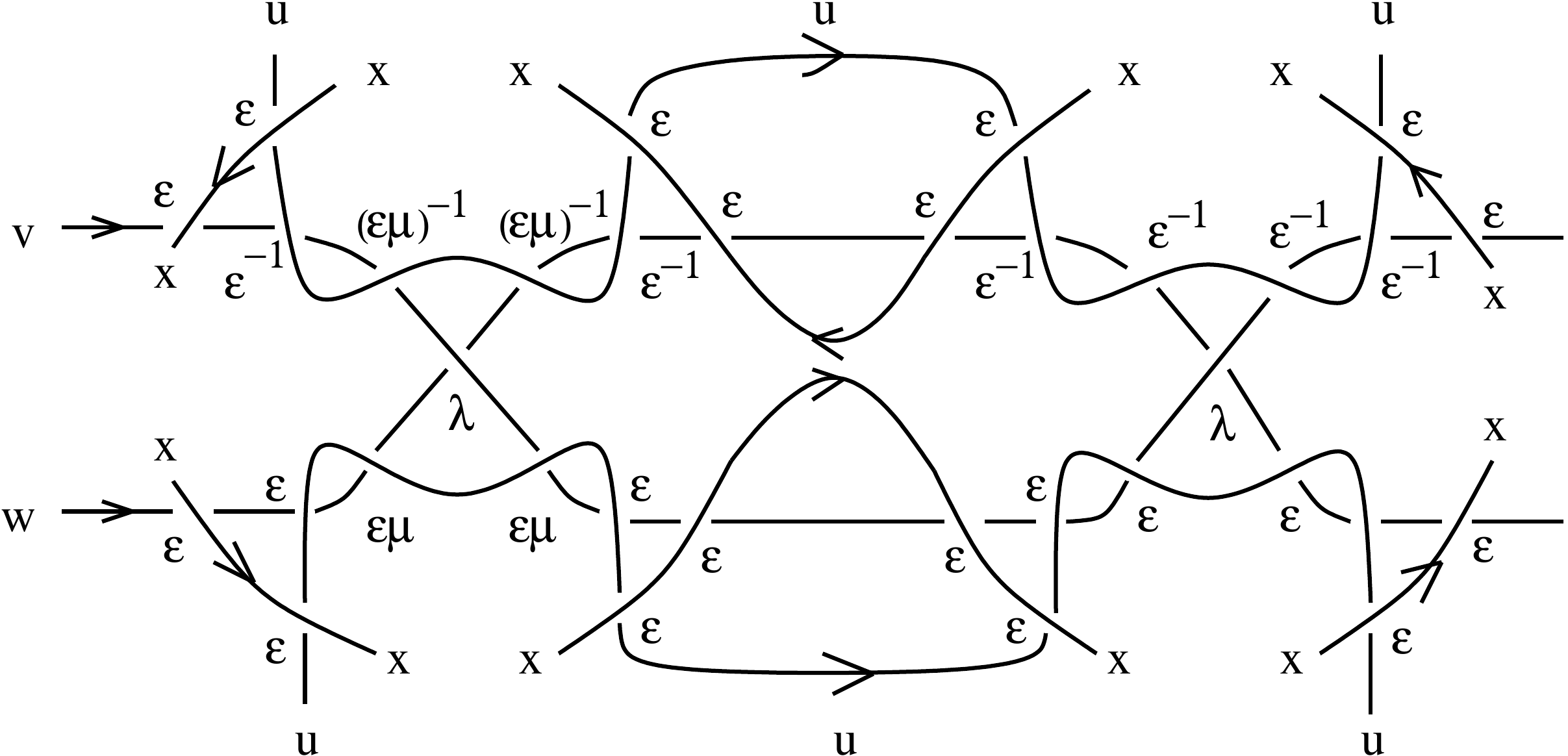}}
%\caption{Carrier inside chora. It decomposes as four carriers.}
%\label{carrierinchora}

\vspace{.5cm}

What we see is that inside the $x,\varepsilon$ chora we compose two elementary 
choroi:  the 
transport  of the original dilation structure by $\displaystyle
\delta^{\delta^{x}_{\varepsilon} u}_{(\varepsilon \mu)^{-1}}$  (at the left) 
with the inverse of the transport of the original dilation structure 
by $\displaystyle \delta^{\delta^{x}_{\varepsilon} u}_{\varepsilon^{-1}}$ (at
the right).  
Otherwise said, we compare the map of the territory at scale $\varepsilon \mu$
with the one at scale $\varepsilon$, as $\varepsilon$ goes to zero. 
I explained previously that the axiom A4 of dilation structures (or the last 
axiom of uniform $\Gamma$-irqs, more generally) is an expression of the 
stability of the "ideal foveal maps" provided by the dilations, which translates
into the fact that this composition converges to the identity map. There is 
though a supplementary observation to make: the basepoint of these maps is not 
$x$, but $\displaystyle \delta^{x}_{\varepsilon} u$, which moves, as
$\varepsilon$ goes to zero, to $x$. Here is where I use the uniformity of the
convergence, which implies that the convergence to the identity still holds,
even for this situation where the basepoint itself converges to $x$. \hfill
$\square$

\section{Appendix I: From maps to dilation structures}

Imagine that  the metric space $(X,d)$ represents a territory. We want to make 
maps  of $(X,d)$ in the metric space $(Y,D)$ (a piece of paper, or a
scaled model). 

In fact, in order to understand the territory $(X,d)$, we need many maps, at
many scales. For any point $x \in X$ and any scale $\varepsilon > 0$ we shall 
make a map of a neighbourhood of $x$, ideally. In practice, a good knowledge 
of a territory amounts to have, for each of several scales $\displaystyle
\varepsilon_{1} > \varepsilon_{2} > ... > \varepsilon_{n}$ an atlas of maps 
of overlapping parts of  $X$ (which together form a cover of the territory $X$). 
All the maps from all the atlasses have to be compatible one with another. 

The ideal model of such a body of knowledge is embodied into the notion of 
a manifold. To have $X$ as a manifold over the model space $Y$ means exactly 
this. 

Examples from metric geometry (like sub-riemannian spaces) show that the 
manifold idea could be too rigid in some situations. We shall replace it 
with the idea of a dilation structure. 

We shall see that a dilation structure (the right generalization of a 
smooth space, like a manifold), represents an idealization 
of the more realistic situation of having at our disposal many maps, at many
scales, of the territory, with the property that the accuracy, precision and 
resolution of such maps, and of relative maps deduced from them, are controlled 
by the scale (as the scale goes to zero, to infinitesimal). 

There are two facts which I need to stress. First is that such a generalization
is necessary. Indeed, by looking at the large gallery of metric spaces which we
now know, the metric spaces with a manifold structure form a tiny and very very
particular class. Second is that we tend to take for granted the body of
knowledge represented by a manifold structure (or by a dilation structure). 
Think as an example at the manifold structure of the Earth. It is an
idealization of the collection of all cartographic maps of parts of the Earth. 
This is a huge data basis and it required a huge amount of time and energy in
order to be constructed. To know, understand the territory is a huge task, 
largely neglected. We "have" a manifold, "let $X$ be a manifold". And even if 
we do not doubt that the physical space (whatever that means) is a boring 
$\displaystyle \mathbb{R}^{3}$, it is nevertheless another task to determine 
with the best accuracy possible a certain point in that physical space, based on
the knowledge of the coordinates. For example GPS costs money and time to build 
and use. Or, it is rather easy to collide protons, but to understand and keep 
the territory fixed (more or less) with respect to the map, that is where most
of the effort goes. 

A model of such a map of $(X,d)$ in $(Y,D)$ is  a relation  
$\rho \subset X \times Y$,  a subset of 
a cartesian product $X \times Y$ of two sets. A particular type of relation 
is the graph of a function $\displaystyle f: X \rightarrow Y$, defined as the
relation 
$$\rho \, = \, \left\{ (x, f(x)) \mbox{ : } x \in X \right\}$$ 
but there are many relations which cannot be described as graphs of functions.

%\vspace{.5cm}

\centerline{\includegraphics[angle=270, width=0.5\textwidth]{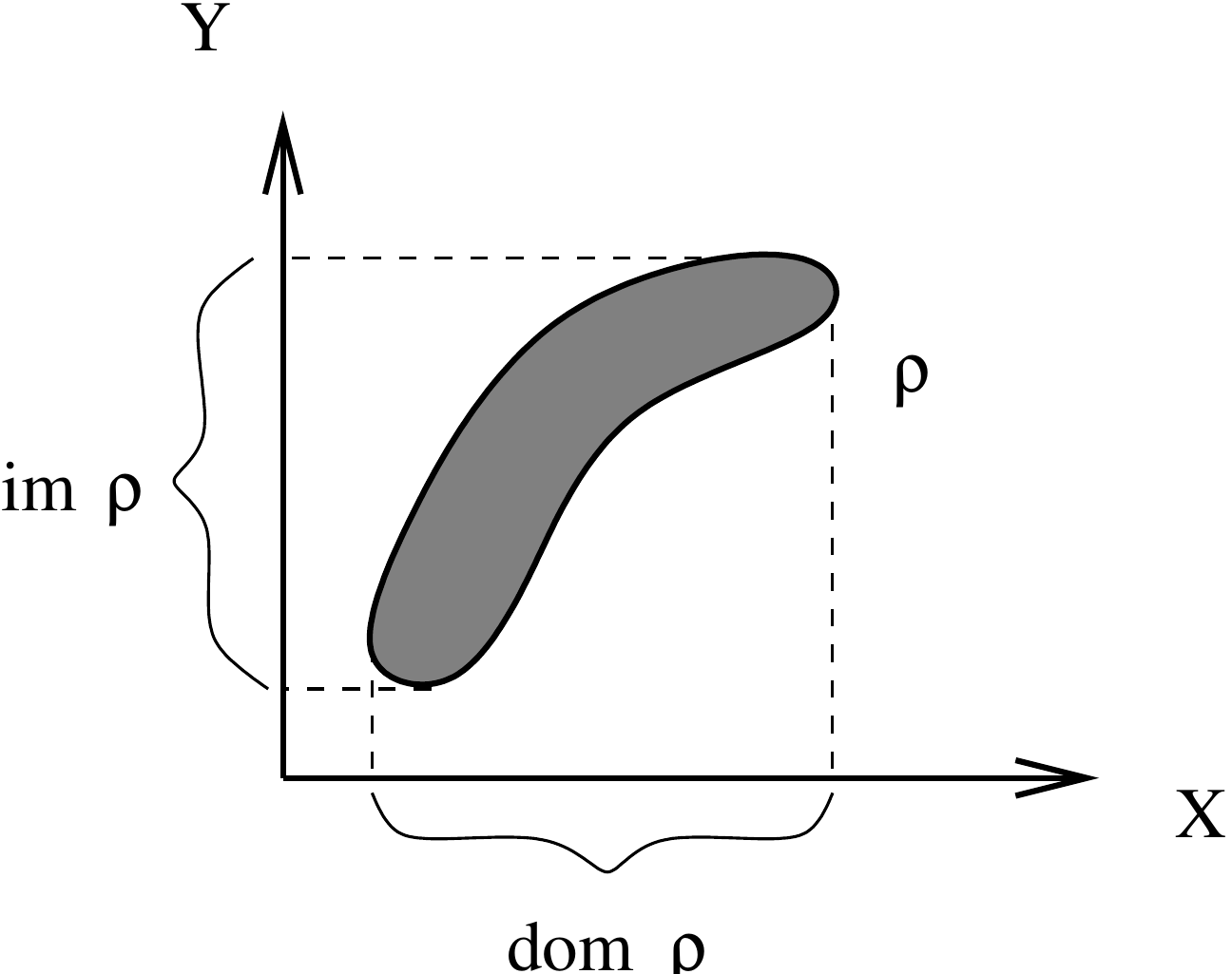}}
%\caption{}
%\label{fovea1}

%\vspace{.5cm}

Imagine that pairs $(u,u') \in \rho$ are pairs 

\vspace{.5cm}

\centerline{(point in the space $X$, pixel in the "map space" Y)} 

\vspace{.5cm}

\noindent with the meaning that 
the point $u$ in $X$ is represented as the pixel $u'$ in $Y$. 

I don't suppose that there is a one-to-one correspondence between points in 
$X$ and pixels in $Y$, for various reasons, for example: due to repeated measurements 
there is no unique way to associate pixel to a point, or a point to a pixel. 
The relation $\rho$ represents  the cloud of pairs point-pixel which are 
compatible with all measurements. 

I shall use this model of a map for simplicity reasons. A better,  
more realistic model could be one using probability measures, but this model is 
sufficient for the needs of this paper. 

For a given map $\rho$ the point $x \in X$ in the space $X$  is associated
 the set of points $\left\{ y \in Y \mbox{ : } (x,y) \in \rho \right\}$ in the 
 "map space" $Y$. Similarly, to the "pixel" $y \in Y$ in the "map space" $Y$
  is associated the set of points 
  $\left\{ x \in X \mbox{ : } (x,y) \in \rho \right\}$ in the space $X$. 

\vspace{.5cm}

\centerline{\includegraphics[angle=270, width=0.5\textwidth]{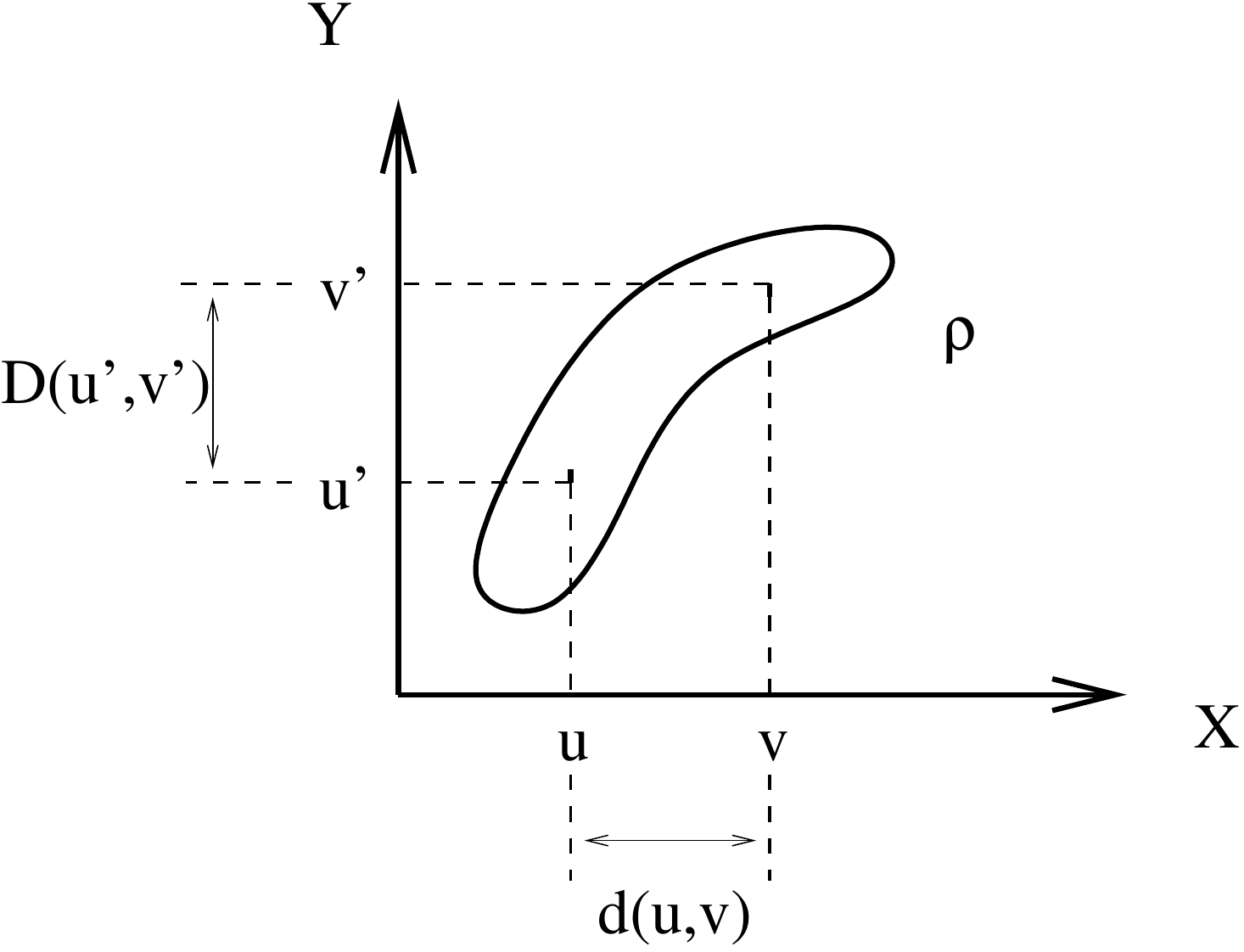}}
%\caption{}
%\label{fovea1}

\vspace{.5cm}

A good map is one which does not distort distances too much. Specifically,
considering any two points $u,v \in X$ and any two pixels $u',v' \in Y$ which 
represent these points, i.e. $(u,u'), (v,v') \in \rho$, the distortion of 
distances between these points is measured by the number

$$\mid d(u,v) - D(u',v') \mid$$

\subsection{Accuracy, precision, resolution, Gromov-Hausdorff distance}

\paragraph{Notations concerning relations.}  
Even if relations are more general than (graphs of) functions, there is no harm to
use, if needed, a functional notation. For any relation  
$\rho  \subset X \times Y$  we shall write $\rho(x) = y$ or $\displaystyle 
\rho^{-1}(y) = x$ if $(x,y) \in \rho$. Therefore we may have $\rho(x) = y$ and 
$\rho(x) = y'$ with $y \not = y'$, if $(x,y) \in f$ and $(x,y') \in f$. In some
drawings, relations will be figured by a large arrow, as shown further.

\vspace{.5cm}

\centerline{\includegraphics[angle=270, width=0.5\textwidth]{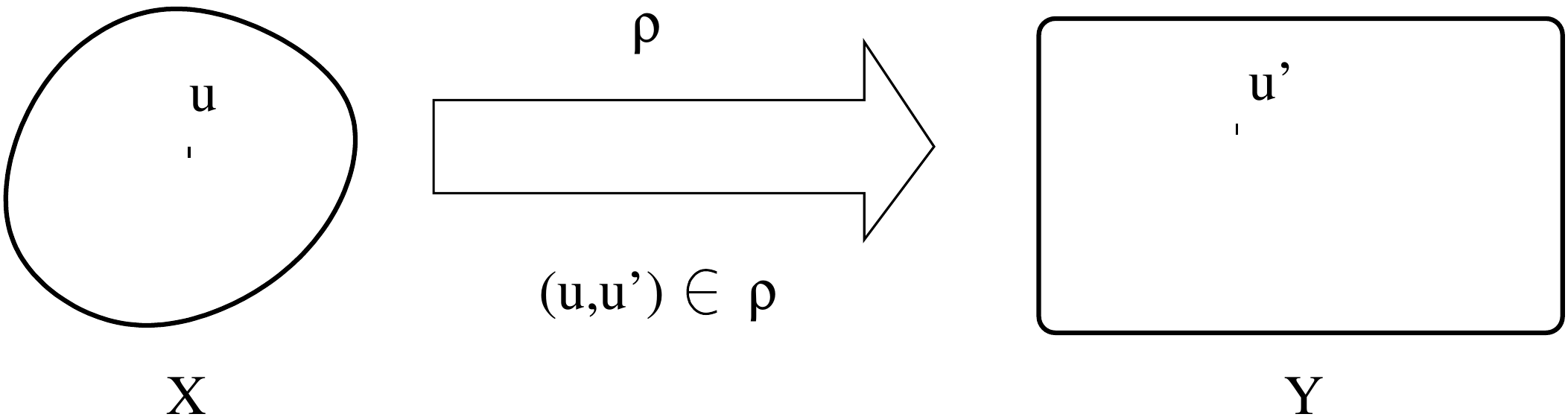}}
%\caption{}
%\label{fovea1}

\vspace{.5cm}

The domain of the relation $\rho$ is the set $dom \ \rho \,  \subset X$ such that for any $x \in \, dom \ \rho$ 
there is $y \in Y$ with $\rho(x) = y$. The image of $\rho$ is the set of $im \ \rho \  \subset Y$ such that for any $y \in \, im \ \rho$ there is $x \in X$ with $\rho(x) =  y$.  
By convention, when we write that a statement $R(f(x), f(y), ...)$ is true, we mean that $R(x',y', ...)$ is 
true for any choice of $x', y', ...$, such that $(x,x'), (y,y'), ... \in f$. 

The inverse of a relation $\rho \subset X \times Y$ is the relation 
$$\rho^{-1} \subset Y \times X \, , \quad \rho^{-1} \, = \, \left\{ (u',u) \mbox{
: } (u,u') \in \rho \right\}$$ 
and if $\rho' \subset X \times Y$, $\rho" \subset Y \times Z$ are two relations,
their composition is defined as 
$$\rho = \rho" \circ \rho' \subset X \times Z $$ 
$$\rho \, = \, \left\{ (u,u") \in X \times Z \mbox{ : } \exists u'\in Y \, (u,u')
\in \rho' \, (u',u") \in \rho" \right\}$$

I shall use the following convenient notation: by $\mathcal{O}(\varepsilon)$ we mean a positive function such that $\displaystyle \lim_{\varepsilon \rightarrow 0} \mathcal{O}(\varepsilon) = 0$.

\vspace{.5cm}

%\subsection{Accuracy, precision, resolution} 

In metrology, by definition, accuracy is \cite{metrology} 2.13 (3.5) "closeness of 
agreement between a measured quantity value and a true quantity value of a
measurand". (Measurement) precision is \cite{metrology} 2.15 "closeness of agreement between 
indications or measured quantity values obtained by replicate
measurements on the same or similar objects under specified conditions".
Resolution  is \cite{metrology} 2.15 "smallest change in a quantity being measured that
causes a perceptible change in the corresponding
indication".

For our model of a map, if $(u,u') \in \rho$ then $u'$ represent the measurement
of $u$. Moreover, because we see a map as a relation, the definition of 
the resolution can be restated as the supremum of distances between points 
in $X$ which are represented by the same pixel. Indeed, if the distance between 
two points in $X$ is bigger than this supremum then they cannot be represented 
by the same pixel.

\begin{definition}
Let $\rho \subset X \times Y$ be a relation which represents a map of 
$dom \ \rho \, \subset (X,d)$ into $ im \ \rho \, \subset (Y,D)$. To this map 
are associated three quantities: accuracy, precision and resolution. 

 The accuracy of $\rho$ is defined by: 
\begin{equation}
acc(\rho) \, = \, \sup \left\{ \mid D(y_{1}, y_{2}) - d(x_{1},x_{2}) \mid \mbox{
: } (x_{1},y_{1}) \in \rho \, , \, (x_{2},y_{2}) \in \rho \right\}
\label{acc1}
\end{equation}
 The resolution of $\rho$ at  $y \in im \ \rho$ is 
\begin{equation}
res(\rho)(y) \, = \, \sup \left\{ d(x_{1},x_{2}) \mbox{
: } (x_{1},y) \in \rho \, , \, (x_{2},y) \in \rho \right\}
\label{resy1}
\end{equation}
and the resolution of $\rho$ is given by: 
\begin{equation}
res(\rho) \, = \, \sup \left\{ res(\rho)(y) \mbox{
: } y \in \, im \ \rho \right\}
\label{res1}
\end{equation}
 The precision of $\rho$ at $x \in dom \ \rho$ is 
\begin{equation}
prec(\rho)(x) \, = \, \sup \left\{ D(y_{1},y_{2}) \mbox{
: } (x,y_{1}) \in \rho \, , \, (x,y_{2}) \in \rho \right\}
\label{precx1}
\end{equation}
and the precision of $\rho$ is given by: 
\begin{equation}
prec(\rho) \, = \, \sup \left\{ prec(\rho)(x) \mbox{
: } x \in \, dom \ \rho \right\}
\label{prec1}
\end{equation}
\label{defacc}
\end{definition}

After measuring (or using other means to deduce) the distances $d(x', x")$ between 
all pairs of points in $X$ (we may have several values for the distance $d(x',x")$), 
we try to represent the collection 
of these distances  in $(Y,D)$. 
When we make a map $\rho$ we are not really measuring the distances between 
all points in $X$, then representing them as accurately as possible in $Y$.

What we do is that  we consider a relation $\rho$, with domain $M = \,
dom(\rho)$ 
 which is $\varepsilon$-dense in $(X,d)$, then we perform a " cartographic generalization"\footnote{\url{http://en.wikipedia.org/wiki/Cartographic_generalization},
"Cartographic generalization is the method whereby information is selected and 
represented on a map in a way that adapts to the scale of the display medium 
of the map, not necessarily preserving all intricate geographical or other 
cartographic details.} of the relation $\rho$ to a relation $\displaystyle \bar{\rho}$, a
 map of $(X,d)$ in $(Y,D)$,  for example as in the following definition. 

\begin{definition}
A subset $M \subset X$ of a metric space $(X,d)$ is $\varepsilon$-dense in 
$X$ if for any $u \in X$ there is $x \in M$ such that $d(x,u) \leq \varepsilon$. 

Let $\rho \subset X \times Y$ be a relation such that $dom \ \rho$ is 
$\varepsilon$-dense in $(X,d)$ and $im \ \rho$  is $\mu$-dense in 
 $(Y,D)$. We define then $\displaystyle \bar{\rho} \subset X \times Y$ by:
 $\displaystyle (x,y) \in \bar{\rho}$ if there is $\displaystyle (x',y') \in
 \rho$ such that $\displaystyle d(x,x')\leq \varepsilon$ and $\displaystyle 
 D(y,y') \leq \mu$.
 \label{defgencart1}
\end{definition}

If $\rho$ is a relation as described in definition \ref{defgencart1} then 
accuracy $acc(\rho)$, $\varepsilon$ and $\mu$ control the precision 
$prec(\rho)$ and resolution $res(\rho)$. Moreover, the accuracy, precision and 
resolution of $\displaystyle \bar{\rho}$ are controlled by those of 
$\rho$ and $\varepsilon$, $\mu$, as well. This is explained in the next
proposition.

\begin{proposition}
Let $\rho$ and $\displaystyle \bar{\rho}$ be as described in definition
\ref{defgencart1}. Then: 
\begin{enumerate}
\item[(a)] $\displaystyle res(\rho) \, \leq \, acc(\rho)$, 
\item[(b)] $\displaystyle prec(\rho) \, \leq \, acc(\rho)$, 
\item[(c)] $\displaystyle res(\rho) + 2 \varepsilon \leq \, res(\bar{\rho}) \leq
\, acc(\rho) + 2(\varepsilon + \mu)$,
\item[(d)] $\displaystyle prec(\rho) + 2 \mu \leq \, prec(\bar{\rho}) \leq
\, acc(\rho) + 2(\varepsilon + \mu)$,
\item[(e)] $\displaystyle \mid acc(\bar{\rho}) - \, acc(\rho) \mid \leq
2(\varepsilon + \mu)$. 
\end{enumerate}
\label{propacc1}
\end{proposition}

\paragraph{Proof.} Remark that (a), (b) are immediate consequences of definition
\ref{defacc} and that 
(c) and (d) must have identical proofs, just by switching $\varepsilon$ with 
$\mu$ and $X$ with $Y$ respectively. I shall prove therefore (c) and (e). 

For proving (c), consider $y \in Y$. By definition of 
$\displaystyle \bar{\rho}$ we write 
$$\left\{ x \in X \mbox{ : } (x,y) \in \bar{\rho} \right\} \, = \, \bigcup_{(x',y')
\in \rho , y' \in \bar{B}(y,\mu)} \bar{B}(x',\varepsilon)$$
Therefore we get 
$$res(\bar{\rho})(y) \, \geq \, 2 \varepsilon + \sup \left\{ res(\rho)(y') 
\mbox{ : } y' \in \,
im(\rho) \cap \bar{B}(y,\mu) \right\} $$
By taking the supremum over all $y \in Y$ we obtain the inequality 
$$res(\rho) + 2 \varepsilon \leq \, res(\bar{\rho})$$
For the other inequality, let us consider $\displaystyle (x_{1},y), (x_{2},y)
\in \bar{\rho}$ and $\displaystyle (x_{1}', y_{1}'), (x_{2}', y_{2}') \in \rho$
such that $\displaystyle d(x_{1},x_{1}') \leq \varepsilon, d(x_{2},x_{2}') 
\leq \varepsilon, D(y_{1}',y) \leq \mu,  D(y_{2}',y) \leq \mu$. Then: 
$$d(x_{1},x_{2}) \leq 2 \varepsilon + d(x_{1}',x_{2}') \leq 2 \varepsilon + \, 
acc(\rho) + d(y_{1}',y_{2}') \leq 2 (\varepsilon + \mu) + \, 
acc(\rho)$$
Take now a supremum and arrive to the desired inequality. 

For the proof of (e)  let us consider for $i = 1,2$   
$\displaystyle (x_{i},y_{i}) \in \bar{\rho}, (x_{i}',y_{i}') \in \rho$ such 
that $\displaystyle d(x_{i}, x_{i}') \leq \varepsilon, D(y_{i},y_{i}') \leq
\mu$.   It is then  enough to take absolute values and transform 
the following equality  
$$d(x_{1},x_{2}) - D(y_{1},y_{2}) = d(x_{1},x_{2}) - d(x_{1}',x_{2}') + 
d(x_{1}',x_{2}') - D(y_{1}',y_{2}') + $$ 
$$+ D(y_{1}',y_{2}') - D(y_{1},y_{2})$$ 
into well chosen, but straightforward, inequalities. \hfill $\square$

\vspace{.5cm}

The  following definition of the  Gromov-Hausdorff distance for metric spaces is natural, 
owing to the fact that the accuracy (as defined in definition \ref{defacc}) controls the 
precision and resolution. 

\begin{definition}
Let $\displaystyle (X, d)$, $(Y,D)$,  be a pair of metric spaces and $\mu > 0$. 
We shall say that $\mu$ is admissible if  there is a relation 
$\displaystyle \rho \subset X \times Y$ such that  
$\displaystyle dom \ \rho = X$, $\displaystyle im \ \rho = Y$, and $acc(\rho) \leq \mu$.
The Gromov-Hausdorff distance  between $\displaystyle (X,d)$ and $\displaystyle  
(Y,D)$  is   the infimum of admissible numbers $\mu$. 
\label{defgh}
\end{definition}

As introduced in definition \ref{defgh}, the Gromov-Hausdorff  (GH) distance is not a true 
distance, because the GH distance between two isometric  metric spaces   
is equal to zero. In fact the GH distance induces a distance on isometry classes of 
compact metric spaces. 

The GH distance thus represents a lower bound on the accuracy of making maps of 
$(X,d)$ into $(Y,D)$. Surprising as it might seem, there are many examples of pairs of 
metric spaces with the property that the 
GH distance between any pair of closed balls from these spaces, 
considered with the distances properly rescaled, is greater than a strictly
positive number, independent of the choice of the balls. Simply put: {\it there 
are pairs of spaces $X$, $Y$ such that is impossible to make maps of parts 
of $X$ into $Y$ with arbitrarily small accuracy.} 

Any measurement is equivalent with making a map, say of $X$ (the territory 
of the phenomenon) into $Y$ (the map space of the laboratory). The possibility
that there might a physical difference (manifested as a strictly positive 
GH distance) between these two spaces, even if they both might be topologically 
the same (and with trivial topology, say of a $\displaystyle \mathbb{R}^{n}$),
 is ignored in physics, apparently. On one
side, there is no experimental way to confirm that a territory is the same 
at any scale (see the section dedicated to the notion of scale), but much of physical
explanations are based on differential calculus, which has as the most basic
assumption that locally and infinitesimally the territory is the same. 
On the other side the imposibility of making maps of the phase space of a 
quantum object into the macroscopic map space of the laboratory might be a 
manifestation of the fact that there is a difference (positive GH distance 
between maps of the territory realised with the help of physical phenomena) 
between "small" and "macroscopic" scale.

\subsection{Scale} Let $\varepsilon >0$. A map of $(X,d)$ into $(Y,D)$, at scale 
$\varepsilon$ is a map of $\displaystyle (X, \frac{1}{\varepsilon}d)$ into $(Y,D)$. 
Indeed, if this map would have accuracy equal to $0$ then a value of a 
distance between points in $X$ equal to  $\displaystyle L$ would correspond to a value 
of the distance between the corresponding points on the map in $(Y,D)$ equal to
$\varepsilon \, L$. 

In cartography, maps  of the same territory done at smaller and smaller scales (smaller 
and smaller $\varepsilon$) must have the property that, at the same resolution, 
the  accuracy and precision (as defined in definition \ref{defacc}) have to become 
smaller and smaller. 

In mathematics, this could serve as the definition of the metric tangent space to 
a point in $(X,d)$, as seen in $(Y,D)$. 

\begin{definition}
We say that $(Y,D, y)$ ($y \in Y$) represents the (pointed unit ball in the) metric 
tangent space at $x \in X$ of $(X,d)$ if there exist a pair formed by: 
\begin{enumerate}
\item[-] a "zoom sequence", that is a sequence 
$$\displaystyle \varepsilon \in (0,1] \, \mapsto \,  \rho_{\varepsilon}^{x} \subset
(\bar{B}(x,\varepsilon), \frac{1}{\varepsilon} d)  \times (Y,D)$$ 
such that $\displaystyle dom \ \rho_{\varepsilon}^{x} = \bar{B}(x,\varepsilon)$, 
$\displaystyle im \ \rho_{\varepsilon}^{x} = Y$, $\displaystyle (x,y) \in \rho_{\varepsilon}^{x}$ 
for any $\displaystyle \varepsilon \in (0,1]$ and 
\item[-] a "zoom modulus" $\displaystyle F: (0,1) \rightarrow [0,+\infty)$ such
that  $\displaystyle \lim_{\varepsilon \rightarrow 0} F(\varepsilon)  =  0$, 
\end{enumerate}
such that for all $\varepsilon \in (0,1)$  we have 
$\displaystyle acc(\rho_{\varepsilon}^{x}) \leq F(\varepsilon)$. 
\label{deftangsp}
\end{definition}

Using the notation proposed previously, we can write 
$\displaystyle F(\varepsilon) =  \mathcal{O}(\varepsilon)$, if there is no need 
to precisely specify a zoom modulus function.

Let us write again the definition of resolution, accuracy, precision, in the presence of
scale. The accuracy of $\displaystyle \rho_{\varepsilon}^{x}$ is defined by: 
\begin{equation}
acc(\rho_{\varepsilon}^{x}) \, = \, \sup \left\{ \mid D(y_{1}, y_{2}) - \frac{1}{\varepsilon} d(x_{1},x_{2}) \mid \mbox{
: } (x_{1},y_{1})  , \, (x_{2},y_{2}) \in \rho_{\varepsilon}^{x} \right\}
\label{acc2}
\end{equation}
 The resolution of $\displaystyle \rho_{\varepsilon}^{x}$ at  $z \in Y$ is 
\begin{equation}
res(\rho_{\varepsilon}^{x})(z) \, = \, \frac{1}{\varepsilon} \, \sup \left\{ d(x_{1},x_{2}) \mbox{
: } (x_{1},z) \in \rho_{\varepsilon}^{x} \, , \, (x_{2},z) \in \rho_{\varepsilon}^{x} \right\}
\label{resy2}
\end{equation}
and the resolution of $\displaystyle \rho_{\varepsilon}^{x}$ is given by: 
\begin{equation}
res(\rho_{\varepsilon}^{x}) \, = \,  \sup \left\{ res(\rho_{\varepsilon}^{x})(y) \mbox{
: } y \in Y \right\}
\label{res2}
\end{equation}
 The precision of $\rho_{\varepsilon}^{x}$ at $\displaystyle u \in
\bar{B}(x,\varepsilon)$ is 
\begin{equation}
prec(\rho_{\varepsilon}^{x})(u) \, = \, \sup \left\{ D(y_{1},y_{2}) \mbox{
: } (u,y_{1}) \in \rho_{\varepsilon}^{x} \, , \, (u,y_{2}) \in \rho_{\varepsilon}^{x} \right\}
\label{precx2}
\end{equation}
and the precision of $\displaystyle \rho_{\varepsilon}^{x}$ is given by: 
\begin{equation}
prec(\rho_{\varepsilon}^{x}) \, = \, \sup \left\{ prec(\rho_{\varepsilon}^{x})(u) \mbox{
: } u \in \bar{B}(x,\varepsilon) \right\}
\label{prec2}
\end{equation}

If $(Y,D, y)$  represents the (pointed unit ball in the) metric 
tangent space at $x \in X$ of $(X,d)$ and $\displaystyle \rho_{\varepsilon}^{x}$ is the
sequence of maps at smaller and smaller scale, then we have: 

\begin{equation}
 \sup \left\{ \mid D(y_{1}, y_{2}) - \frac{1}{\varepsilon} d(x_{1},x_{2}) \mid \mbox{
: } (x_{1},y_{1}) , \, (x_{2},y_{2}) \in \rho_{\varepsilon}^{x}
\right\} \, = \, \mathcal{O}(\varepsilon)
\label{acclim}
\end{equation}

\begin{equation}
 \sup \left\{  D(y_{1}, y_{2})  \mbox{
: } (u,y_{1}) \in \rho_{\varepsilon}^{x} \, , \, (u,y_{2}) \in \rho_{\varepsilon}^{x} \, ,
\,  u \in \bar{B}(x,\varepsilon)
\right\} \, = \, \mathcal{O}(\varepsilon)
\label{preclim}
\end{equation}

\begin{equation}
  \sup \left\{  d(x_{1},x_{2}) \mbox{
: } (x_{1},z) \in \rho_{\varepsilon}^{x} \, , \, (x_{2},z) \in \rho_{\varepsilon}^{x} \, , \, z
\in Y 
\right\} \, = \, \varepsilon \, \mathcal{O}(\varepsilon)
\label{reslim}
\end{equation}

Of course, relation (\ref{acclim}) implies the other two, but it is interesting to 
notice the mechanism of rescaling.

\subsection{Scale stability. Viewpoint stability}

\vspace{.5cm}

\centerline{\includegraphics[angle=270, width=0.7\textwidth]{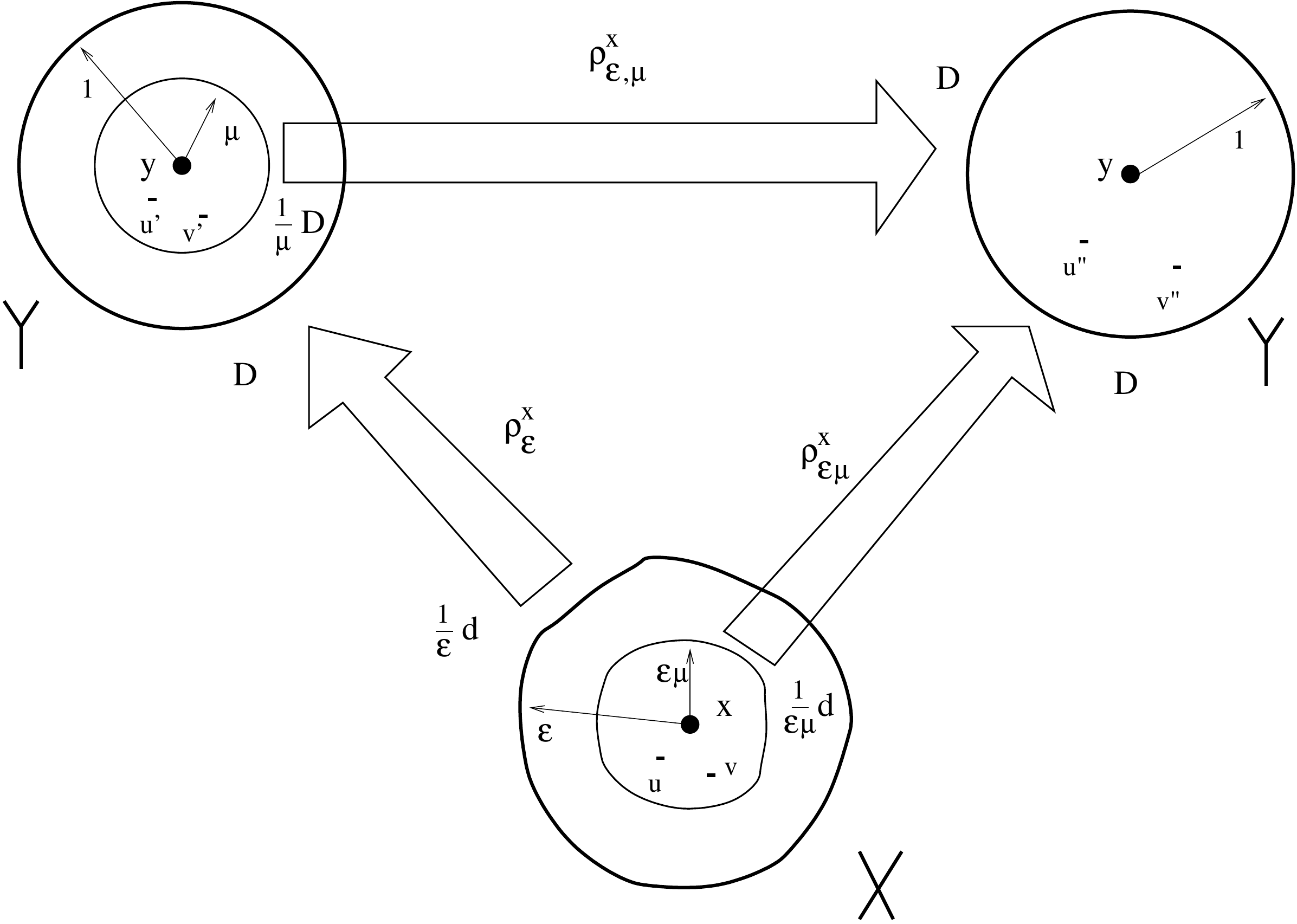}}
%\caption{}
%\label{fovea1}

\vspace{.5cm}

I shall suppose further that there is a metric tangent space at $x \in X$ and I
shall work with a zoom sequence of maps described in definition \ref{deftangsp}. 

Let $\varepsilon, \mu \in (0,1)$ be two scales. Suppose we have the maps 
of the territory $X$, around $x \in X$, 
at scales $\varepsilon$ and $\varepsilon \mu$, 
$$\displaystyle 
\rho_{\varepsilon}^{x} \subset \bar{B}(x,\varepsilon) \times \bar{B}(y,1)$$ 
$$\displaystyle \rho_{\varepsilon \mu}^{x} \subset \bar{B}(x,\varepsilon 
\mu) \times \bar{B}(y,1)$$
made into the tangent space at $x$, $\displaystyle (\bar{B}(y,1), D)$. 
The ball $\displaystyle \bar{B}(x, \varepsilon \mu) \subset X$ has then 
two maps. These maps are at different scales: the first is done at scale $\varepsilon$, 
the second is done at scale $\varepsilon \mu$. 

What are the differences between these two maps? We could find out by defining
a new map 
\begin{equation}
\rho^{x}_{\varepsilon,\mu} \, = \, \left\{ (u',u") \in 
\bar{B}(y,\mu) \times \bar{B}(y,1) \mbox{ : } \right. 
\label{relchart}
\end{equation}
$$\left. \exists u \in \bar{B}(x,\varepsilon \mu) \, (u,u') \in \rho_{\varepsilon}^{x} \,
, \, (u, u") \in \rho^{x}_{\varepsilon \mu} \right\}$$ 
and measuring its accuracy, with respect to the distances $\displaystyle
\frac{1}{\mu} D$ (on the domain) and $\displaystyle D$ (on the image).

Let us consider $\displaystyle (u,u'), (v,v') \in \rho_{\varepsilon}^{x}$ and  
$\displaystyle (u,u"), (v,v") \in \rho_{\varepsilon \mu}^{x}$ such that 
$\displaystyle (u', u"), (v',v")  \in \rho^{x}_{\varepsilon,\mu}$. Then: 
$$\mid D(u",v") - \frac{1}{\mu} D(u',v') \mid \, \leq \, \mid 
\frac{1}{\mu} D(u',v') - \frac{1}{\varepsilon \mu} d(u,v) \mid \, + \, 
\mid \frac{1}{\varepsilon \mu} d(u,v) -  D(u",v") \mid$$ 
We have therefore an estimate for the accuracy of the map 
$\displaystyle \rho^{x}_{\varepsilon,\mu}$, coming from estimate (\ref{acclim})
applied for $\displaystyle \rho_{\varepsilon}^{x}$ and 
$\displaystyle \rho_{\varepsilon \mu}^{x}$: 
\begin{equation}
acc(\rho^{x}_{\varepsilon,\mu}) \, \leq \, \frac{1}{\mu} 
\mathcal{O}(\varepsilon) + \mathcal{O}(\varepsilon \mu)
\label{eqcascade}
\end{equation}
This explains the cascading of errors phenomenon, namely, for fixed $\mu$, as 
$\varepsilon$ goes to $0$ the accuracy of the map $\displaystyle
\rho^{x}_{\varepsilon,\mu}$ becomes smaller and smaller, meaning that the 
maps of the ball $\displaystyle \bar{B}(x, \varepsilon \mu) \subset X$ at the 
scales $\varepsilon, \varepsilon \mu$ (properly rescaled) are  more and more
alike. On the contrary, for fixed $\varepsilon$, as $\mu$ goes to $0$, the bound
on the accuracy becomes bigger and bigger, meaning that by using only the map 
at scale $\varepsilon$,  magnifications of  a smaller
scale region of this map may be less accurate than the map of this smaller 
region done at the smaller scale.

I shall add a supplementary hypothesis to the one 
concerning the existence of the metric tangent space. It is somehow natural to
suppose that as $\varepsilon$ converges to $0$ the map $\displaystyle
\rho^{x}_{\varepsilon,\mu}$ converges to a map $\displaystyle
\bar{\rho}^{x}_{\mu}$. This is described further. 

\vspace{.5cm}

\centerline{\includegraphics[angle=270, width=0.7\textwidth]{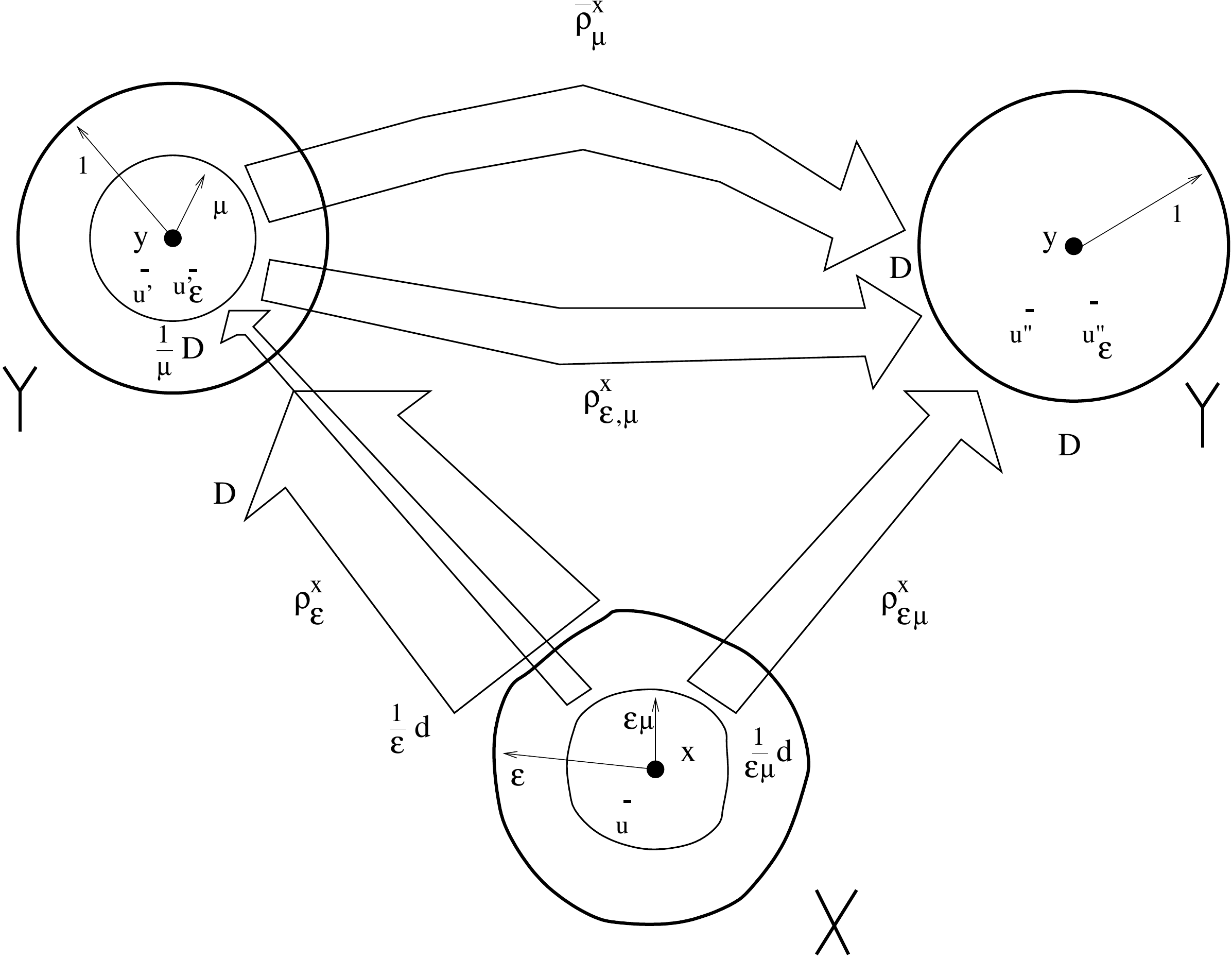}}
%\caption{}
%\label{fovea2}

\vspace{.5cm}

\begin{definition}
Let the zoom sequence $\displaystyle \rho^{x}_{\varepsilon}$ be as in definition
\ref{deftangsp} and for given $\mu \in (0,1)$, the map $\displaystyle
\rho^{x}_{\varepsilon, \mu}$ be defined as in (\ref{relchart}). We say that 
the zoom sequence $\displaystyle \rho^{x}_{\varepsilon}$ is  scale stable at 
scale $\mu$ if there is a relation $\displaystyle \bar{\rho}^{x}_{\mu} \subset 
\bar{B}(y,\mu) \times \bar{B}(y,1)$ such that the Haussorff distance between 
$\displaystyle \rho^{x}_{\varepsilon, \mu}$ and $\displaystyle
\bar{\rho}^{x}_{\mu}$, in the metric space $\displaystyle \bar{B}(y,\mu) \times
\bar{B}(y,1)$ 
with the distance 
$$D_{\mu} \left( (u',u"), (v',v") \right) \, = \, \frac{1}{\mu} D(u',v') +
D(u",v") $$
can be estimated as: 
$$D_{\mu}^{Hausdorff}\left( \rho^{x}_{\varepsilon, \mu} , \bar{\rho}^{x}_{\mu}
\right)  \,
\leq \, F_{\mu}(\varepsilon)$$ 
with $\displaystyle F_{\mu}(\varepsilon) = \mathcal{O}_{\mu}(\varepsilon)$.
Such a function $\displaystyle F_{\mu}(\cdot)$ is called a scale stability modulus of 
the zoom sequence  $\displaystyle \rho^{x}_{\varepsilon}$.
\label{defstable}
\end{definition}

This means that for any $\displaystyle (u',u") \in \bar{\rho}^{x}_{\mu}$ there
is a sequence $\displaystyle (u'_{\varepsilon}, u"_{\varepsilon}) \in
\rho^{x}_{\varepsilon, \mu}$ such that 
$$\lim_{\varepsilon \rightarrow 0} u'_{\varepsilon} \, = \, u' \quad
\lim_{\varepsilon \rightarrow 0} u"_{\varepsilon} \, = \, u"$$

\begin{proposition}
If there is a scale stable zoom sequence   $\displaystyle \rho^{x}_{\varepsilon}$ 
 as in definitions \ref{deftangsp} and \ref{defstable} then the space 
 $(Y,D)$ is self-similar in a neighbourhood of point $y \in Y$, namely for any 
 $\displaystyle (u',u"), (v',v") \in \bar{\rho}^{x}_{\mu}$ we have: 
 $$D(u",v") \, = \, \frac{1}{\mu} D(u',v')$$
 In particular $\displaystyle \bar{\rho}^{x}_{\mu}$ is the graph of a function 
 (the precision and resolution are respectively equal to $0$). 
 \label{pstable}
 \end{proposition}
 
 \paragraph{Proof.} 
 Indeed, for any $\varepsilon \in (0,1)$ let us consider $\displaystyle 
 (u'_{\varepsilon}, u"_{\varepsilon}), (v'_{\varepsilon},v"_{\varepsilon}) 
 \in \rho^{x}_{\varepsilon,\mu}$ such that 
 $$\frac{1}{\mu} D(u',u'_{\varepsilon}) + D(u",u"_{\varepsilon}) \, \leq \,
 \mathcal{O}_{\mu}(\varepsilon)$$  
 $$\frac{1}{\mu} D(v',v'_{\varepsilon}) + D(v",v"_{\varepsilon}) \, \leq \,
 \mathcal{O}_{\mu}(\varepsilon)$$  
 Then we get the following inequality, using also the cascading of errors
 inequality (\ref{eqcascade}),  
 $$\mid D(u",v") - \frac{1}{\mu} D(u',v') \mid \, \leq \, 2
 \mathcal{O}_{\mu}(\varepsilon) + \frac{1}{\mu} 
\mathcal{O}(\varepsilon) + \mathcal{O}(\varepsilon \mu)$$ 
We pass with $\varepsilon$ to $0$ in order to obtain the conclusion. \hfill
$\square$

\vspace{.5cm}

Instead of changing the scale (i.e. understanding the scale stability of 
the zoom sequence), we could explore what happens when we change the point of
view. 

\vspace{.5cm}

\centerline{\includegraphics[angle=270, width=0.7\textwidth]{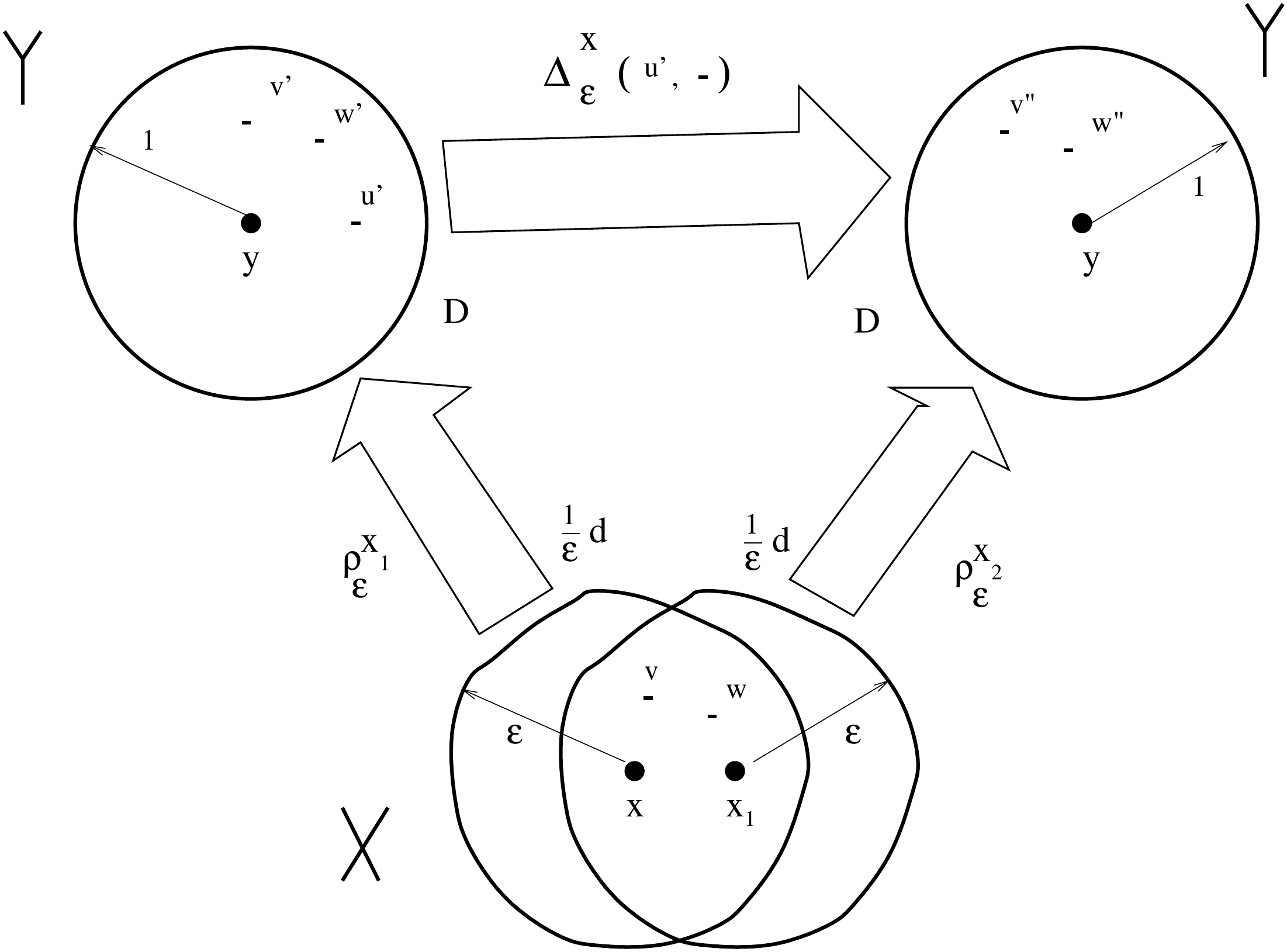}}
%\caption{}
%\label{fovea1}

\vspace{.5cm}

This time we have a zoom sequence, a scale $\varepsilon \in (0,1)$ and two 
points: $\displaystyle x \in X$ and $\displaystyle u' \in \bar{B}(y,1)$. To the 
point $u'$ from the map space $Y$ corresponds a point $\displaystyle 
x_{1} \in \bar{B}(x,\varepsilon)$ such that 
$$(x_{1}, u') \in \rho_{\varepsilon}^{x}$$

The points $\displaystyle x, x_{1}$ are neighbours, in the sense that 
$\displaystyle d(x,x_{1}) < \varepsilon$. The points of $X$ which are in the 
intersection 
$$\bar{B}(x,\varepsilon) \cap \bar{B}(x_{1},\varepsilon)$$ 
are represented by both maps, $\displaystyle \rho^{x}_{\varepsilon}$ and 
$\displaystyle \rho^{x_{1}}_{\varepsilon}$. These maps are different; the
relative map between them is defined as: 
\begin{equation}
\Delta^{x}_{\varepsilon}(u', \cdot) \, = \, 
\left\{ (v',v") \in \bar{B}(y,1) \mbox{ : } \exists v \in
\bar{B}(x,\varepsilon) \cap \bar{B}(x_{1},\varepsilon) \right.
\label{diffrel}
\end{equation}
$$ \left. \quad \quad \, 
(v,v') \in \rho^{x}_{\varepsilon} \, , \, (v,v") \in 
\rho^{x_{1}}_{\varepsilon} \right\}$$
and it is called "difference at scale $\varepsilon$, from $\displaystyle x$  
to $\displaystyle x_{1}$, as seen from $\displaystyle u'$". 

The viewpoint stability of the zoom sequence is expressed as the scale
stability: the zoom sequence is stable if the difference at scale 
$\varepsilon$ converges in the sense of Hausdorff distance, as $\varepsilon$
goes to $0$. 

\begin{definition}
Let the zoom sequence $\displaystyle \rho^{x}_{\varepsilon}$ be as in definition
\ref{deftangsp} and for any $\displaystyle u' \in \bar{B}(y,1)$, 
the map $\displaystyle \Delta^{x}_{\varepsilon}(u', \cdot)$ be defined as 
in (\ref{diffrel}). The zoom sequence $\displaystyle \rho^{x}_{\varepsilon}$ 
is  viewpoint stable if there is a relation $\displaystyle \Delta^{x}(u', \cdot) \subset 
\bar{B}(y,1) \times \bar{B}(y,1)$ such that the Haussorff distance
can be estimated as: 
$$D_{\mu}^{Hausdorff}\left( \Delta^{x}_{\varepsilon}(u', \cdot) , 
\Delta^{x}(u', \cdot)
\right)  \,
\leq \, F_{diff}(\varepsilon)$$ 
with $\displaystyle F_{diff}(\varepsilon) = \mathcal{O}(\varepsilon)$.
Such a function $\displaystyle F_{diff}(\cdot)$ is called a viewpoint stability modulus of 
the zoom sequence  $\displaystyle \rho^{x}_{\varepsilon}$.
\label{defstablev}
\end{definition}

There is a proposition analoguous with proposition \ref{pstable}, stating that 
the difference relation $\displaystyle \Delta^{x}(u', \cdot)$ is the graph of 
an isometry of $(Y,D)$. 

\subsection{Foveal maps} 
The following proposition shows that if we have a scale stable zoom sequence of maps  $\displaystyle \rho^{x}_{\varepsilon}$ 
 as in definitions \ref{deftangsp} and \ref{defstable} then we can improve every
 member of the sequence such that all maps from the new zoom sequence have better
 accuracy near the "center" of the map $x \in X$, which justifies the name 
 "foveal maps".
 
\begin{definition}
Let   $\displaystyle \rho^{x}_{\varepsilon}$ be a scale stable zoom sequence. We define for any $\varepsilon \in (0,1)$ the 
 $\mu$-foveal    map $\displaystyle \phi^{x}_{\varepsilon}$  made of all pairs 
$\displaystyle (u,u') \in  \bar{B}(x,\varepsilon) \times \bar{B}(y,1)$ such 
that
\begin{enumerate}
\item[-] if $\displaystyle u \in \bar{B}(x,\varepsilon \mu) $ then
$\displaystyle (u, \bar{\rho}_{\mu}^{x}(u')) \in \rho_{\varepsilon \mu}^{x}$, 
\item[-] or else $\displaystyle (u,u') \in \rho^{x}_{\varepsilon}$.
\end{enumerate}
\label{deffoveal}
\end{definition}

\begin{proposition}
Let   $\displaystyle \rho^{x}_{\varepsilon}$ be a scale stable zoom sequence with
associated zoom modulus $F(\cdot)$ and scale stability modulus $\displaystyle
F_{\mu}(\cdot)$. The sequence of $\mu$-foveal maps $\displaystyle
\phi^{x}_{\varepsilon}$  is then a scale stable 
zoom sequence with zoom modulus $\displaystyle F(\cdot) + \mu F_{\mu}(\cdot)$. 
Moreover, the accuracy of the restricted foveal map $\displaystyle
\phi^{x}_{\varepsilon} \cap \left(\bar{B}(x,\varepsilon \mu) \times 
\bar{B}(y,\mu)\right)$ 
is bounded by $\mu F(\varepsilon \mu)$, therefore the right hand side term in
the cascading of errors inequality (\ref{eqcascade}), applied for the restricted
foveal map, can be improved to $2 F(\varepsilon \mu)$.
\label{pfoveal}
\end{proposition}

\paragraph{Proof.}
Let $\displaystyle u \in \bar{B}(x,\varepsilon \mu)$. Then there are 
$\displaystyle u', u'_{\varepsilon} \in \bar{B}(y,\mu)$ and 
$\displaystyle u",u"_{\varepsilon} \in \bar{B}(y,1)$ such that 
$\displaystyle (u,u')\in \phi^{x}_{\varepsilon}$, $\displaystyle 
(u,u") \in \rho^{x}_{\varepsilon} \mu)$, $\displaystyle (u',u") 
\in \bar{\rho}^{x}_{\mu}$, $\displaystyle (u'_{\varepsilon},u"_{\varepsilon})
\in \rho^{x}_{\varepsilon,\mu}$ and 
$$ \frac{1}{\mu} D(u',u'_{\varepsilon}) + D(u",u"_{\varepsilon}) \, \leq \,
F_{\mu}(\varepsilon)$$

Let $\displaystyle u,v \in \bar{B}(x,\varepsilon \mu)$ and $u',v' \in 
\bar{B}(y,\mu)$ such that $\displaystyle (u,u'), (v,v') \in
\phi^{x}_{\varepsilon}$. According to the definition of 
$\displaystyle \phi^{x}_{\varepsilon}$, it follows that there are 
uniquely defined $u",v" \in \bar{B}(y,1)$ such that $\displaystyle 
(u,u"), (v,v") \in \rho^{x}_{\varepsilon \mu}$ and $\displaystyle 
(u', u"), (v',v") \in \bar{\rho}^{x}_{\mu}$. We then have: 
$$\mid \frac{1}{\varepsilon} d(u,v) - D(u',v') \mid \, = \, $$
$$= \, \mid \frac{1}{\varepsilon} d(u,v) - \mu D(u",v") \mid \, = \,$$
$$= \, \mu \mid \frac{1}{\varepsilon \mu} d(u,v) - D(u",v") \mid \, \leq \, 
\mu F(\varepsilon \mu)$$
Thus we proved that the accuracy of the restricted foveal map 
$$\displaystyle
\phi^{x}_{\varepsilon} \cap \left(\bar{B}(x,\varepsilon \mu) \times
\bar{B}(y,\mu)\right)$$
is bounded by $\mu F(\varepsilon \mu)$: 
\begin{equation}
\mid \frac{1}{\varepsilon} d(u,v) - D(u',v') \mid \, \leq \, 
\mu F(\varepsilon \mu)
\label{nd1}
\end{equation}
 If $\displaystyle u, v \in \bar{B}(x,\varepsilon) \setminus \bar{B}(x,\mu)$
  and  $\displaystyle (u,u'), (v,v') \in \phi^{x}_{\varepsilon}$ then  $\displaystyle (u,u'), (v,v') \in
\rho^{x}_{\varepsilon}$, therefore 
$$\mid \frac{1}{\varepsilon} d(u,v) - D(u',v') \mid \, \leq \, F(\varepsilon)$$
Suppose now that $\displaystyle (u,u'), (v,v') \in
\phi^{x}_{\varepsilon}$ and $\displaystyle u \in \bar{B}(x,\varepsilon \mu)$ but
$\displaystyle v 
\in \bar{B}(x,\varepsilon) \setminus \bar{B}(x,\mu)$. 
 We have then: 
$$\mid \frac{1}{\varepsilon} d(u,v) - D(u',v') \mid \, \leq \, $$
$$\leq \, \mid \frac{1}{\varepsilon} d(u,v) - D(u'_{\varepsilon},v') \mid + 
D(u',u'_{\varepsilon}) \, \leq \, F(\varepsilon) + \mu F_{\mu}(\varepsilon)$$ 
We proved that the sequence of $\mu$-foveal maps $\displaystyle
\phi^{x}_{\varepsilon}$  is a  
zoom sequence with zoom modulus $\displaystyle F(\cdot) + \mu F_{\mu}(\cdot)$. 

In order to prove that the sequence is scale stable, we have to compute 
$\displaystyle \phi^{x}_{\varepsilon, \mu}$, graphically shown in the next
figure.

\vspace{.5cm}

\centerline{\includegraphics[angle=270, width=0.7\textwidth]{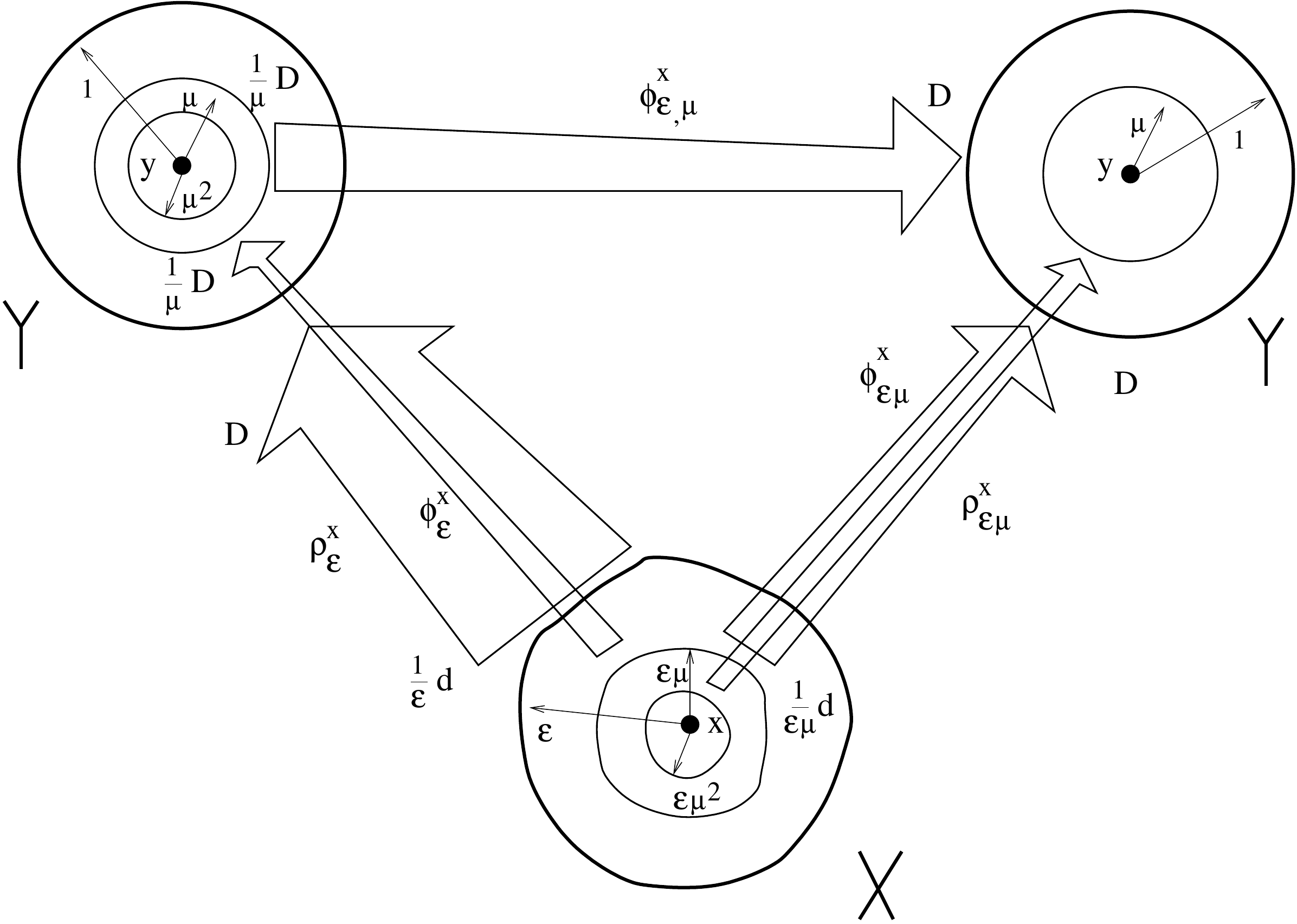}}

\vspace{.5cm}

We see that $\displaystyle (u',u") \in \phi^{x}_{\varepsilon, \mu}$ implies 
that $\displaystyle (u',u") \in \rho^{x}_{\varepsilon, \mu}$ or 
$\displaystyle (u',u") \in \rho^{x}_{\varepsilon \mu, \mu}$. From here we deduce
that the sequence of foveal maps is scale stable and that 
$$\varepsilon \mapsto \max \left\{ F_{\mu}(\varepsilon) , \mu
F_{\mu}(\varepsilon \mu) \right\}$$ 
is a scale stability modulus for the foveal sequence.

The improvement of the right hand side for the cascading of errors
inequality (\ref{eqcascade}), applied for the restricted
foveal map is then straightforward if we use (\ref{nd1}).  \hfill $\square$

\vspace{.5cm}

\section{Appendix II: Dilation structures}

From  definition \ref{deffoveal} we see that 
\begin{equation}
\bar{\rho}^{x}_{\mu} \circ \phi^{x}_{\varepsilon} \, = \, 
\rho^{x}_{\varepsilon \mu}
\label{pre1param}
\end{equation}

Remark that if the $\mu$-foveal map $\displaystyle \phi^{x}_{\varepsilon}$ 
coincides with the chart $\displaystyle \rho^{x}_{\varepsilon}$ for every 
$\varepsilon$ (that is, if 
the zoom sequence $\displaystyle \rho^{x}_{\varepsilon}$ is already so good 
that it cannot be improved by the construction of foveal maps), then relation 
(\ref{pre1param}) becomes 

\begin{equation}
\bar{\rho}^{x}_{\mu} \circ \phi^{x}_{\varepsilon} \, = \, 
\phi^{x}_{\varepsilon \mu}
\label{pre2param}
\end{equation}

By proposition \ref{pstable}, it follows that $\mu$-foveal map at scale 
$\varepsilon \mu$ is just a $1/\mu$ dilation of a part of the $\mu$-foveal map at scale 
$\varepsilon$. 

An idealization of these "perfect", stable zoom sequences which cannot be
improved by the $\mu$-foveal map construction for any $\mu \in (0,1)$, are
dilation structures. 

There are several further assumptions, which clearly amount to yet other 
idealizations. These are the following:
\begin{enumerate}
\item[-] the "map is the territory assumption", namely $\displaystyle 
Y = U(x)$, the "map space" 
is included in $X$, the "territory", and $y = x$. 
\item[-] "functions instead relations", that is the perfect stable zoom
sequences $\displaystyle \rho^{x}_{\varepsilon} = \phi^{x}_{\varepsilon}$ 
are graphs of functions, called dilations. That means: 
$$\rho^{x}_{\varepsilon} \subset \, \left\{ (\delta^{x}_{\varepsilon} u' , u') 
\mbox{ : } u' \in Y = V_{\varepsilon}(x) \right\}$$
\item[-] "hidden uniformity", that is: in order to pass to the limit in various 
situations, we could choose the  zoom modulus  and stability
modulus to not depend on $x \in X$. This innocuous assumption is the least
obvious, but necessary one. 
\end{enumerate}

With these idealizations in force, remember that we want our dilations to form a stable zoom sequence and 
we want also the subtler viewpoint stability, which consists in 
being  able to change the point of view in a coherent way, 
as the scale goes to zero. These are the axioms of a dilation structure. 

We shall use here a slightly particular version of dilation structures. 
For the general definition of a dilation structure see \cite{buligadil1}. 
More about this, as well as about length dilation structures, see  
\cite{buligadil3}.

\begin{definition}
Let $(X,d)$ be a complete metric space such that for any $x  \in X$ the 
closed ball $\bar{B}(x,3)$ is compact. A dilation structure $(X,d, \delta)$ 
over $(X,d)$ is the assignment to any $x \in X$  and $\varepsilon \in (0,+\infty)$ 
of a  homeomorphism, defined as: if 
$\displaystyle   \varepsilon \in (0, 1]$ then  $\displaystyle 
 \delta^{x}_{\varepsilon} : U(x)
\rightarrow V_{\varepsilon}(x)$, else 
$\displaystyle  \delta^{x}_{\varepsilon} : 
W_{\varepsilon}(x) \rightarrow U(x)$,  with the following properties.  
\begin{enumerate}
\item[{\bf A0.}]  For any $x \in X$ the sets $ \displaystyle U(x), V_{\varepsilon}(x), 
W_{\varepsilon}(x)$ are open neighbourhoods of $x$.  There are   $1<A<B$ such that for any $x \in X$  and any 
$\varepsilon \in (0,1)$ we have: 
$$\displaystyle  B_{d}(x, \varepsilon) \subset \delta^{x}_{\varepsilon}  B_{d}(x, A) 
\subset V_{\varepsilon}(x) \subset $$ 
$$\subset W_{\varepsilon^{-1}}(x) \subset \delta_{\varepsilon}^{x}  B_{d}(x, B)$$
Moreover for any compact set $K \subset X$ there are $R=R(K) > 0$ and 
$\displaystyle \varepsilon_{0}= \varepsilon(K) \in (0,1)$  such that  
for all $\displaystyle u,v \in \bar{B}_{d}(x,R)$ and all 
$\displaystyle \varepsilon  \in (0,\varepsilon_{0})$,  we have 
$\displaystyle 
\delta_{\varepsilon}^{x} v \in W_{\varepsilon^{-1}}(
\delta^{x}_{\varepsilon}u)$. 

\item[{\bf A1.}]  For any $x \in X$ 
$\displaystyle  \delta^{x}_{\varepsilon} x = x $ and $\displaystyle \delta^{x}_{1} = id$. 
Consider  the closure $\displaystyle Cl(dom \, \delta)$ of the set 
$$ dom \, \delta = \left\{ (\varepsilon, x, y) \in (0,+\infty) \times X 
\times X \mbox{ : } \right.$$ 
$$\left. \mbox{ if } \varepsilon \leq 1 \mbox{ then } y 
\in U(x) \,
\, , \mbox{  else } y \in W_{\varepsilon}(x) \right\} $$ 
seen in  $[0,+\infty) \times X \times X$ endowed with
 the product topology. The function $\displaystyle \delta : dom \, \delta 
\rightarrow  X$,  $\displaystyle \delta (\varepsilon,  x, y)  = 
\delta^{x}_{\varepsilon} y$ is continuous, admits a continuous extension 
over $\displaystyle Cl(dom \, \delta)$ and we have 
$\displaystyle \lim_{\varepsilon\rightarrow 0} \delta_{\varepsilon}^{x} y \, =
\, x$. 

\item[{\bf A2.}] For any  $x, \in X$, $\displaystyle \varepsilon, \mu \in (0,+\infty)$
 and $\displaystyle u \in U(x)$, whenever one of the sides are well defined
   we have the equality 
$\displaystyle  \delta_{\varepsilon}^{x} \delta_{\mu}^{x} u  =
\delta_{\varepsilon \mu}^{x} u$.

\item[{\bf A3.}]  For any $x$ there is a distance  function $\displaystyle (u,v) \mapsto d^{x}(u,v)$, defined for any $u,v$ in the closed ball (in distance d) $\displaystyle 
\bar{B}(x,A)$, such that uniformly with respect to $x$ in compact set we have
the limit: 
$$\lim_{\varepsilon \rightarrow 0} \quad \sup  \left\{  \mid 
\frac{1}{\varepsilon} d(\delta^{x}_{\varepsilon} u, \delta^{x}_{\varepsilon} v) \ - \ d^{x}(u,v) \mid \mbox{ :  } u,v \in \bar{B}_{d}(x,A)\right\} \ =  \ 0$$

\item[{\bf A4.}] Let us define 
$\displaystyle \Delta^{x}_{\varepsilon}(u,v) =
\delta_{\varepsilon^{-1}}^{\delta^{x}_{\varepsilon} u} \delta^{x}_{\varepsilon} v$. 
Then we have the limit, uniformly with respect to $x, u, v$ in compact set,  
$$\lim_{\varepsilon \rightarrow 0}  \Delta^{x}_{\varepsilon}(u,v) =  \Delta^{x}(u, v)  $$
 
\end{enumerate}
\label{defweakstrong}
\end{definition}

\subsection{Dilations as morphisms: towards the chora}

It is algebraically straightforward to transport a dilation structure: given 
$(X,d,  \delta)$ a dilation structure and $f: X \rightarrow Z$ a uniformly
continuous homeomorphism from $X$ (as a topological space) to another topological
space $Z$ (actually more than a topological space, it should be a space endowed
with an uniformity), we can define the transport of $(X,d,\delta)$ by $f$ as 
the dilation structure $(Z,f*d, f*\delta)$. The distance $f*d$ is defined as 
$$\left(f*d\right) (u,v) = d(f(u),f(v))$$
which is a true distance, because we supposed $f$ to be a homeomorphism. 
For any $u,v \in X$ and $\varepsilon > 0$, we define the new dilation based 
at $f(u) \in Z$, of coefficient $\varepsilon$, applied to $f(v) \in Z$ as
$$\left(f*\delta\right)^{f(u)}_{\varepsilon} f(v) \, = \, f \left(
\delta^{u}_{\varepsilon} v \right)$$
It is easy to check that this is indeed a dilation structure. 

In particular we may consider to transport a dilation structure by one of its 
dilations. Visually, this corresponds to transporting the atlas representing 
a dilation structure on $X$ to a neighbourhood of one of its points. It is like 
a scale reduction of the whole territory $(X,d)$ to a smaller set. 

Inversely, we may transport the (restriction of the) dilation structure 
$(X,d,\delta)$ from $\displaystyle V_{\varepsilon}(x)$ to $U(x)$, by using 
$\displaystyle \delta^{x}_{\varepsilon^{-1}}$ as the transport function $f$. 
This is like a magnification of the "infinitesimal neighbourhood" 
$\displaystyle V_{\varepsilon}(x)$. (This neighbourhood is infinitesimal in the
sense that we may consider $\varepsilon$ as a variable, going to $0$ when
needed. Thus, instead of one neighbourhood $\displaystyle V_{\varepsilon}(x)$, 
there is a sequence of them, smaller and smaller). 

This is useful, because it allows us to make "infinitesimal statements", i.e. 
statements concerning this sequence of magnifications, as $\varepsilon
\rightarrow 0$. 

Let us compute then the magnified dilation structure. We should also rescale 
the distance on $\displaystyle V_{\varepsilon}(x)$ by a factor $1/\varepsilon$. 
Let us compute this magnified dilation structure: 
\begin{enumerate}
\item[-] the space is $U(x)$
\item[-] for any $u,v \in U(x)$ the (transported) distance between them is 
$$d^{x}_{\varepsilon}(u,v) \, = \, \frac{1}{\varepsilon}
d(\delta^{x}_{\varepsilon} u , \delta^{x}_{\varepsilon} v ) $$
\item[-] for any $u,v \in U(x)$ and scale parameter $\mu \in (0,1)$ 
(we could take $\mu > 0$ but then we have to be careful with the domains and
codomains of these new dilations), the transported dilation based at $u$, of 
coefficient $\mu$, applied to $v$, is 
\begin{equation}
 \delta^{x}_{\varepsilon^{-1}} \, \delta_{\varepsilon}^{\delta^{x}_{\varepsilon} u} \, 
\delta^{x}_{\varepsilon} v
\label{firstchora}
\end{equation}
\end{enumerate}

It is visible that working with such combinations of dilations becomes quickly
difficult. This is one of the reasons of looking for more graphical notations. 

\paragraph{Important remark.}  This sequence of magnified dilation structures, around $x$, at scale 
$\varepsilon$, with the parameter $\varepsilon$ seen as a variable 
converging to $0$, is related to the "chora", or place $x$ at scale 
varepsilon.

In particular, the transported dilation (\ref{firstchora}) 
appears in relation with the "elementary chora".

\begin{definition}
Let  $(X,d, \delta)$ be a dilation structure. A property 
$$\displaystyle \mathcal{P}(x_{1},x_{2},x_{3}, ...)$$ 
is true  for $\displaystyle x_{1}, x_{2}, x_{3}, ... \in X$  sufficiently 
close if for any compact, non empty set $K \subset X$, there
is a positive constant $C(K)> 0$ such that $\displaystyle \mathcal{P}(x_{1},x_{2},
x_{3}, ...)$ is true for any $\displaystyle x_{1},x_{2},
x_{3}, ... \in K$ with $\displaystyle d(x_{i}, x_{j}) \leq C(K)$.
\end{definition}

For a dilation structure the metric tangent spaces  have the algebraic
structure of a normed group with dilations.

 We shall work further with local groups, which are spaces endowed with 
a locally defined  operation which satisfies  the conditions of a uniform group. 
 See section 3.3 \cite{buligadil1} for details about the 
definition of local groups.

\subsection{Normed conical groups}

This name has been introduced in section 8.2 \cite{buligadil1}, but these
objects appear more or less in the same form under the name "contractible
group" or "homogeneous group". Essentially these are groups endowed with a 
family of "dilations". They were also studied  in section 4 
\cite{buligadil2}. 

In the following general definition appear a topological 
commutative group $\Gamma$ endowed with a continuous morphism $\nu: \Gamma \rightarrow (0, +\infty)$ from $\Gamma$ 
to the group $(0, +\infty)$ with multiplication.  
The morphism $\nu$ induces an invariant topological filter on $\Gamma$ (other names for such an invariant filter are "absolute" or "end").  The convergence 
of a variable $\varepsilon \in \Gamma$ to this filter is denoted by $\varepsilon \rightarrow 0$ and 
it means simply $\nu(\varepsilon) \rightarrow 0$ in $\mathbb{R}$. 

Particular, interesting examples of pairs $(\Gamma, \nu)$ are: $(0, +\infty)$ with identity, which is 
the case interesting for this paper, $\displaystyle \mathbb{C}^{*}$ with the modulus of complex numbers, 
or $\mathbb{N}$ (with addition) with the exponential, which is relevant for the case of normed 
contractible groups, section 4.3 \cite{buligadil2}. 

\begin{definition}
A  normed group with dilations $(G, \delta, \| \cdot \|)$ is a local  group 
$G$  with  a local action of $\Gamma$ (denoted by $\delta$), on $G$ such that
\begin{enumerate}
\item[H0.] the limit  $\displaystyle \lim_{\varepsilon \rightarrow 0}
\delta_{\varepsilon} x  =  e$ exists and is uniform with respect to $x$ in a compact neighbourhood of the identity $e$.
\item[H1.] the limit
$\displaystyle \beta(x,y)  =  \lim_{\varepsilon \rightarrow 0} \delta_{\varepsilon}^{-1}
\left((\delta_{\varepsilon}x) (\delta_{\varepsilon}y ) \right)$
is well defined in a compact neighbourhood of $e$ and the limit is uniform with
respect to $x$, $y$.
\item[H2.] the following relation holds: 
$\displaystyle \lim_{\varepsilon \rightarrow 0} \delta_{\varepsilon}^{-1}
\left( ( \delta_{\varepsilon}x)^{-1}\right)  =  x^{-1}$, 
where the limit from the left hand side exists in a neighbourhood $U \subset G$  of 
$e$ and is uniform with respect to $x \in U$.
\end{enumerate}

Moreover the group is endowed with a continuous norm
function $\displaystyle \|\cdot \| : G \rightarrow \mathbb{R}$ which satisfies
(locally, in a neighbourhood  of the neutral element $e$) the properties:
 \begin{enumerate}
 \item[(a)] for any $x$ we have $\| x\| \geq 0$; if $\| x\| = 0$ then $x=e$,
 \item[(b)] for any $x,y$ we have $\|xy\| \leq \|x\| + \|y\|$,
 \item[(c)] for any $x$ we have $\displaystyle \| x^{-1}\| = \|x\|$,
 \item[(d)] the limit
$\displaystyle \lim_{\varepsilon \rightarrow 0} \frac{1}{\nu(\varepsilon)} \| \delta_{\varepsilon} x \| = \| x\|^{N}$
 exists, is uniform with respect to $x$ in compact set,
 \item[(e)] if $\displaystyle \| x\|^{N} = 0$ then $x=e$.
  \end{enumerate}
  \label{dnco}
  \end{definition}

\begin{theorem}
(Thm. 15 \cite{buligadil1}) Let $(G, \delta, \| \cdot \|)$ be  a locally compact  normed local group with dilations. 
Then $(G, d, \delta)$ is a  dilation structure, where the dilations 
$\delta$ and the distance $d$ are defined by: 
$\displaystyle  \delta^{x}_{\varepsilon} u = x \delta_{\varepsilon} ( x^{-1}u)  \quad , \quad
d(x,y) = \| x^{-1}y\|$. 

Moreover $(G, d, \delta)$ is linear, in the sense of definition \ref{dline}. 
\label{tgrd}
\end{theorem}

\begin{definition}
A normed conical group $N$ is a normed  group with dilations  such that
for any $\varepsilon \in \Gamma$  the dilation
 $\delta_{\varepsilon}$ is a group morphism  and such that for any $\varepsilon >0$
  $\displaystyle  \| \delta_{\varepsilon} x \| = \nu(\varepsilon) \| x \|$.
\end{definition}

A normed conical group is the infinitesimal version of a normed group with
dilations (\cite{buligadil1} proposition 2).

\begin{proposition}
Let $(G, \delta, \| \cdot \|)$ be  a locally compact  normed local group with
dilations. Then  $\displaystyle (G,\beta, \delta, \| \cdot \|^{N})$ is a locally compact, 
local normed conical group, with operation $\beta$,  dilations $\delta$ and homogeneous norm $\displaystyle \| \cdot \|^{N}$.
\label{here3.4}
\end{proposition}

\subsection{Tangent bundle of a dilation structure}
\label{induced}

The most important metric and algebraic first order  properties of a dilation 
structure are  presented here as condensed statements, available in full length as theorems 7, 8,
10 in \cite{buligadil1}.

\begin{theorem}
Let $(X,d,\delta)$ be a  dilation structure. Then the metric space $(X,d)$ 
admits a metric tangent space at $x$, for any point $x\in X$. 
More precisely we have  the following limit: 
$$\lim_{\varepsilon \rightarrow 0} \ \frac{1}{\varepsilon} \sup \left\{  \mid d(u,v) - d^{x}(u,v) \mid \mbox{ : } d(x,u) \leq \varepsilon \ , \ d(x,v) \leq \varepsilon \right\} \ = \ 0 \ .$$
\label{thcone}
\end{theorem}

\begin{theorem}
  If $(X,d,\delta)$ is a  dilation structure 
 then for any $x \in X$ the triple $\displaystyle (U^{x}, \delta^{x}, d^{x})$ is
 a locally compact normed conical group, with operation 
 $\displaystyle \Sigma^{x}(\cdot, \cdot)$, neutral element $x$ and inverse 
 $\displaystyle inv^{x}(y) = \Delta^{x}(y, x)$.  
\label{tgene}
\end{theorem}

The conical group $\displaystyle (U(x), \Sigma^{x}, \delta^{x})$ can be seen as the tangent space 
of $(X,d, \delta)$ at $x$. We shall  denote it by  
$\displaystyle T_{x} (X, d, \delta) =  (U(x), \Sigma^{x}, \delta^{x})$, or by $\displaystyle T_{x} X$ if 
$(d,\delta)$ are clear from the context.

The following proposition is  corollary 
6.3 from \cite{buligadil2}, which gives a more precise
description of the conical group 
$\displaystyle (U(x), \Sigma^{x}, \delta^{x})$.  In the proof of that corollary
there is a gap pointed by S. Vodopyanov, namely that  
 Siebert' proposition 5.4 \cite{siebert}, which 
is true for conical groups (in our language), is used for local  
conical groups. Fortunately, this gap was filled  by the 
theorem 1.1 \cite{recent}, which states that a locally compact, 
locally connected, contractible (with  Siebert' wording) group is 
locally isomorphic to a contractive Lie group.

\begin{proposition}
Let $(X,d,\delta)$ be a  dilation structure. 
Then for any $x \in X$ the local group
 $\displaystyle (U(x), \Sigma^{x})$ is locally a simply connected Lie group
 whose Lie algebra admits a positive graduation (a homogeneous group), given by
 the eigenspaces of $\displaystyle \delta^{x}_{\varepsilon}$ for an 
 arbitrary $\varepsilon \in (0,1)$.
\label{cor63}
\end{proposition}

There is a bijection between linear (in the sense of definition 
\ref{dline})   dilation structures and normed conical groups. Any 
normed conical group induces a linear  dilation structure, by theorem 
\ref{tgrd}. Conversely, we have the following result (see  theorem 6.1 
\cite{buligairq} for a more general statement).

\begin{theorem}
Let $(G,d, \delta)$ be a linear  dilation structure. Then, with the 
notations from  theorem \ref{tgene}, for any 
$x \in G$, the dilation structure $(U(x), d, \delta)$ coincides with 
the dilation structure of the conical group $\displaystyle (U(x), \Sigma^{x},
\delta^{x})$. 
\label{pgroudlin}
\end{theorem}

\subsection{Differentiability with respect to dilation structures}

For any  dilation structure or there is an associated  notion  of 
differentiability (section 7.2 \cite{buligadil1}). 
For defining differentiability with respect to dilation structures we need 
first the  definition of a morphism of conical groups. 

\begin{definition}
 Let $(N,\delta)$ and $(M,\bar{\delta})$ be two  conical groups. A function $f:N\rightarrow M$ is a conical group morphism if $f$ is a group morphism and for any $\varepsilon>0$ and $u\in N$ we have 
 $\displaystyle f(\delta_{\varepsilon} u) = \bar{\delta}_{\varepsilon} f(u)$. 
\label{defmorph}
\end{definition}

The definition of the derivative, or differential,  with respect to dilations 
structures is a straightforward generalization of  the 
definition of the Pansu derivative  \cite{pansu}.

 \begin{definition}
 Let $(X, d, \delta)$ and $(Y, \overline{d}, \overline{\delta})$ be two 
  dilation structures   and $f:X \rightarrow Y$ be a continuous function. The function $f$ is differentiable in $x$ if there exists a 
 conical group morphism  $\displaystyle D \, f(x):T_{x}X\rightarrow T_{f(x)}Y$, defined on a neighbourhood of $x$ with values in  a neighbourhood  of $f(x)$ such that 
\begin{equation}
\lim_{\varepsilon \rightarrow 0} \sup \left\{  \frac{1}{\varepsilon} \overline{d} \left( f\left( \delta^{x}_{\varepsilon} u\right) ,  \overline{\delta}^{f(x)}_{\varepsilon} D \, f(x)  (u) \right) \mbox{ : } d(x,u) \leq \varepsilon \right\}Ê  = 0 , 
\label{edefdif}
\end{equation}
The morphism $\displaystyle D \, f(x) $ is called the derivative, or differential,  of $f$ at $x$.
\label{defdiffer}
\end{definition}

The definition also makes sense if the function $f$ is defined on a 
open subset of $(X,d)$.

%\newpage

\section{Supplementary material}

\subsection{Gregory Bateson on maps and difference}

 The following are excerpts from the Nineteenth Annual Korzybski Memorial
 Lecture "Form, Substance and Difference", delivered by Bateson 
 on January 9, 1970, under the auspices of the Institute of
 General Semantics, re-printed from the {\em General Semantics} Bulletin, no.
 37, 1970,  in Steps to an Ecology of Mind (1972).

"Let us go back to the original statement for which Korzybski is most 
famous -- the statement that {\it the map is not the territory}. This statement 
came out of a very wide range of philosophic thinking, going back to Greece, and
wriggling through the history of European thought over the last 2000 years. In
this history, there has been a sort of rough dichotomy and often deep
controversy. There has been a violent enmity and bloodshed. It all starts, I
suppose, with the Pythagoreans versus their predecessors, and the argument took
the shape of "Do you ask what it's made of -- earth, fire, water, etc?" Or do
you ask "What is its pattern?" Pythagoras stood for inquiry into pattern rather
than inquiry into substance. That controversy has gone through the ages, and the
Pythagorean half of it has, until recently, been on the whole submerged half."

 "What is in the territory that gets onto the map?" [...] What gets onto 
the map, in fact, is {\it difference},[...] A difference is a very peculiar and
obscure concept. It is certainly not a thing or an event. This piece of paper is
different from the wood of this lectern. There are many differences between them
-- of color, texture, shape, etc. But if we start to ask about the localization
of those differences, we get into trouble. Obviously the difference between the
paper and the wood is not in the paper; it is obviously not in the wood; it is
obviously not in the space between them [...] 

A difference, then, is an abstract matter. [...]

Difference travels from the wood and paper into my retina. It then gets picked
up and worked on by this fancy piece of computing machinery in my head.

I suggest to you now that the word "idea", in its most elementary sense, is
synonymous with "difference". Kant, in the {\it Critique of Jugment} -- if I
understand him correctly -- asserts that the most elementary aestethic act is
the selection of a fact. He argues that in a piece of chalk there are an
infinite number of potential facts. The {\it Ding an sich}, the piece of chalk,
can never enter into communication or mental process because of this infinitude.
The sensory receptors cannot accept it; they filter it out. What they do is to
select certain {\it facts} out of the piece of chalk, which then become, in
modern terminology, information. 

I suggest that Kant's statement can be modified to say that there is an
infinite number of {\it differences} around and within the piece of chalk. Of
this infinitude, we select a very limited number, which become information. In
fact, what we mean by information -- the elementary unit of information -- is a
{\it difference which makes a difference}."

" We say the map is different from the territory. But what is the territory? 
Operationally, somebody went out with a retina or a measuring stick and made 
representations which were then put upon paper. What is on the paper map is a 
representation of what was in the retinal representation of the man who made 
the map; and as you push the question back, what you find is an infinite 
regress, an infinite series of maps. The territory never gets in at all.
 The territory is {\it Ding an sich} and you can't do anything with it. 
 Always, the process of representation will filter it out so that the mental 
 world is only maps of maps, ad infinitum. 
  Or we may spell the matter out and say that at every step, as a
  difference is transformed and propagated along its pathway, the embodiment of
  the difference before the step is a "territory" of which the embodiment after
  the step is a "map". The map-territory relation obtains at every step.
   All "phenomena" are literary appearances.
  
  Or we can follow the chain forward. I receive various sorts of mappings which
  I call data or information. Upon receipt of these I act. But my actions, my
  muscular contractions, are transforms of differences in the input material.
  And I receive again data which are transforms of my actions."

\subsection{Plato about chora}

Here are some relevant quotations (source \cite{plato}): 

 "This new beginning of our discussion of the universe requires a fuller 
 division than the former; for then we made two classes, now a third must 
 be revealed. The two sufficed for the former discussion: one, 
 which we assumed, was a pattern intelligible and always the same; 
 and the second was only the imitation of the pattern, generated and visible. 
 There is also a third kind which we did not distinguish at the time, 
 conceiving that the two would be enough. But now the argument seems to 
 require that we should set forth in words another kind, which is 
 difficult of explanation and dimly seen. 

What nature are we to attribute to this new kind of being? We reply, that it 
is the receptacle, and in a manner the nurse, of all generation.

 I have spoken the truth; but I must express myself in clearer language, 
 and this will be an arduous task for many reasons, and in particular 
 because I must first raise questions concerning fire and the other elements, 
 and determine what each of them is;"

"...as the several elements never present themselves in the same form, how 
can any one have the assurance to assert positively that any of them, 
whatever it may be, is one thing rather than another? No one can. But 
much the safest plan is to speak of them as follows:-Anything which we see 
to be continually changing, as, for example, fire, we must not call "this" 
or "that," but rather say that it is "of such a nature"..."

"That in which the elements severally grow up, and appear, and decay, is 
alone to be called by the name "this" or "that"; but that which is of a 
certain nature, hot or white, or anything which admits of opposite 
equalities, and all things that are compounded of them, ought not to be so 
denominated."  

""Let me make another attempt to explain my meaning more clearly. Suppose a 
person to make all kinds of figures of gold and to be always transmuting one 
form into all the rest-somebody points to one of them and asks what it is. 
By far the safest and truest answer is, That is gold; and not to call the 
triangle or any other figures which are formed in the gold "these," as though 
they had existence, since they are in process of change while he is making 
the assertion; but if the questioner be willing to take the safe and 
indefinite expression, "such," we should be satisfied. And the same argument 
applies to the universal nature which receives all bodies-that must be always 
called the same; for, while receiving all things, she never departs at all 
from her own nature, and never in any way, or at any time, assumes a form 
like that of any of the things which enter into her; she is the natural 
recipient of all impressions, and is stirred and informed by them, and 
appears different from time to time by reason of them."

"For the present we have only to conceive of three natures: first, that which 
is in process of generation; secondly, that in which the generation takes 
place; and thirdly, that of which the thing generated is a resemblance. And 
we may liken the receiving principle to a mother, and the source or spring to 
a father, and the intermediate nature to a child; and may remark further, 
that if the model is to take every variety of form, then the matter in which 
the model is fashioned will not be duly prepared, unless it is formless, 
and free from the impress of any of these shapes which it is hereafter to 
receive from without."

"Wherefore, the mother and receptacle of all created and visible and in 
any way sensible things, is not to be termed earth, or air, or fire, or 
water, or any of their compounds or any of the elements from which these 
are derived, but is an invisible and formless being which receives all 
things and in some mysterious way partakes of the intelligible, and is 
most incomprehensible. In saying this we shall not be far wrong; as far, 
however, as we can attain to a knowledge of her from the previous 
considerations, we may truly say that fire is that part of her nature 
which from time to time is inflamed, and water that which is moistened, 
and that the mother substance becomes earth and air, in so far as she 
receives the impressions of them."

After that, Plato continues:

"...we must acknowledge that there is one kind of being which is always the 
same, uncreated and indestructible, never receiving anything into itself 
from without, nor itself going out to any other, but invisible and 
imperceptible by any sense, and of which the contemplation is granted to 
intelligence only. 

And there is another nature of the same name with it, and like to it, 
perceived by sense, created, always in motion, becoming in place and 
again vanishing out of place, which is apprehended by opinion and sense.

And there is a third nature, which is {\em space}, and is eternal, and 
admits not of destruction and provides a home for all created things, 
and is apprehended without the help of sense, by a kind of spurious reason, 
and is hardly real; which we beholding as in a dream, say of all existence 
that it must of necessity be in some place and occupy a space, but that what 
is neither in heaven nor in earth has no existence. Of these and other 
things of the same kind, relating to the true and waking reality of nature, 
we have only this dreamlike sense, and we are unable to cast off sleep and 
determine the truth about them."

The "dreamlike sense" comes from the fact that "true and exact reason"'  
leads to contradiction (52c): 

"For an image, since the reality, after which it is modeled, does not 
belong to it, and it exists ever as the fleeting shadow of some other, 
must be inferred to be in another... But true and exact reason, 
vindicating the nature of true being, maintains that while two things 
[i.e. the image and space] are different they cannot exist one of them 
in the other and so be one and also two at the same time."


\begin{thebibliography}{99}

\bibitem{wikimap} \url{http://en.wikipedia.org/wiki/Map-territory_relation#.22The_map_is_not_the_territory.22}


\bibitem{buligadil1} M. Buliga, Dilatation structures I. Fundamentals, {\it 
J. Gen. Lie Theory Appl.},  {\bf 1} (2007),  2, 65-95. 

\bibitem{buligadil2} M. Buliga, Infinitesimal affine geometry of metric spaces 
endowed with a dilatation structure ,  {\it Houston Journal 
of Math.} 36, 1 (2010), 91-136,  \url{http://arxiv.org/abs/0804.0135}

\bibitem{buligasr} M. Buliga, Dilatation structures in sub-riemannian geometry, 
 in: Contemporary Geometry and Topology and Related Topics. 
Cluj-Napoca, Cluj-Napoca, Cluj University Press (2008), 89-105

\bibitem{buligadil3} M. Buliga, A characterization of sub-riemannian spaces as length dilation 
structures constructed via coherent projections, {\it Commun. Math. Anal.} 
{\bf 11} (2011), No. 2, pp. 70-111,   \\ 
\url{http://arxiv.org/abs/0810.5042}

\bibitem{buligairq} M. Buliga, Emergent algebras as generalizations of
differentiable algebras, with applications, (2009), \url{http://arxiv.org/abs/0907.1520}

\bibitem{buligabraided} M. Buliga, Braided spaces with dilations and sub-riemannian symmetric spaces, (2010),  
\url{http://arxiv.org/abs/1005.5031}

\bibitem{buligafrontend} M. Buliga, What is a space? Computations in emergent algebras and the front end 
visual system (2010), \url{http://arxiv.org/abs/1009.5028}

\bibitem{buligamore} M. Buliga, More than discrete or continuous: a bird's view
(2010), \url{http://arxiv.org/abs/arXiv:1011.4485}

\bibitem{buligaintro} M, Buliga, Introduction to metric spaces with dilations
(2010), \url{http://arxiv.org/abs/1007.2362}

\bibitem{buligagr} M. Buliga, Deformations of normed groupoids and differential calculus. 
First part, (2009), \url{http://arxiv.org/abs/0911.1300}

\bibitem{buliga2} M. Buliga, Tangent bundles to sub-Riemannian groups, (2003),  \\ 
 \url{http://xxx.arxiv.org/abs/math.MG/0307342}


\bibitem{duzhin} S. Chmutov, S. Duzhin, J. Mostovoy, Introduction to Vassiliev
Knot Invariants, (2011), \\ \url{http://arxiv.org/abs/1103.5628}

\bibitem{connes} A. Connes, Sur la theorie non commutative de l'integration, in: Alg\'ebres d'Op\'erateurs, S\'eminaire sur les Alg\'ebres d'Op\'erateurs, Les Plans-sur-Bex, Suisse, 13-18 mars 1978, Lecture Notes in Mathematics 725, ed. by 
A. Dold and B. Eckmann, Springer-Verlag 1079, p. 19-143



\bibitem{recent} L. van den Dries, I. Goldbring, Locally compact contractive local groups, 
(2009), http://arxiv.org/abs/0909.4565


\bibitem{fly} N. Franceschini, J.M. Pichon, C. Blanes, From insect vision to 
robot vision, {\it Phil. Trans.: Biological Sciences},  {\bf 337},   1281 (1992), 
Natural and Artificial Low-Level Seeing Systems, 283-294


\bibitem{fennrourke} R. Fenn, C. Rourke, Racks and Links in codimension two, 
 {\it J. Knot Theory Ramifications},  {\bf 1}  (1992),  no. 4, 343--406




\bibitem{hil1} Hillier B, Penn A (2004) Rejoinder to Carlo Ratti in Environment and
 Planning B: Planning and Design 31 512-511 ISSN 0265 8135





\bibitem{joyce} D. Joyce, A classifying invariant of knots; the knot quandle, 
{\it J. Pure Appl. Alg.}, {\bf 23} (1982), 37-65 

\bibitem{kassel} C. Kassel, V. Turaev, Chord diagram invariants of tangles and
graphs, {\it Duke Math. J.}, {\bf 92} (1998), no. 3. 497-552


\bibitem{koen} J. Koenderink, The brain a geometry engine, {\it Psychol. Res.} {\bf 52} (1990), 122-127

\bibitem{koen2} J.. Koenderink, A. Kappers, A. van Doorn, Local Operations :The Embodiment of Geometry. Basic Research Series, (1992), 1-23

\bibitem{korzybski} A. Korzybski, A Non-Aristotelian System and its Necessity 
for Rigour in Mathematics and Physics, a paper presented before the American 
Mathematical Society at the New Orleans, Louisiana, meeting of the American 
Association for the Advancement of Science, December 28, 1931. 
Reprinted in Science and Sanity, 1933, p. 747-761. 




\bibitem{meredith} L.G. Meredith, D.F. Snyder, Knots as processes: a new kind of
invariant, (2010) \\  
\url{http://arxiv.org/abs/1009.2107} 

\bibitem{metrology} JCGM 200:2008 International vocabulary of metrology - Basic
and general concepts and associated terms (VIM), 
 \url{http://www.bipm.org/utils/common/documents/jcgm/JCGM_200_2008.pdf}




\bibitem{pansu} P. Pansu, M\'etriques de Carnot-Carath\'eodory et
quasi-isom\'etries des espaces sym\'etriques de rang un, Ann. of Math., (2) 
{\bf 129}, (1989), 1-60

\bibitem{plato} Plato, Timaeus 48e - 53c,  \url{http://www.ellopos.net/elpenor/physis/plato-timaeus/space.asp}



\bibitem{siebert} E. Siebert, Contractive automorphisms on locally compact 
groups, {\it Math. Z.}, 191, 73-90, (1986)




\bibitem{toffoli} E. Fredkin,  T. Toffoli, Conservative logic, {\it Int. J. of  Theoretical Physics}, {\bf 21} (1982), %no. 3-4, 219-253


\bibitem{hil2} Turner A, Hillier B, Penn A (2005) An algorithmic definition of the
 axial map Environment and Planning B 32-3, 425-444





\end{thebibliography}
\end{document}